\documentclass{article}
\usepackage[utf8]{inputenc}

\usepackage{amsfonts}
\usepackage{amsmath}
\usepackage{mathrsfs}
\usepackage{amssymb}
\usepackage{authblk}
\usepackage{graphics}
\usepackage{grffile}
\usepackage{epstopdf}
\usepackage{hyperref}
\usepackage{subcaption}
\usepackage{tikz}
\usepackage{comment}
\usepackage{soul}
\usepackage{color}

\usepackage[square,numbers]{natbib}
\usepackage{geometry}
\newcommand{\bm}[1]{\boldsymbol{#1}}
\newcommand{\wtilde}{\widetilde}
\newcommand{\what}{\widehat}
\newcommand{\dd}{\mathrm{d}}
\newcommand{\pd}{\partial}
\newcommand{\Tr}{\intercal}

\DeclareMathAlphabet{\mathsfit}{T1}{\sfdefault}{\mddefault}{\sldefault}
\SetMathAlphabet{\mathsfit}{bold}{T1}{\sfdefault}{\bfdefault}{\sldefault}

\providecommand{\keywords}[1]{\textbf{Keywords: } #1}

\title{A coupled discontinuous Galerkin-Finite Volume framework for solving gas dynamics over embedded geometries}
\author[1]{Vincenzo Gulizzi}
\author[1]{Ann S. Almgren}
\author[1]{John B. Bell}
\affil[1]{Center for Computational Sciences and Engineering (CCSE), Lawrence Berkeley National Laboratory MS 50A-3111, Berkeley, CA 94720, USA}

\begin{document}

\date{}
\maketitle

\begin{abstract}
We present a computational framework for solving the equations of inviscid gas dynamics using structured grids with embedded geometries.
The novelty of the proposed approach is the use of high-order discontinuous Galerkin (dG) schemes and a shock-capturing Finite Volume (FV) scheme coupled via an $hp$ adaptive mesh refinement ($hp$-AMR) strategy that offers high-order accurate resolution of the embedded geometries.
The $hp$-AMR strategy is based on a multi-level block-structured domain partition in which each level is represented by block-structured Cartesian grids and the embedded geometry is represented implicitly by a level set function.
The intersection of the embedded geometry with the grids produces the implicitly-defined mesh that consists of a collection of regular rectangular cells plus a relatively small number of irregular curved elements in the vicinity of the embedded boundaries.
High-order quadrature rules for implicitly-defined domains enable high-order accuracy resolution of the curved elements with a cell-merging strategy to address the small-cell problem.
The $hp$-AMR algorithm treats the system with a second-order finite volume scheme at the finest level to dynamically track the evolution of solution discontinuities while using dG schemes at coarser levels to provide high-order accuracy in smooth regions of the flow.
On the dG levels, the methodology supports different orders of basis functions on different levels.
The space-discretized governing equations are then advanced explicitly in time using high-order Runge-Kutta algorithms.
Numerical tests are presented for two-dimensional and three-dimensional problems involving an ideal gas.
The results are compared with both analytical solutions and experimental observations and demonstrate that the framework provides high-order accuracy for smooth flows and accurately captures solution discontinuities.\\

\end{abstract}


\section{Introduction}\label{sec:INTRO}
In spite of the enormous amount of work on solving the equations of gas dynamics, there is still considerable interest in developing higher-order methods in complex domains.
In general, computational methods for problems in irregular domains can be classified into those employing body-fitted meshes (either unstructured or structured with curvilinear coordinates) and those employing structured meshes with embedded boundaries (EB).
Body-fitted mesh approaches use irregular cells where the boundaries of the mesh elements coincide with the boundaries of the domain; although very flexible for modeling complex geometrical features, these methods often require significant effort to generate high-quality elements.
On the other hand, EB methods represent the domain as a distinguished interface in a uniform background grid.
The EB approach offers a number of benefits,
such as automatic grid generation, data storage and adaptive mesh refinement \cite{aftosmis2000parallel}, while the regularity of the mesh elements allows one to use methods designed for uniform grids in most of the domain.
However, additional considerations need to be addressed in the presence of an embedded geometry.

The are two main approaches to designing EB methods, diffuse-interface methods \cite{peskin2002immersed,mokbel2018phase,kemm2020simple}, where the embedded boundaries or interfaces are smeared over a finite-thickness region that is captured with a sufficiently fine mesh, and sharp-interface methods \cite{moes1999finite,qin2013discontinuous,alauzet2016nitsche,saye2017implicit,saye2017implicitII,milazzo2018extended}, where the embedded boundaries or interfaces are represented as a lower dimensional surface.
The latter is often the recommended approach if the flow details in proximity of the geometries are of interest.
However, the sharp-interface methods suffer from the so-called \emph{small-cell} problem, i.e.~the presence of cells cut by the embedded geometry that have unacceptably small volume fractions, which induces overly restrictive time steps and ill-conditioned discrete operators if not properly treated \cite{berger2017cut}.

Most EB discretizations are based on Finite Volume (FV) schemes \cite{leveque2002finite,versteeg2007introduction} that are relatively easy to implement and offer robust shock-tracking capabilities.
In embedded-boundary FV approaches, several strategies have been explored to ameliorate the small-cell problem.  Notable examples include cell-merging \cite{clarke1986euler,gaffney1987euler,hartmann2008adaptive,ji2010numerical,hartmann2011strictly,cecere2014immersed,muralidharan2016high}, flux-redistribution  \cite{pember1995adaptive,almgren1997cartesian,colella2006cartesian,graves2013cartesian},
the $h$-box method \cite{helzel2005high,berger2012simplified}, the dimensionally-split Cartesian cut cell method 
\cite{gokhale2018dimensionally} and state redistribution \cite{berger2021state}.
However, embedded-boundary FV methods are typically limited to second-order accuracy in space and, if high-order accuracy is desired, see e.g.~\cite{muralidharan2016high}, they require large stencil to perform the solution reconstruction.

Recently, the discontinuous Galerkin (dG) methods \cite{cockburn1990runge,cockburn1998runge,arnold2002unified,cockburn2018discontinuous} have gained popularity thanks to the flexible variational framework on which they are based.
The dG methods represent the unknown fields in a suitably chosen space of basis functions, whose coefficients become the unknowns of the discrete problem.
This approach has several highly desirable features for parallel computing such as high-order accuracy via compact stencils, a natural way to handle $hp$ refinement, more general (e.g.~polytopic) element shapes \cite{cangiani2017hp,antonietti2019v} and, unlike high-order schemes based on continuous approximations, block-structured mass matrices.
As a result, the dG method, typically in combination with a cell-merging technique, has been successfully employed in embedded-boundary numerical schemes for the solution of several scientific and engineering problems including elliptic partial differential equations \cite{bastian2009unfitted,johansson2013high,brandstetter2015high,gurkan2019stabilized}, reaction-diffusion equations \cite{lew2008discontinuous}, two-phase flow \cite{heimann2013unfitted,kummer2017extended}, two-phase flow with moving boundaries \cite{krause2017incompressible}, surface tension dynamics, free-surface flow and rigid body-fluid interaction in the incompressible regime \cite{saye2017implicit,saye2017implicitII,saye2020fast}, and solid mechanics of thin structures with cutouts \cite{gulizzi2020implicit,gulizzi2020high}.

In the context of compressible flow, embedded-boundary dG methods are less thoroughly investigated.
For two-dimensional (2D) problems, Qin and Krivodonova \cite{qin2013discontinuous} presented a dG methods for inviscid gas dynamics, M{\"u}ller et al.\cite{muller2017high} developed a dG method that includes viscous terms and, recently, Geisenhofer et al.\cite{geisenhofer2019discontinuous} presented a dG method for inviscid compressible flow where the formation of solution discontinuities was captured by adding a suitably-chosen artificial viscosity to the governing equations.
The EB approaches mentioned above use Cartesian background grids.
Fidkowski and Darmofal \cite{fidkowski2007triangular} presented a 2D dG method using triangular cut cells for compressible flow, whereas, to the best of our knowledge, the only EB method involving both two- and three-dimensional computations of compressible fluids was recently proposed by Xiao et al.\cite{xiao2019immersed}, who used triangular and tetrahedral meshes as their background grid.

As seen in the brief literature review provided above, most embedded-boundary dG schemes are for elliptic and incompressible fluid flow equations, which do not generally involve the formation and evolution of solution discontinuities, unless sharp corners or cracks are present in the geometry.
Conversely, the development of shock waves and contact discontinuities are key features of inviscid compressible flow.
When discontinuities are present, higher-order methods, including dG methods, produce nonphysical oscillations, referred to as the Gibbs phenomenon.
Controlling oscillations in higher-order methods requires the introduction of a \emph{limiting} strategy.
Although various shock-capturing approaches for dG methods exist in the literature, such as artificial-viscosity methods \cite{persson2006sub,chaudhuri2017explicit}, moment limiters
\cite{krivodonova2007limiters} or shock-fitting deforming-mesh approaches \cite{zahr2020implicit}, among others, shock-capturing methods for embedded-boundary dG schemes is still not a mature area.

The aim of this work is to present a computational framework for embedded-boundary approaches that provides high-accuracy resolution of the geometry, high-order accuracy in regions of smooth flows and shock-capturing capabilities in the presence of solution discontinuities for two- and three-dimensional problems.
This is achieved by coupling high-order dG schemes with a shock-capturing FV scheme through an $hp$ adaptive mesh refinement ($hp$-AMR) strategy featuring high-accuracy resolution of embedded geometries.

The paper is organized as follows:
Sec.(\ref{sec:IBVP}) introduces the considered initial-boundary value problem for inviscid gas dynamics;
Sec.(\ref{sec:FRAMEWORK}) describes the proposed framework that includes the construction of the implicitly-defined mesh, the formulation of the corresponding high-order dG schemes and the shock-capturing FV scheme, and the $hp$-AMR strategy;
Sec.(\ref{sec:RESULTS}) presents a set of two- and three-dimensional numerical tests.
Finally, Sec.(\ref{sec:CONCLUSIONS}) concludes the work and discusses potential further developments.

\section{Governing equations}\label{sec:IBVP}
Inviscid gas dynamics is described by a system of first-order hyperbolic conservation laws given by
\begin{equation}\label{eq:governing equations - vector form}
\pd_t\bm{U}+\pd_k\bm{F}_k = \bm{0},
\end{equation}
where $\pd_t$ denotes the partial derivative with respect to the time $t$ and $\pd_k$ denotes the partial derivative with respect to the $k$-th coordinate $x_k$ of the vector $\bm{x} = \{x_k\}$ in $d$-dimensional space.
In Eq.(\ref{eq:governing equations - vector form}), $\bm{U}$ is the $(d+2)$-dimensional vector of conserved variables and $\bm{F}_k$ is the $(d+2)$-dimensional vector of the flux in the $k$ direction  defined as 
\begin{equation}\label{eq:governing equations - vectors}
\bm{U}\equiv
\left\{\begin{array}{c}
\rho\\
\rho\bm{v}\\
\rho E
\end{array}
\right\}
\quad\mathrm{and}\quad
\bm{F}_k\equiv
\left\{\begin{array}{c}
\rho v_k\\
\rho v_k\bm{v}+p\bm{\delta}_k\\
(\rho E+p)v_k
\end{array}
\right\},
\end{equation}
where $\bm{\delta}_k$ is a $d$-dimensional vector containing 1 in its $k$-th component and 0 in all the others, $\rho$ is the density, $\bm{v} = \{v_k\}$ is the velocity vector, $E$ is the total energy and $p$ is the pressure.
The governing equations are closed by an equation of state.
Here, we assume an ideal gas with ratio of specific heats $\gamma$, such that 
\begin{equation}\label{eq:total energy}
p = (\gamma-1)\rho\left(E-\frac{1}{2}v_kv_k\right).
\end{equation}
The governing equations given in Eq.(\ref{eq:governing equations - vector form}) are assumed to be valid for $\{t,\bm{x}\} \in \mathscr{I}_T\times\mathscr{D}$, where $\mathscr{I}_T \equiv [0,T]$ is the time interval, $T$ is the final time, and $\mathscr{D} \subset \mathbb{R}^d$ is the spatial domain.
The governing equations are supplemented by initial conditions and boundary conditions, which must be enforced on the boundary $\pd\mathscr{D}$ of the domain $\mathscr{D}$.

In Eqs.(\ref{eq:governing equations - vector form}-\ref{eq:total energy}) and in the remainder of the paper, the indices $k$ and $l$ are used as subscripts to denote the component of a vector or an array; they take values in $\{1,\dots,d\}$ and imply summation when repeated unless it is explicitly stated otherwise.
Finally, we note that the symbol $p$ will also be used to denote the order of the basis functions of dG schemes but its meaning will be clear from the context.

\section{Coupled discontinuous Galerkin-Finite Volume framework}\label{sec:FRAMEWORK}

\subsection{Implicitly-defined mesh}\label{ssec:IMPLICIT MESH}
\begin{figure}
\centering
\includegraphics[width = \textwidth]{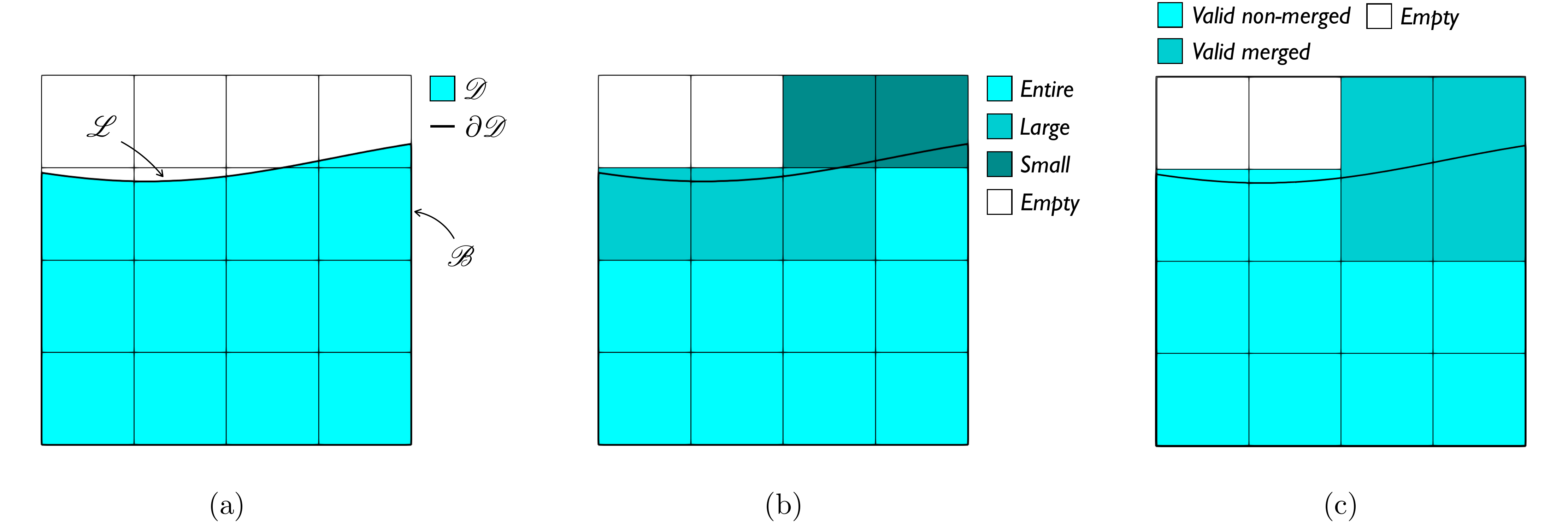}
\caption
{
    (a) Domain $\mathscr{D}$ and its boundary $\pd\mathscr{D} \equiv \mathscr{B} \cup \mathscr{L}$ implicitly defined in a $4\times 4$ background grid.
    (b) Classification of the background cells of figure (a) on the basis of their volume fraction.
    (c) Classification of the background cells of figure (a) after cell-merging.
}
\label{fig:IBVP - GEOMETRY AND CLASSIFICATION}
\end{figure}

\begin{figure}
\centering
\includegraphics[width = \textwidth]{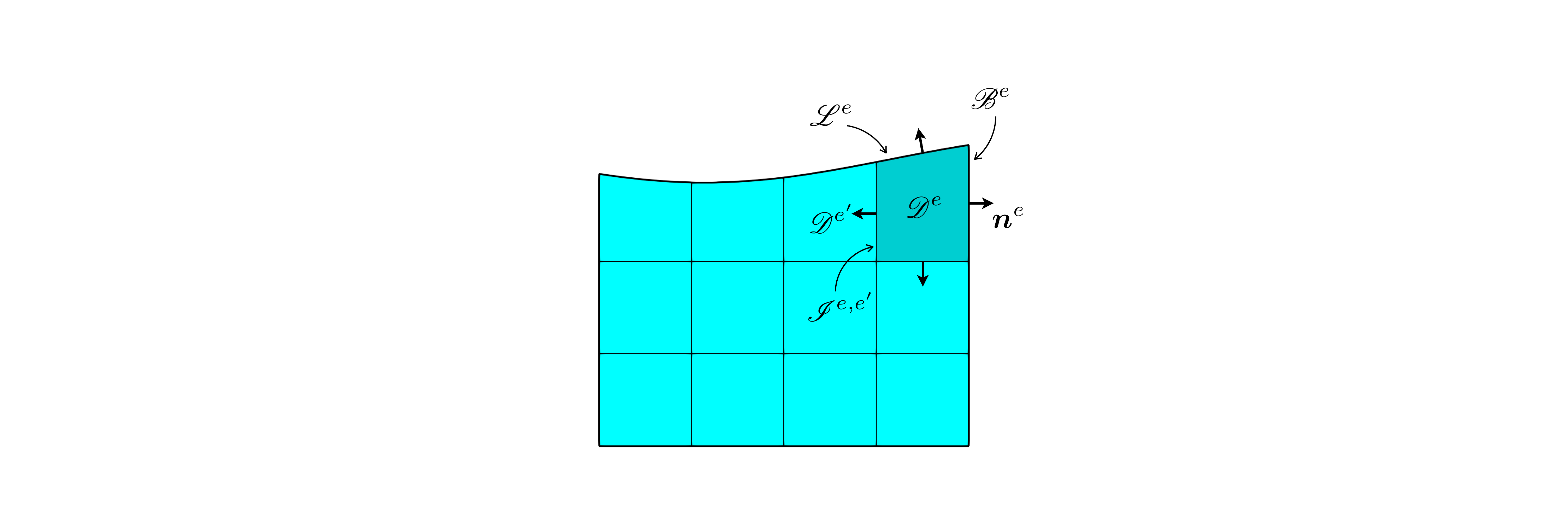}
\caption
{
    Implicitly-defined mesh of the domain shown in Fig.(\ref{fig:IBVP - GEOMETRY AND CLASSIFICATION}); the figure also highlights a mesh element $\mathscr{D}^{e}$ (in darker color) obtained from the merging between a small cell and an entire cell.
    See Figs.(\ref{fig:IBVP - GEOMETRY AND CLASSIFICATION}b) and (\ref{fig:IBVP - GEOMETRY AND CLASSIFICATION}c).
}
\label{fig:IBVP - MESH}
\end{figure}

\begin{figure}
\centering
\includegraphics[width = \textwidth]{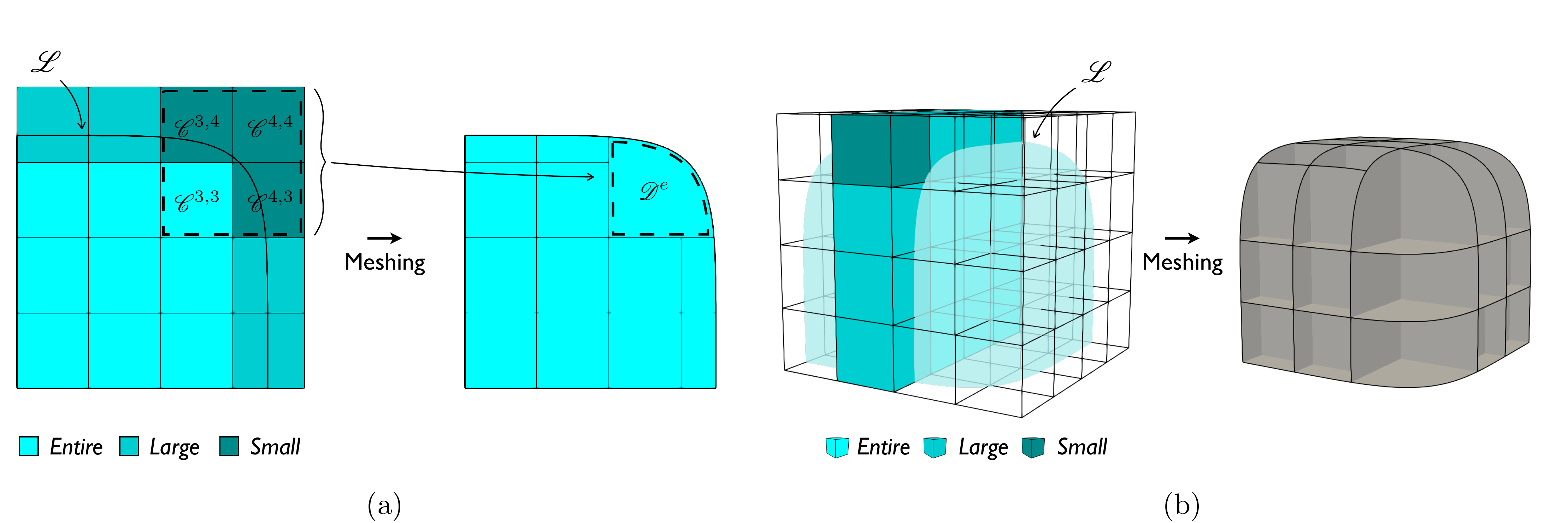}
\caption
{
    Examples of (a) two-dimensional and (b) three-dimensional implicitly-defined meshes where multiple small cells are merged with the same valid cell.
}
\label{fig:IBVP - MERGING DETAIL}
\end{figure}

Let $\mathscr{D}$ be the domain, which we assume to be enclosed in a background rectangle $\mathscr{R}\supseteq\mathscr{D}$.
We then let $\Phi:\mathscr{R}\rightarrow\mathbb{R}$ be a level set function such that $\mathscr{D}$ is implicitly defined as the region where $\Phi$ is negative, i.e.~$\mathscr{D} \equiv \{\bm{x} \in \mathscr{R} : \Phi(\bm{x}) < 0\}$, and its boundary $\pd\mathscr{D}$ is defined as $\pd\mathscr{D} \equiv \mathscr{B} \cup \mathscr{L}$, where $\mathscr{B} \equiv \{\bm{x}\in\pd\mathscr{R} : \Phi(\bm{x}) < 0\}$ is the portion of the rectangle's outer boundary $\pd\mathscr{R}$ where $\Phi$ is negative, and $\mathscr{L} \equiv \{\bm{x}\in\mathscr{R} : \Phi(\bm{x}) = 0\}$ is the zero level of $\Phi$.
We then introduce a structured grid $\mathscr{G}_h \equiv \cup_{\bm{i}}\mathscr{C}^{\bm{i}}\subseteq\mathscr{R}$, where $h$ denotes a characteristic mesh size, $\mathscr{C}^{\bm{i}}$ is a generic rectangular cell and $\bm{i} \equiv \{i_k\}$ is the $d$-tuple identifying the location of the cell within the grid.
The cell $\mathscr{C}^{\bm{i}}$ is defined as $[x_1^{\bm{i}},x_1^{\bm{i}}+h_1]\times\dots\times [x_d^{\bm{i}},x_d^{\bm{i}}+h_d]$, where $x_k^{\bm{i}}$ is the cell's lower corner in the $k$-th direction and $h_k$ is the grid's mesh size in that direction.

When intersected with the domain $\mathscr{D}$, cells are characterized by their volume fraction $\nu$:
\emph{entire} cells fall entirely inside $\mathscr{D}$ and have volume fraction $\nu = 1$;
\emph{empty} cells fall entirely outside $\mathscr{D}$ and have volume fraction $\nu = 0$;
\emph{partial} cells are cut by $\mathscr{L}$ and have volume fraction $0 < \nu < 1$.
We then introduce a volume fraction threshold $\overline{\nu}$ to label each partial cell as either \emph{large} if $\overline{\nu} < \nu < 1$ or \emph{small} if $0 < \nu \le \overline{\nu}$.
As an example, Fig.(\ref{fig:IBVP - GEOMETRY AND CLASSIFICATION}a) shows a two-dimensional domain embedded in a $4\times 4$ grid and Fig.(\ref{fig:IBVP - GEOMETRY AND CLASSIFICATION}b) shows the corresponding cell classification.

The volume fraction threshold triggers the cell-merging strategy, whereby all cells labelled as small are merged with nearby cells.
Several possibilities exist for implementing the cell-merging strategy.
Here, a small cell is merged with one of its entire or large neighboring cells, henceforth collectively referred to as \emph{valid} cells.  The neighbor targeted for merging is selected according to its location with respect to the small cell and its volume fraction.
Specifically, neighboring cells are queried in the following order: in 3D, we first consider the set of valid cells sharing a face with the small cell, second we consider the set of valid cells sharing an edge with the small cell and finally we consider the set of valid cells sharing a corner with the small cell. 
Within each set, the neighbors are ordered according to their volume fraction.
Then, the target neighbor is the first cell in the first non-empty set.
It follows that, in the choice of the target neighbor, a large cell sharing a face with the small cell has higher priority than an entire cell sharing an edge with the small cell.
In 2D, the search starts from the cells that share an edge with a small cell.

Once the cell-merging procedure is applied, the cells of the grid $\mathscr{G}_h$ consist of \emph{valid non-merged} cells, i.e.~the entire and large cells that have not been targeted during the merging process, \emph{valid merged} cells, i.e.~the union of small cells and their merging large or entire neighbors, and empty cells.
Such a classification for the background cells of Fig.(\ref{fig:IBVP - GEOMETRY AND CLASSIFICATION}b) is shown in Fig.(\ref{fig:IBVP - GEOMETRY AND CLASSIFICATION}c).
Then, the implicitly-defined mesh $\mathcal{M}_h$ associated with $\mathscr{G}_h$ is defined by the intersection of the valid merged and non-merged cells of $\mathscr{G}_h$ and the domain $\mathscr{D}$, such that $\mathcal{M}_h \equiv \{\mathscr{D}^{e}\}$ and $\mathscr{D}^{e}$ is the $e$-th implicitly-defined mesh element.
Figure (\ref{fig:IBVP - MESH}) shows the implicitly-defined mesh associated with the grid of Fig.(\ref{fig:IBVP - GEOMETRY AND CLASSIFICATION}c) and highlights (in darker color) an implicitly-defined element $\mathscr{D}^e$ resulting from the merging of a small cell and an entire cell.
The figure also shows the outer boundaries $\mathscr{B}^e$ and $\mathscr{L}^e$ of the $e$-th element, which stem from the intersection of the grid with $\mathscr{B}$ and $\mathscr{L}$, respectively, and the internal boundaries $\mathscr{I}^{e,e'}$ that the $e$-th element shares with its generic neighboring element $e'$.
A few comments are worth pointing out:
\begin{itemize}
          \item In the implementation of the implicit mesh, we do not explicitly merge cells.  When a small cell and a valid cell have been merged, both the dG and the reconstructed FV representations of the solution on that merged cell are specified in terms of coefficients associated with the valid cell.
          Moreover, the geometry of the curved elements is never explicitly constructed or parameterized and is only inferred by the quadrature rules for domains and boundaries implicitly defined by the level set function.
    \item The cell-merging procedure needs to traverses the cells only once and is fully parallelizable.
          In fact, the volume fraction of the cells is computed prior to and does not change during the merging process; for example, after being merged, a small cell does not become a large cell eligible for merging with another small cell.
          Each small cell queries its neighboring cells and selects the most appropriate neighbor according to the procedure described above.
          Note that multiple small cells can target the same cell for merging; as a 2D example, Fig.(\ref{fig:IBVP - MERGING DETAIL}a) shows a 4$\times$4 grid where the small cells $\mathscr{C}^{4,3}$, $\mathscr{C}^{3,4}$ and $\mathscr{C}^{4,4}$ all merge with the cell $\mathscr{C}^{3,3}$ to define the implicitly-defined mesh element $\mathscr{D}^e$.
          A similar example in 3D is shown in Fig.(\ref{fig:IBVP - MERGING DETAIL}b).
    \item In general, and especially for high volume fraction thresholds, it is possible that the merging algorithm does not find a suitable neighbor, i.e., given a small cell, all the cells in its 3$\times$3 (in 2D) or its 3$\times$3$\times$3 (in 3D) neighborhood are small or empty; if this occurs, the simulation will terminate and require a smaller volume fraction threshold.
          However, in all the numerical tests considered here, a volume fraction threshold of 0.3 in 2D and a volume fraction threshold of 0.15 in 3D was found to be satisfactory.
          Other approaches could be implemented to deal with small cells that do not find suitable neighbors for merging.
          For example, one could consider merging multiple small cells together, thus increasing the effective mesh size of the elements, or one could lower the volume fraction threshold that triggers the merging and introduce an additional higher volume fraction thresholds that instead triggers a reduced polynomial order representation, thus improving the stability of the scheme.
          These approaches are, however, not considered here.
\end{itemize}

\begin{figure}
\centering
\includegraphics[width = \textwidth]{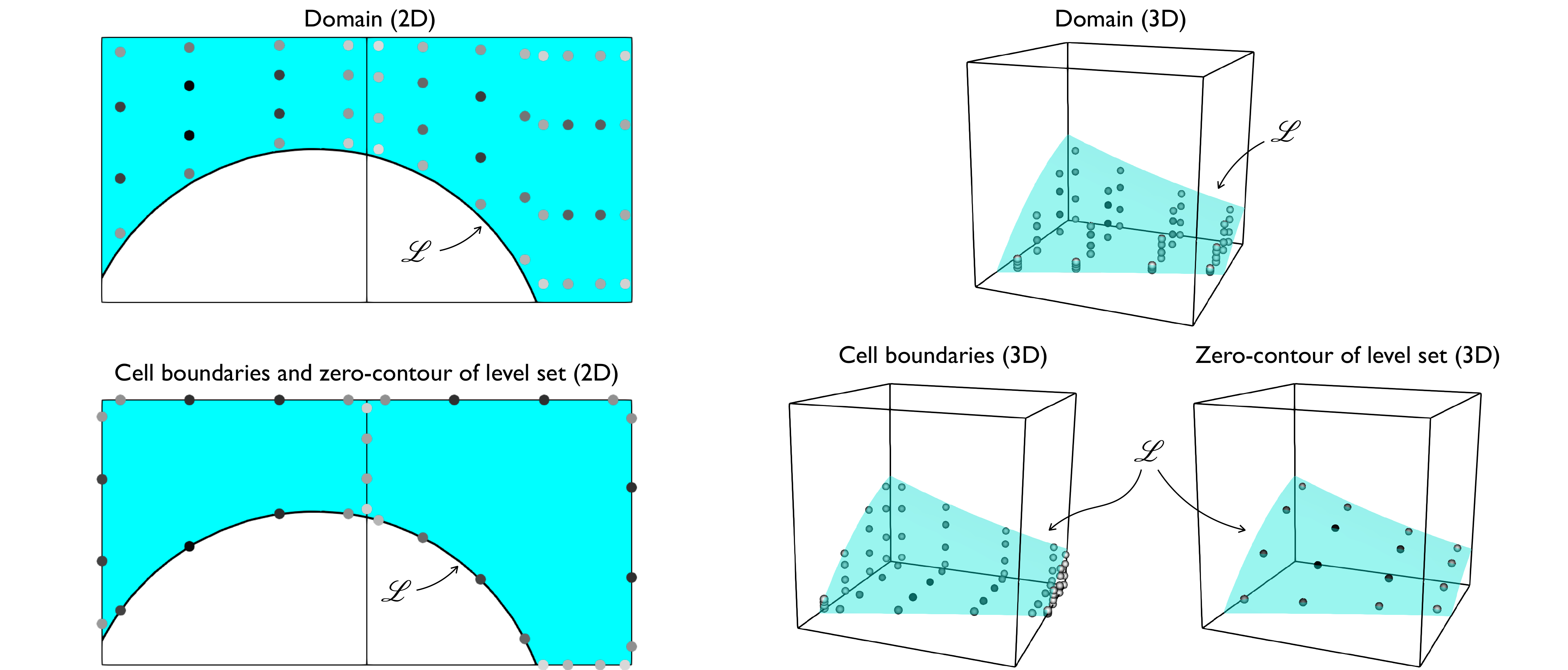}
\caption
{
    Examples of the location of the quadrature points obtained using the algorithm of Ref.\cite{saye2015high};
    in the images, the quadrature points are coloured according to their weight: light grays denote small weights while dark grays denotes large weights.
    Note that the quadrature weights are always strictly positive.
}
\label{fig:IBVP - QUADRATURE}
\end{figure}

\subsection{Discontinuous Galerkin schemes}\label{sec:DG}
Discontinuous Galerkin schemes are based on a weak form of the governing equations.
For each element in the mesh we define a local polynomial basis on that element and extend it to vanish
on all of the other elements. The collection of these functions for each element then forms a discontinuous basis.
Here, it is natural to define the local basis functions as tensor-product polynomials since the discretization is constructed using Cartesian grids and the majority of the mesh elements are rectangular cells.

Let $\mathscr{D}^{e}$ be the implicitly-defined element associated with the grid cell $\mathscr{C}^{\bm{i}}$ (and all the small cells that are merged with $\mathscr{C}^{\bm{i}}$) and let $\mathcal{P}_{hp}^e$ be the space of tensor-product polynomials of degree $p$ in the rectangular domain occupied by $\mathscr{C}^{\bm{i}}$.
Then, the space $\mathcal{V}_{hp}$ of discontinuous basis functions for the mesh $\mathcal{M}_h$ based on the grid $\mathscr{G}_h$ is
\begin{equation}\label{eq:space}
\mathcal{V}_{hp} \equiv \left\{v:\mathscr{G}_h\rightarrow\mathbb{R}~|~v|_{\mathscr{D}^e}\in\mathcal{P}_{hp}^e,~\forall \mathscr{D}^e\in\mathcal{M}_h\right\},
\end{equation}
while the related space $\mathcal{V}_{hp}^{N}$ of discontinuous polynomials vector fields is $\mathcal{V}_{hp}^N \equiv (\mathcal{V}_{hp})^N$.

Following Cockburn et al.\cite{cockburn1990runge,cockburn1998runge}, we integrate by parts in space to obtain the weak form of Eq.(\ref{eq:governing equations - vector form}) as
\begin{equation}\label{eq:governing equations - dG weak form}
\int_{\mathscr{D}^e}\bm{V}^\Tr\frac{\pd\bm{U}}{\pd t}\dd{V} = 
\int_{\mathscr{D}^e}\frac{\pd\bm{V}^\Tr}{\pd x_k}\bm{F}_k\dd{V}-
\int_{\mathscr{B}^{e}\cup\mathscr{L}^e}\bm{V}^\Tr\what{\bm{F}}_n\dd{S}-
\sum_{e'\in\mathcal{N}^{e}}\int_{\mathscr{I}^{e,e'}}\bm{V}^\Tr\what{\bm{F}}_n\dd{S} \;\;\; \forall \mathscr{D}^e\in\mathcal{M}_h \; \mathrm{and} \;\bm{V}\in\mathcal{V}_{hp}^{d+2}.
\end{equation}
In Eq.(\ref{eq:governing equations - dG weak form}), $\mathcal{N}^{e}$ denotes the set of mesh elements that are neighbors of $\mathscr{D}^e$, and $\what{\bm{F}}_n$ denotes the so-called \emph{numerical flux}, whose expression for hyperbolic equations is typically computed using an exact or an approximate Riemann solver \cite{colella1985efficient,toro2013riemann}.

We now approximate $\bm{U}$ on each element as a linear combination of basis functions with time dependent coefficients. This can be written compactly as
\begin{equation}
\bm{U}|_{\mathscr{D}^e}(\bm{x}) = \bm{B}^e(\bm{x})\bm{X}^e(t)
\label{eq:ansatz}
\end{equation}
where $\bm{B}^e$ is an $N_u \times N_u N_p$ matrix where each column corresponds to a basis function in $\mathcal{V}_{hp}$ and $\bm{X}^e(t)$ is a vector of coefficients of length $N_u N_p$. Here $N_u =d+2$ is the number of unknowns in the system and $N_p = (1+p)^d$ is the number of polynomials in $\mathcal{P}_{hp}^e$.
Substituting Eq. (\ref{eq:ansatz}) into Eq.(\ref{eq:governing equations - dG weak form})
and letting $\bm{V}$ range over the basis functions, we obtain the final semidiscrete evolution equation
\begin{equation}\label{eq:governing equations - dG weak form - discrete}
\bm{M}^e\dot{\bm{X}}{}^e = 
\int_{\mathscr{D}^e}\frac{\pd\bm{B}^{e\Tr}}{\pd x_k}\bm{F}_k\dd{V}-
\int_{\mathscr{B}^{e}\cup\mathscr{L}^e}\bm{B}^{e\Tr}\what{\bm{F}}_n\dd{S}-
\sum_{e'\in\mathcal{N}^{e}}\int_{\mathscr{I}^{e,e'}}\bm{B}^{e\Tr}\what{\bm{F}}_n\dd{S},
\end{equation}
where the superimposed dot denotes the derivative with respect to time, and $\bm{M}^e$ is the mass matrix of $e$-th element given by
\begin{equation}\label{eq:mass matrix}
\boldsymbol{M}^{e} \equiv \int_{\mathscr{D}^{e}}\bm{B}^{e\Tr}\bm{B}^{e}~\dd{V}.
\end{equation}

The domain and boundary integrals appearing in Eqs.(\ref{eq:governing equations - dG weak form - discrete}) and (\ref{eq:mass matrix}) are evaluated using high-order quadrature rules for implicitly-defined domains and boundaries.
These quadrature rules are generated using the algorithm developed by Saye \cite{saye2015high} and provide high-order accuracy in the presence of smooth embedded geometries.
This is a key ingredient for the accuracy of the overall scheme.
For example, Qin and Krivodonova \cite{qin2013discontinuous} discuss the importance of an accurate representation of embedded geometries and use curvature-aware boundary conditions to avoid the loss of accuracy due to the piecewise representation of the embedded boundaries.
An example of the quadrature rules used here is provided in Fig.(\ref{fig:IBVP - QUADRATURE}), which shows the distribution of the quadrature points for the domain and boundary integrals of 2D and 3D implicitly-defined elements.

\subsection{Finite volume schemes}\label{sec:FV}
In a FV setting, the weak form of the governing equations can be derived by integrating Eq.(\ref{eq:governing equations - vector form}) over a generic implicitly-defined element $\mathscr{D}^{e} \in \mathcal{M}_h$ and applying the divergence theorem with respect to the space variables to obtain
\begin{equation}\label{eq:governing equations - FV weak form}
m^e\dot{\bm{U}}{}^e = -
\int_{\mathscr{B}^{e}\cup\mathscr{L}^e}\what{\bm{F}}_n\dd{S}-
\sum_{e'\in\mathcal{N}^{e}}\int_{\mathscr{I}^{e,e'}}\what{\bm{F}}_n\dd{S},
\end{equation}
where the superimposed dot again denotes the derivative with respect to time, $m^e \equiv \int_{\mathscr{D}^{e}}1~\dd{V}$ is the volume of $\mathscr{D}^{e}$, $\bm{U}^e$ represents the average of the conserved variables over $\mathscr{D}^{e}$, i.e.~$\bm{U}^{e} \equiv \frac{1}{m^e}\int_{\mathscr{D}^{e}} \bm{U}\dd{V}$, and the remaining symbols have the same meaning as those appearing in Eq.(\ref{eq:governing equations - dG weak form}).

Similar to Eq.(\ref{eq:governing equations - dG weak form - discrete}), high-order quadrature rules are employed to evaluate the right-hand side of Eq.(\ref{eq:governing equations - FV weak form}) and maintain the high-order accuracy representation of the geometry for the FV scheme.
However, unlike dG schemes, where the solution at any space location within the element is computed via its polynomial representation, Eq.(\ref{eq:governing equations - FV weak form}) requires a suitable reconstruction and limiting strategy, which allows one to evaluate the solution, and thus the numerical flux $\what{\bm{F}}_n$, at the element's boundaries.

Although high-order reconstruction methods are available in the literature, the present approach is based on a linear representation of the solution, which leads to a scheme that is at most second-order accurate in space.
This choice is mainly motivated by the fact that, within the adaptive mesh refinement algorithm, the FV grids will be fine-tuned to track and resolve the regions containing solution discontinuities, where high-order schemes would be reduced to first-order schemes by the limiter.
The details of the limiter employed here are given in Sec.(\ref{ssec:RECONSTRUCTION}).

\begin{figure}
\centering
\includegraphics[width = \textwidth]{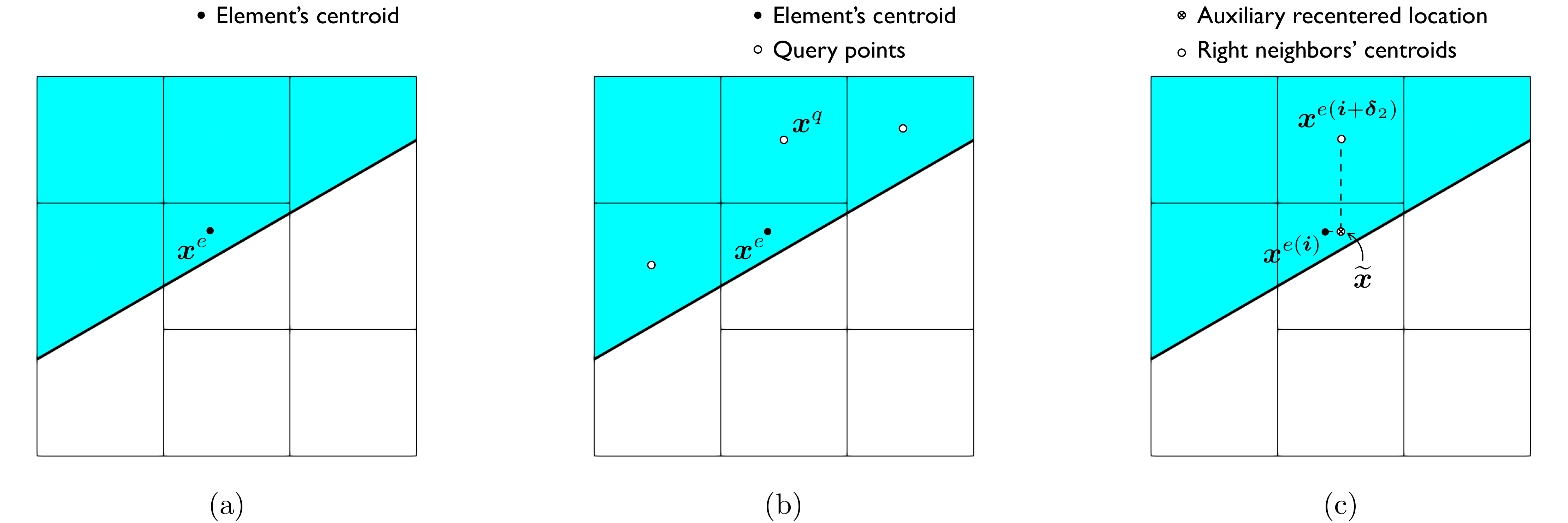}
\caption
{
    Reconstruction and limiting for an element belonging to an implicitly-defined mesh in a 3$\times$3 grid:
    (a) Location of the centroid $\bm{x}^{e}$ of the element where the reconstruction is to be performed;
    (b) Cloud of query points used in the least-square problem;
    (c) Auxiliary recentered location $\wtilde{\bm{x}}$ of the centroid $\bm{x}^{e(\bm{i})}$ such that $\wtilde{\bm{x}}$ is aligned with the centroid of the element $e(\bm{i}+\bm{\delta}_2)$.
}
\label{fig:IBVP - RECONSTRUCTION AND LIMITING}
\end{figure}

\subsubsection{Reconstruction and limiting}\label{ssec:RECONSTRUCTION}
In FV schemes, the typical reconstruction and limiting strategy is to compute the slopes of the primitive variables, perform the limiting of the slopes using characteristic variables and then compute the flux $\hat{\bm{F}}_n$ appearing in Eq.(\ref{eq:governing equations - FV weak form}) using the reconstructed primitive variables.
Here, we follow a slightly different approach, whereby the limiting is still performed using characteristics variables but the final reconstruction is performed on the conserved variables.
This preserves the robustness of characteristics-based limiters but enables a more natural coupling with the dG schemes, where conserved variables are readily available at the elements' boundaries.

The vector $\wtilde{\bm{U}}{}^{e}$ of the reconstructed conserved variables at any point $\bm{x}\in\mathscr{D}^{e}$ is assumed to be given by
\begin{equation}\label{eq:reconstructed solution}
\wtilde{\bm{U}}{}^{e} = 
\bm{U}^{e}+\pd_k\bm{U}^{e}\cdot(x_k-x_k^{e}),
\end{equation}
where $\bm{x}^{e} \equiv \{x_k^{e}\}$ is the centroid of $\mathscr{D}^{e}$, $\bm{U}^{e}$ is the vector of the conserved variable averages and $\pd_k\bm{U}^{e}$ are the slopes to be determined.
It is easy to see that the expression given in Eq.(\ref{eq:reconstructed solution}) preserves conservation, since $\int_{\mathscr{D}^e} \wtilde{\bm{U}}{}^{e} \dd{V}= m^e\bm{U}^{e}$, and that the conserved variable averages $\bm{U}^{e}$ coincide with the values of the conserved variables evaluated at the centroid of the element. 
As discussed next, $\pd_k\bm{U}^{e}$ are determined by reconstructing and limiting the slopes of the primitive variables.

Let $\bm{Q}^{e} \equiv \{\rho,\bm{v},p\}$ be the $(d+2)$-dimensional vector of the primitive variables corresponding to $\bm{U}^{e}$ and $Q^e$ be a generic component of $\bm{Q}^{e}$.
The initial slopes of $Q^e$ are obtained by the solution of the normal equations arising from the least-square interpolation problem given by
\begin{equation}\label{eq:least-square system}
\bm{S}\bm{\Delta} = \bm{b},
\end{equation}
where
\begin{equation}\label{eq:least-square components}
S_{kl} \equiv \sum_{q\in\mathcal{Q}}(x_k^q-x_k^e)(x_l^q-x_l^e),
\quad
\Delta_{k} \equiv \pd_k Q^e,
\quad\mathrm{and}\quad
b_{k} \equiv \sum_{q\in\mathcal{Q}}(x_k^q-x_k^e)(Q^q-Q^e).
\end{equation}
In (\ref{eq:least-square components}), $\mathcal{Q}$ denotes the set of neighboring query elements, whose centroid $\bm{x}^q$ and primitive variables $\bm{Q}^q$ are selected to perform the least-square interpolation.
Here, for an implicitly-defined element associated with the cell $\mathscr{C}^{\bm{i}}$, the set $\mathcal{Q}$ consists of the elements associated with the valid cells in the standard 4-cell, in 2D, or 6-cell, in 3D, neighborhood of $\mathscr{C}^{\bm{i}}$ plus the elements associated with the valid cells that extend into the small cells of the same neighborhood of $\mathscr{C}^{\bm{i}}$.
Figure (\ref{fig:IBVP - RECONSTRUCTION AND LIMITING}a) shows the case of a mesh implicitly defined on a $3\times{3}$ grid and Fig.(\ref{fig:IBVP - RECONSTRUCTION AND LIMITING}b) shows the corresponding cloud of query points (denoted by the open circles) that are selected to perform the least-square reconstruction over the element associated to the central cell.
It is worth noting that it is possible that the number of query points is not sufficient for the definition of the least-square problem given in Eq.(\ref{eq:least-square system}); in these cases, the slopes are set to zero.

Once the reconstruction is obtained a slope limiter is required to avoid non-physical oscillations and ensure monotonicity in the presence of solution discontinuities.
The limiter employed in this work consists of a two-step process.
The first step is based on a suitably modified version of the dimensionally-split characteristics-based Van Leer limiter \cite{van1979towards} combined with the recentering strategy suggested by Berger et al.\cite{berger2005analysis} to handle nonaligned centroids of adjacent elements.
More specifically, let $e(\bm{i})$ be the superscript denoting a quantity defined on the element $\mathscr{D}^e$ associated with the cell $\mathscr{C}^{\bm{i}}$, and let $\pd_k\bm{Q}^{e(\bm{i})}$ be the unlimited slopes in $k$-th direction obtained via Eq.(\ref{eq:least-square system}).
In the $k$-th direction, let $e(\bm{i}+\bm{\delta}_k)$ be the index of the \emph{right} neighbor of the element $e(\bm{i})$, and $e(\bm{i}-\bm{\delta}_k)$ be the index of the \emph{left} neighbor of element $e(\bm{i})$.
The values of the primitive variables in the right and left neighbors are then used to compute the corresponding right slopes $\pd_k^+\bm{Q}^{e(\bm{i})}$ and the left slopes $\pd_k^-\bm{Q}^{e(\bm{i})}$.
Consider for example $\pd_k^+\bm{Q}^{e(\bm{i})}$ and an auxiliary recentered location $\wtilde{\bm{x}}$ of the centroid $\bm{x}^{e(\bm{i})}$ obtained by translating $\bm{x}^{e(\bm{i})}$ in the plane perpendicular to the $k$-th direction and aligning it with the centroid of the right neighbor $e(\bm{i}+\bm{\delta}_k)$.
To illustrate, Fig.(\ref{fig:IBVP - RECONSTRUCTION AND LIMITING}c) sketches the location of $\wtilde{\bm{x}}$ when the right slopes are to be computed along the $x_2$ direction for the element associated to the central cell.
Using the unlimited slopes, the primitive variables $\wtilde{\bm{Q}}$ evaluated at $\wtilde{\bm{x}}$ are given by
\begin{equation}\label{eq:primitive variables at recentered centroid}
\wtilde{\bm{Q}} = \bm{Q}^{e(\bm{i})}+\pd_l\bm{Q}^{e(\bm{i})}\cdot(\wtilde{x}_l-x_l^{e(\bm{i})}) = \bm{Q}^{e(\bm{i})}+\sum_{l \ne k}\pd_l\bm{Q}^{e(\bm{i})}\cdot(x_l^{e(\bm{i}+\bm{\delta}_k)}-x_l^{e(\bm{i})})
\end{equation}
where the second equality follows by noting that $\wtilde{x}_k = x_k^{e(\bm{i})}$ and $\wtilde{x}_l = x_l^{e(\bm{i}+\bm{\delta}_k)}$, for $l \ne k$.
Then, the right slopes $\pd_k^+\bm{Q}^{e(\bm{i})}$ are defined as
\begin{equation}\label{eq:right slopes}
\pd_k^+\bm{Q}^{e(\bm{i})} \equiv \frac{\bm{Q}^{e(\bm{i}+\bm{\delta}_k)}-\wtilde{\bm{Q}}}{x_k^{e(\bm{i}+\bm{\delta}_k)}-x_k^{e(\bm{i})}}.
\end{equation}
Substituting Eq.(\ref{eq:primitive variables at recentered centroid}) into Eq.(\ref{eq:right slopes}), and following the same derivation for the left slopes, one obtains
\begin{equation}\label{eq:left and right derivatives}
\pd_k^{\pm}\bm{Q}^{e(\bm{i})} \equiv \frac{\bm{Q}^{e(\bm{i}\pm\bm{\delta}_k)}-\bm{Q}^{e(\bm{i})}}{x_k^{e(\bm{i}\pm\bm{\delta}_k)}-x_k^{e(\bm{i})}}-
\frac{\sum_{l \ne k}\pd_l\bm{Q}^{e(\bm{i})}\cdot(x_l^{e(\bm{i}\pm\bm{\delta}_k)}-x_l^{e(\bm{i})})}{x_k^{e(\bm{i}\pm\bm{\delta}_k)}-x_k^{e(\bm{i})}}.
\end{equation}
In Eqs.(\ref{eq:right slopes}) and (\ref{eq:left and right derivatives}), no summation is implied over repeated subscripts.
Moreover, it is worth noting that, far from the embedded boundaries, where the elements and their neighbors are regular entire cells, the slopes given in Eq.(\ref{eq:left and right derivatives}) coincide with the slopes of standard structured grid methods since, along the $k$-th direction, one has $x_l^{e(\bm{i}\pm\bm{\delta}_k)} = x_l^{e(\bm{i})}$, for $l \ne k$.
This is not the case near the embedded boundary where the term $x_l^{e(\bm{i}\pm\bm{\delta}_k)}-x_l^{e(\bm{i})}$ accounts for the difference in location of the elements' centroids.
However, as suggested by Fig.(\ref{fig:IBVP - RECONSTRUCTION AND LIMITING}c), it is possible that some of the neighboring elements might not be available because the neighboring cells are small or empty according to the classification introduced in Sec.(\ref{ssec:IMPLICIT MESH}).
In these cases, the left or the right derivative cannot be computed and will not be considered in the evaluation of the limited slopes.

The slopes of the primitive variables are used to compute the corresponding slopes of the characteristics variables as
\begin{equation}\label{eq:left and right derivatives - characteristics}
\pd_k\bm{C}^{e(\bm{i})} \equiv \bm{L}^{k}\pd_k\bm{Q}^{e(\bm{i})}
\quad\mathrm{and}\quad
\pd_k^{\pm}\bm{C}^{e(\bm{i})} \equiv \bm{L}^k\pd_k^{\pm}\bm{Q}^{e(\bm{i})},
\end{equation}
where $\bm{L}^{k}$ is the matrix of left eigenvectors of $(\pd\bm{Q}/\pd\bm{U})(\pd\bm{F}_k/\pd\bm{U})(\pd\bm{U}/\pd\bm{Q})$ evaluated at $\bm{U}^{e(\bm{i})}$.
Then, the vector $\overline{\pd}_k\bm{C}^{e(\bm{i})}$ of limited slopes of the characteristic variables is obtained by applying the Var Leer limiter 
\begin{equation}\label{eq:Van Leer limiter}
\overline{\pd}_kC^{e(\bm{i})} \equiv
\begin{cases}
s\min\{\theta^-|\delta^-|,|\delta|,\theta^+|\delta^+|\},&\mathrm{if}~\mathscr{C}^{\bm{i}\pm\bm{\delta}_k}~\mathrm{are~valid~and}~\delta^-,\delta,\delta^+~\mathrm{have~the~same~sign};\\
s\min\{\theta^-|\delta^-|,|\delta|\},&\mathrm{if}~\mathscr{C}^{\bm{i}-\bm{\delta}_k}~\mathrm{is~valid,}~\mathscr{C}^{\bm{i}+\bm{\delta}_k}~\mathrm{is~not~valid~and}~\delta^-,\delta~\mathrm{have~the~same~sign};\\
s\min\{|\delta|,\theta^+|\delta^+|\},&\mathrm{if}~\mathscr{C}^{\bm{i}+\bm{\delta}_k}~\mathrm{is~valid,}~\mathscr{C}^{\bm{i}-\bm{\delta}_k}~\mathrm{is~not~valid~and}~\delta,\delta^+~\mathrm{have~the~same~sign};\\
0,&\mathrm{otherwise}
\end{cases}
\end{equation}
component-wise.
In Eq.(\ref{eq:Van Leer limiter}), $\overline{\pd}_kC^{e(\bm{i})}$ denotes a generic component of $\overline{\pd}_k\bm{C}^{e(\bm{i})}$, $\delta \equiv \pd_kC^{e(\bm{i})}$, $\delta^{\pm} \equiv \pd_k^{\pm}C^{e(\bm{i})}$, $s \equiv \mathrm{sign}(\delta)$, $\theta^- \equiv (x_k^{e(\bm{i})}-x_k^{e(\bm{i}-\bm{\delta}_k)})/(x_k^{e(\bm{i})}-x_k^{\bm{i}})$, $\theta^+ \equiv (x_k^{e(\bm{i}+\bm{\delta}_k)}-x_k^{e(\bm{i})})/(x_k^{\bm{i}}+h_k-x_k^{e(\bm{i})})$, with no summation implied over the repeated index $k$.
Finally, the limited slopes of the primitive variables and the corresponding slopes of the conserved variables are obtained as
\begin{equation}\label{eq:limited slopes}
\overline{\pd}_k\bm{Q}^{e} = (\bm{L}^{k})^{-1}\overline{\pd}_k\bm{C}^{e}
\quad\mathrm{and}\quad
\overline{\pd}_k\bm{U}^{e} = \left(\frac{\pd\bm{U}}{\pd\bm{Q}}\right)_{\bm{U}^e}\overline{\pd}_k\bm{Q}^{e}
\end{equation}
respectively.

The slopes computed using Eq.(\ref{eq:limited slopes}) do not ensure that over- or under-shoots of the solution are avoided at the boundaries of the mesh elements, especially in the case of extended elements where left or right derivatives are not available.
For this reason, the second step of the limiter applies the Barth-Jespersen limiter \cite{barth1989design}, which is a scalar limiter that reduces the slopes along all directions by the same multiplicative factor $\alpha$.
Thus, the slopes entering Eq.(\ref{eq:reconstructed solution}) are defined as $\pd_k\bm{U}^{e} \equiv \alpha^e\overline{\pd}_k\bm{U}^{e}$, where $\alpha^e$ is computed via a Barth-Jespersen limiter procedure that uses the reconstructed solution at quadrature point on the boundaries $\mathscr{I}^{e,e'}$ $\forall e'\in\mathcal{N}^e$.

\subsection{Block-structured adaptive mesh refinement}\label{ssec:AMR}
\begin{figure}
\centering
\includegraphics[width = \textwidth]{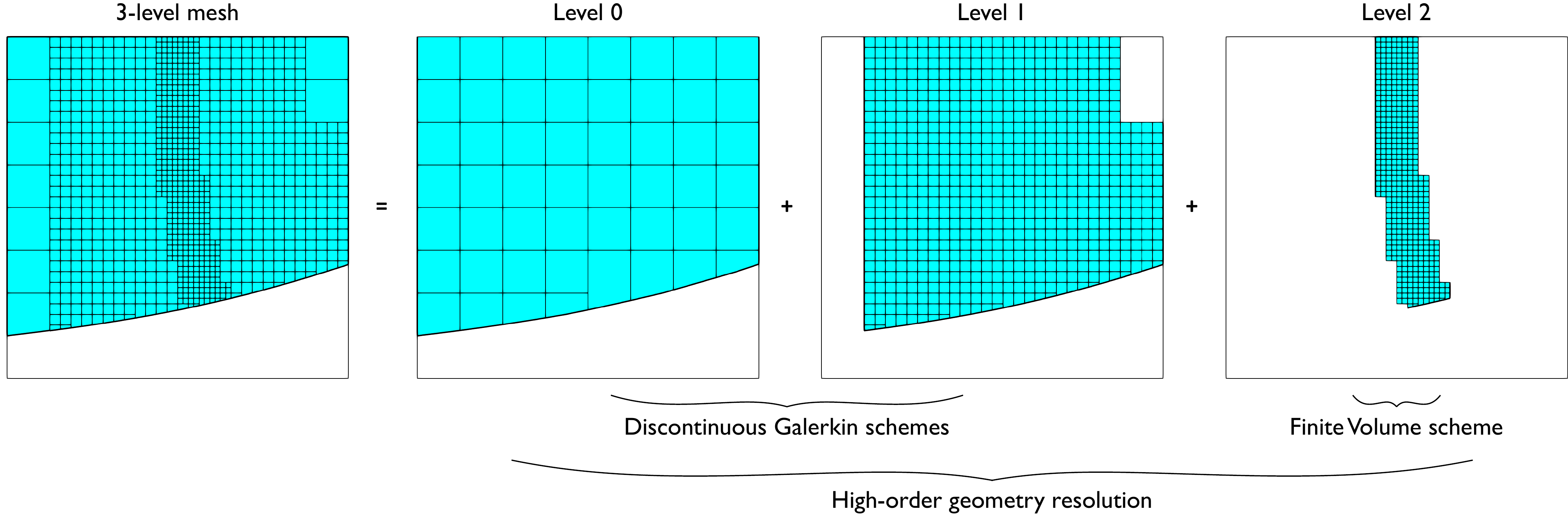}
\caption
{
    Example of a three-level implicitly-defined mesh in which the levels 0 and 1 use discontinuous Galerkin schemes and the level 2 uses a finite-volume scheme.
    At all levels, the geometry is resolved with high-order accuracy.
}
\label{fig:IBVP - AMR - EXAMPLE}
\end{figure}
\begin{figure}
\centering
\includegraphics[width = \textwidth]{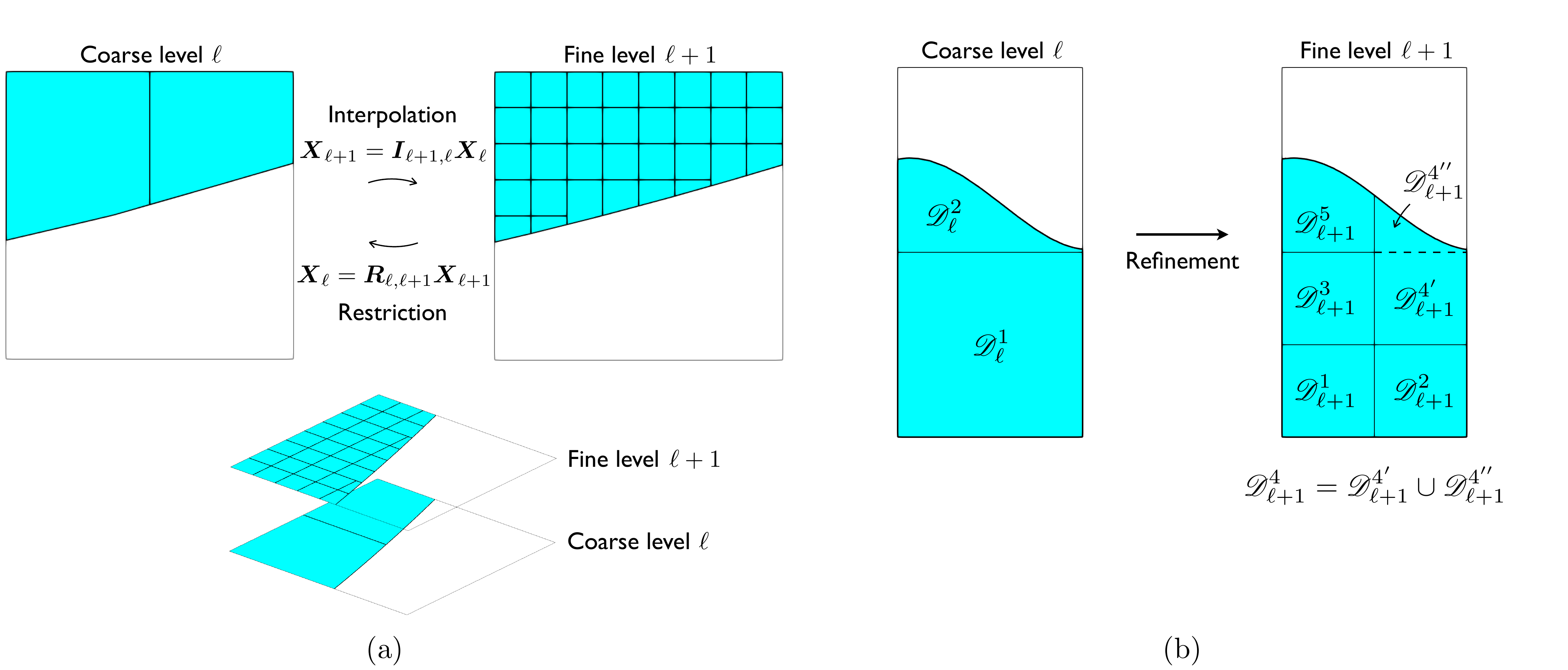}
\caption
{
    (a) Illustration of the operations performed by the interpolation operator $\bm{I}_{\ell+1,\ell}$ and the restriction operator $\bm{R}_{\ell+1,\ell}$.
    (b) Example of an AMR configuration where the fine element $\mathscr{D}_{\ell+1}^4$ covers portions of the two distinct coarse elements $\mathscr{D}_{\ell}^1$ and $\mathscr{D}_{\ell}^2$.
}
\label{fig:IBVP - AMR - DETAILS}
\end{figure}

The dG schemes and FV schemes are integrated into a block-structured adaptive mesh refinement algorithm (AMR) \cite{berger1989local,bell1994three,zhang2019amrex} that represents the solution as a hierarchy of levels of refinement ranging from coarsest ($\ell = 0$) to finest ($\ell = {\ell}_{\rm finest}$).
Each level is represented by the union of non-overlapping logically-rectangular grids of a given resolution.

Here, the grids at level $\ell = 0$ are generated at the beginning of the simulation and cover the entire domain; the corresponding implicitly-defined mesh is also generated at startup and does not change during the simulation.
Grids at levels $\ell > 0$ are created and destroyed dynamically during the simulation, and the corresponding implicitly-defined meshes follow the same dynamic evolution.
We consider AMR configurations where dG schemes are used at the levels $0 \le \ell < {\ell}_{\rm finest}$ while the FV scheme is always used at the finest level $\ell = {\ell}_{\rm finest}$.
This includes the case where AMR uses dG schemes only or the case where AMR consists of a single-level FV scheme.
The present AMR strategy also supports the use of different order basis functions at different dG levels and, as such, is referred to as $hp$-AMR.
Figure (\ref{fig:IBVP - AMR - EXAMPLE}) shows an example of a 3-level $hp$-AMR configuration with a refinement ratio of 4 between level 0 and level 1 and a refinement ratio of 2 between level 1 and level 2.
The figure also shows that dG schemes are used at levels 0 and 1 and a FV scheme is used at the level 2.
This configuration is typical of what one would expect for a discontinuity attached to a curved boundary.

The dynamic update of the AMR levels is governed by tagging and un-tagging operations, whereby a cell at a given level is assigned a two-value variable, or tag, that determines whether a finer resolution is required at its location or whether the existing set of finer grid cells can be replaced by a coarse cell.
These operations require coarsening and refinement operations that transfer the solution between levels.
For instance, when a coarse cell is tagged for refinement and a finer grid is overlaid on it, the solution from the coarse level needs to be interpolated to the newly created fine grid.
Conversely, the solution at a fine level must be restricted to the corresponding coarser level to determine whether or not the fine grid is still required at that location.
These operations must preserve the high-order accuracy of the method and must be available for all the possible combinations of the numerical schemes used between consecutive levels.
Nevertheless, these are linear operations local to elements and therefore can be efficiently implemented using matrix-vector products taking advantage of parallelization offered by dG and FV schemes.

To illustrate how interpolation and restriction are performed here, we introduce the vector $\bm{X}_{\ell}$ of the solution at the level $\ell$.
For a level that uses a dG scheme, $\bm{X}_{\ell}$ contains the coefficients of the basis functions representing the conserved variables, whereas, for a level that uses the FV scheme, $\bm{X}_{\ell}$ collects the element averages of the conserved variables.
For the interpolation and restriction operations, we consider the FV scheme as a dG scheme of discontinuous constant basis functions.
This means that interpolation to the finest level returns the element averages, which are subsequently used to perform the reconstruction and limiting process presented in Sec.(\ref{ssec:RECONSTRUCTION}).
Note that interpolation is never performed from the finest level and therefore the reconstructed slopes are
not used during the interpolation and restriction operations.

As sketched in Fig.(\ref{fig:IBVP - AMR - DETAILS}a), the interpolation operator $\bm{I}_{\ell+1,\ell}$ computes the solution $\bm{X}_{\ell+1}$ at the fine level $\ell+1$ from the solution $\bm{X}_{\ell}$ at the coarse level $\ell$, and is formally written as
\begin{equation}\label{eq:interpolation}
\boldsymbol{X}_{\ell+1} = \boldsymbol{I}_{\ell+1,\ell}\boldsymbol{X}_{\ell}.
\end{equation}
In Eq.(\ref{eq:interpolation}), $\bm{I}_{\ell+1,\ell}$ is a block-structured matrix computed using a Galerkin projection, so that, given an element $\mathscr{D}_{\ell+1}^{e'}$ of the fine level $\ell+1$ and an element $\mathscr{D}_{\ell}^{e}$ of the coarse level $\ell$ for which the interpolation operation is meaningful, the block $\boldsymbol{I}_{\ell+1,\ell}^{e',e}$ is given by
\begin{equation}\label{eq:interpolation - Galerkin projection}
\boldsymbol{I}_{\ell+1,\ell}^{e',e} \equiv (\boldsymbol{M}_{\ell+1}^{e'})^{-1}\int_{\mathscr{D}_{\ell+1}^{e'}\cap\mathscr{D}_{\ell}^{e}}\bm{B}_{\ell+1}^{e'\Tr}\bm{B}_{\ell}^{e}~\dd{V},
\end{equation}
where, for an element $e$ on level $\ell$, $\bm{B}_{\ell}^{e}$ contains the basis functions and $\boldsymbol{M}_{\ell}^{e}$ is the mass matrix as defined in Eq.(\ref{eq:mass matrix}).

It is worth giving a practical example of how the interpolation operator $\boldsymbol{I}_{\ell+1,\ell}$ acts on $\bm{X}_\ell$.
Consider Fig.(\ref{fig:IBVP - AMR - DETAILS}b), which shows a two-element coarse mesh $\mathcal{M}_{h_\ell} \equiv \{\mathscr{D}_{\ell}^1,\mathscr{D}_{\ell}^2\}$ associated with a $1\times 2$ coarse grid, and a five-element fine mesh $\mathcal{M}_{h_{\ell+1}} \equiv \{\mathscr{D}_{\ell+1}^1,\mathscr{D}_{\ell+1}^2,\dots,\mathscr{D}_{\ell+1}^5\}$ associated with a $2\times 4$ grid that is obtained by refining the coarse grid with a refinement ratio of 2.
Then, $\bm{X}_\ell$ consists of the basis functions coefficients of the elements $\mathscr{D}_\ell^1$ and $\mathscr{D}_\ell^2$, $\bm{X}_{\ell+1}$ consists of the basis functions coefficients of the elements $\mathscr{D}_{\ell+1}^1$ to $\mathscr{D}_{\ell+1}^5$.
The interpolation operation given in Eq.(\ref{eq:interpolation}) becomes
\begin{equation}\label{eq:interpolation - example}
\left\{\begin{array}{c}
\bm{X}_{\ell+1}^{1}\\
\bm{X}_{\ell+1}^{2}\\
\bm{X}_{\ell+1}^{3}\\
\bm{X}_{\ell+1}^{4}\\
\bm{X}_{\ell+1}^{5}
\end{array}\right\}
=
\left[\begin{array}{cc}
\bm{I}_{\ell+1,\ell}^{1,1}&\bm{0}\\
\bm{I}_{\ell+1,\ell}^{2,1}&\bm{0}\\
\bm{I}_{\ell+1,\ell}^{3,1}&\bm{0}\\
\bm{I}_{\ell+1,\ell}^{4,1}&\bm{I}_{\ell+1,\ell}^{4,2}\\
\bm{0}&\bm{I}_{\ell+1,\ell}^{5,2}
\end{array}\right]
\left\{\begin{array}{c}
\bm{X}_{\ell}^{1}\\
\bm{X}_{\ell}^{2}
\end{array}\right\},
\end{equation}
where each block $\bm{I}_{\ell+1,\ell}^{e,e'}$, with $e = 1,\dots,5$ and $e' = 1,2$, is computed via Eq.(\ref{eq:interpolation - Galerkin projection}).
While the first, second, third and fifth rows of Eq.(\ref{eq:interpolation - example}) are the result of a trivial application of the Galerkin projection, it is interesting to derive the interpolation operation given in the fourth row of Eq.(\ref{eq:interpolation - example}).
First, with reference to Fig.(\ref{fig:IBVP - AMR - DETAILS}b), note that the fine element $\mathscr{D}_{\ell+1}^4$ is the result of a merging operation and covers portions of the two distinct coarse elements $\mathscr{D}_{\ell}^1$ and $\mathscr{D}_{\ell}^2$.
Then, assuming for simplicity that $\bm{X}_\ell$ contains the basis functions coefficients of a discontinuous scalar field $u: \mathcal{M}_{h_{\ell}} \rightarrow \mathbb{R}$ at the coarse level $\ell$, the Galerkin projection of $u$ onto the fine element $\mathscr{D}_{\ell+1}^4$ is given by
\begin{multline}\label{eq:interpolation - example - non-trivial}
\bm{X}_{\ell+1}^4 = (\bm{M}_{\ell+1}^4)^{-1}\int_{\mathscr{D}_{\ell+1}^4}\bm{B}_{\ell+1}^{4\Tr}u\dd V = \\
= (\bm{M}_{\ell+1}^4)^{-1}\left(\int_{\mathscr{D}_{\ell+1}^{4'}}\bm{B}_{\ell+1}^{4\Tr}\bm{B}_{\ell}^{1}\dd V\right) \bm{X}_{\ell}^1+
(\bm{M}_{\ell+1}^4)^{-1}\left(\int_{\mathscr{D}_{\ell+1}^{4''}}\bm{B}_{\ell+1}^{4\Tr}\bm{B}_{\ell}^{2}\dd V\right) \bm{X}_{\ell}^2
= \bm{I}_{\ell+1,\ell}^{4,1} \bm{X}_{\ell}^1+
\bm{I}_{\ell+1,\ell}^{4,2} \bm{X}_{\ell}^2
\end{multline}
where the last equality is obtained by noting that $\mathscr{D}_{\ell+1}^{4'} \equiv \mathscr{D}_{\ell+1}^{4}\cap\mathscr{D}_{\ell}^{1}$ and $\mathscr{D}_{\ell+1}^{4''} \equiv \mathscr{D}_{\ell+1}^{4}\cap\mathscr{D}_{\ell}^{2}$, and by using Eq.(\ref{eq:interpolation - Galerkin projection}).
The structure of $\bm{I}_{\ell+1,\ell}$ given in Eq.(\ref{eq:interpolation - example}) reflects the configurations of the coarse and the fine meshes displayed in Fig.(\ref{fig:IBVP - AMR - DETAILS}b), including the merging between cells; moreover, if $\mathscr{D}_{\ell+1}^{e'}\cap\mathscr{D}_{\ell}^{e} = \emptyset$, then the corresponding block $\bm{I}_{\ell+1,\ell}^{e,e'}$ is a zero matrix.

Figure (\ref{fig:IBVP - AMR - DETAILS}b) also sketches the inverse of the interpolation operation, namely the restriction operation, where a solution $\bm{X}_\ell$ at a coarse level $\ell$ is computed from the solution $\bm{X}_{\ell+1}$ at a fine level $\ell+1$ using the restriction operator $\bm{R}_{\ell,\ell+1}$
\begin{equation}\label{eq:restriction}
\boldsymbol{X}_{\ell} = \boldsymbol{R}_{\ell,\ell+1}\boldsymbol{X}_{\ell+1}.
\end{equation}
Following Fortunato et al.\cite{fortunato2019efficient}, the restriction operator can be defined as the adjoint of the interpolation operator and can evaluated as follows
\begin{equation}\label{eq:restriction - from interpolation}
\boldsymbol{R}_{\ell,\ell+1} = \bm{M}_{\ell}^{-1}\boldsymbol{I}_{\ell+1,\ell}^{\Tr}\bm{M}_{\ell+1},
\end{equation}
where $\bm{M}_{\ell}$ denotes the block-diagonal mass matrix associated to the entire level $\ell$, i.e.~$\bm{M}_{\ell} \equiv \mathrm{diag}(\{\bm{M}_\ell^e\})$.
Note that Eq.(\ref{eq:restriction - from interpolation}) can also be obtained by applying the Galerkin projection from the fine level $\ell+1$ onto the coarse level $\ell$ and using Eq.(\ref{eq:interpolation - Galerkin projection}).
For the case shown in Fig.(\ref{fig:IBVP - AMR - DETAILS}b), Eq.(\ref{eq:restriction - from interpolation}) becomes
\begin{equation}\label{eq:restriction - example}
\left\{\begin{array}{c}
\bm{X}_{\ell}^{1}\\
\bm{X}_{\ell}^{2}
\end{array}\right\} = \bm{M}_\ell^{-1}
\left[\begin{array}{ccccc}
\bm{I}_{\ell+1,\ell}^{1,1\Tr}&\bm{I}_{\ell+1,\ell}^{2,1\Tr}&\bm{I}_{\ell+1,\ell}^{3,1\Tr}&\bm{I}_{\ell+1,\ell}^{4,1\Tr}&\bm{0}\\
\bm{0}&\bm{0}&\bm{0}&\bm{I}_{\ell+1,\ell}^{4,2\Tr}&\bm{I}_{\ell+1,\ell}^{5,2\Tr}
\end{array}\right]\bm{M}_{\ell+1}
\left\{\begin{array}{c}
\bm{X}_{\ell+1}^{1}\\
\bm{X}_{\ell+1}^{2}\\
\bm{X}_{\ell+1}^{3}\\
\bm{X}_{\ell+1}^{4}\\
\bm{X}_{\ell+1}^{5}
\end{array}\right\}.
\end{equation}

A few comments are worth pointing out:
\begin{itemize}
    \item The interpolation and restriction operators defined via Eqs.(\ref{eq:interpolation - Galerkin projection}) and (\ref{eq:restriction}), respectively, are valid regardless of the choice of  basis functions and can accommodate AMR levels using different mesh sizes, polynomial orders or both.
    In particular, they are valid even if the same mesh uses the same space of basis functions with different polynomial orders; this feature will be used to define the shock sensor discussed in Sec.(\ref{ssec:CELL TAGGING}).
    \item From an implementation perspective, given an $hp$-AMR configuration, one only needs to evaluate the interpolation operators and the Cholesky decomposition of the mass matrices; then the restriction operator can be applied on-the-fly using Eq.(\ref{eq:restriction - from interpolation}).
    Moreover, all the standard rectangular elements (which represent most of the mesh elements) share the same mass matrix, which can precomputed and stored at the beginning of the simulation; the same applies to the interpolation operator between two standard rectangular elements of two different AMR levels. 
    For the remaining implicitly-defined mesh elements, the mass matrices and the interpolation operators are in general unique and are computed using high-order quadrature rules \cite{saye2015high}.
    \item The mass matrix and interpolation operator are formally given by $\bm{M}_{\ell}$ and $\boldsymbol{I}_{\ell+1,\ell}$ that dynamically change as the $hp$-AMR levels evolve.
    However, these matrices are composed of localized blocks that can be computed and stored after each regridding.
    \item The interpolation and restriction operations are well-defined in all configurations of the implicitly-defined mesh obtained from the cell-merging algorithm as long as tagged cells are not merged with non-tagged cells. 
    This situation can always be avoided by making sure that if an extended cell is tagged for refinement, so are all the small cells merged with it.
\end{itemize}

As the last comment on the present $hp$-AMR strategy, recall that when the right-hand side of Eq.(\ref{eq:governing equations - dG weak form - discrete}) or (\ref{eq:governing equations - FV weak form}) is to be evaluated, the coupling between neighboring elements $e$ and $e'$ occurs solely through the numerical flux $\what{\bm{F}}_n$ appearing in the boundary integrals over $\mathscr{I}^{e,e'}$.
This applies equally to the cases when $e$ and $e'$ belong to the same $hp$-AMR level or when $e$ and $e'$ belong to different $hp$-AMR levels, including the level using the FV scheme.
In fact, at a generic quadrature point $\bm{x}$ of $\mathscr{I}^{e,e'}$, $\what{\bm{F}}_n$ is computed as $\what{\bm{F}}_n = \what{\bm{F}}_n(\wtilde{\bm{U}}{}^{e},\wtilde{\bm{U}}{}^{e'})$, where $\wtilde{\bm{U}}{}^{e}$ is the solution evaluated at $\bm{x}$ using either its representation in terms of basis functions, if the element $e$ belongs to a level using a dG scheme, or the reconstruction provided by Eq.(\ref{eq:reconstructed solution}), if the element $e$ belongs to a level using the FV scheme.
However, is it pointed out that at the boundary between two levels using different schemes, the order of the quadrature rules must be high enough to provide an accurate integration of the contributions from the higher-order scheme.

\subsubsection{Cell tagging and shock sensor}\label{ssec:CELL TAGGING}
Two main criteria are employed to tag the elements for refinement and update the AMR levels during the simulations.
The first one is a simple criterion that uses the magnitude of the density gradient.
In particular, an element $e$ of the AMR level $\ell$ is tagged for refinement if the following condition is satisfied
\begin{equation}\label{eq:density gradient refinement}
\frac{1}{m^e}\int_{\mathscr{D}^e}||\nabla\rho||~\dd{V} > \kappa_\ell^\rho
\end{equation}
where $m^e$ is the volume of $\mathscr{D}^e$, $||\nabla\rho||$ denotes the magnitude of the density gradient and $\kappa_\ell^\rho$ is a threshold parameter that determines which elements of the $\ell$-th level are tagged for refinement.
In general, a more advanced tagging criterion, involving for example the entire set of conserved variables and/or higher order derivatives, is recommended to capture a wider variety of flow configurations, such as discontinuous shear flows.
However, for the simulations considered here, the criterion given in Eq.(\ref{eq:density gradient refinement}) was tuned to use a low-order dG scheme in regions of constant flow and high-order dG schemes in regions of smooth flows.
Then, cells at the dG level $\ell = \ell_{\mathrm{finest}}-1$ are tagged for refinement and replaced by the cells at the FV level $\ell = \ell_{\mathrm{finest}}$ using the second criterion discussed next.

The second criterion is the so-called shock-sensor, which is used to determine whether the solution inside an element is behaving like a discontinuous field.
Here, we employ the shock sensor introduced by Persson and Peraire \cite{persson2006sub} suitably modified to account for the implicitly-defined mesh.
In Ref.\cite{persson2006sub}, the Authors defined the following smoothness indicator $s^e$ for the $e$-th mesh element
\begin{equation}\label{eq:smoothness indicator}
s^e \equiv \frac{\int_{\mathscr{D}^e}(U-\wtilde{U})^2~\dd{V}}{\int_{\mathscr{D}^e}U^2~\dd{V}}
\end{equation}
where $U$ is a field of the conserved variables, e.g.~the density, and $\wtilde{U}$ is the same field expressed using a lower-order set of basis functions.
For standard rectangular elements, $\wtilde{U}$ can be computed by considering a truncated series of the basis functions containing the terms up to order $p-1$, where $p$ is the order of the basis functions used to represent $U$.
However, for implicitly-defined elements, the truncated series would not preserve the average of the field $U$ in $\wtilde{U}$ and would degrade the efficacy of Eq.(\ref{eq:smoothness indicator}).
Therefore, at the level $\ell$ using the space $\mathcal{V}_{hp}$, $\wtilde{U}$ is evaluated using the Galerkin projection of the field $U$ onto the space $\mathcal{V}_{h(p-1)}$.
A discontinuity is then assumed to be present if an element $e$ if the following condition is satisfied
\begin{equation}\label{eq:shock sensor}
\log_{10}(s^e) > -4\log_{10}(p)+\kappa^s,
\end{equation}
where $p$ is the order of the polynomial basis functions employed at the level containing the element $e$, and $\kappa^s$ is a threshold parameter that determines which elements are tagged for refinement.
The effect of the parameter $\kappa^s$ on the evolution of the $hp$-AMR levels will be discussed in the numerical examples.
It is worth noting that in Ref.\cite{persson2006sub} the Authors employed the smoothness sensor given in Eq.(\ref{eq:smoothness indicator}) to inject a suitably-chosen artificial viscosity to the governing equations.
Here, it is employed only to activate the AMR level that uses the FV scheme.
We note that the refinement conditions given in (\ref{eq:density gradient refinement}) and (\ref{eq:shock sensor}) are used not only to tag the elements for refinement but also to remove the finer grids on regions where they are no longer needed.
\section{Numerical tests}\label{sec:RESULTS}
In this section, the capabilities of the proposed method are assessed for two- and three-dimensional test cases.
In all computations, an approximate two-shock Riemann solver \cite{colella1985efficient} is used to compute the numerical flux $\what{\bm{F}}_n$ appearing in Eqs.(\ref{eq:governing equations - dG weak form - discrete}) and (\ref{eq:governing equations - FV weak form}).
Tensor-product Legendre polynomials of degree $p$ are used to define the spaces $\mathcal{P}_{hp}^e$ and $\mathcal{V}_{hp}$ and the corresponding dG scheme is labelled as dG$_p$.

Whether a single level or an $hp$-AMR is used, the time evolution of the solution unknowns associated to the generic mesh element can be written as
\begin{equation}\label{eq:ode}
\bm{M}^e\dot{\bm{X}}{}^e = \bm{A}^e(t,\bm{X}),
\end{equation}
where
$\bm{A}^e(t,\bm{X})$ is the result of the evaluation of the right-hand side of Eq.(\ref{eq:governing equations - dG weak form - discrete}) (for dG schemes) or Eq.(\ref{eq:governing equations - FV weak form}) (for FV schemes), 
$\bm{X}$ contains the solution unknowns of the whole problem, and,
if the element $e$ belongs to a level implementing the FV scheme, $\bm{M}^e \equiv \mathrm{diag}(m^e)$ and $\dot{\bm{X}}{}^e \equiv \dot{\bm{U}}{}^e$.
$\bm{X}$ formally contains the solution unknowns of all the mesh elements but only the neighbors of the $e$-th element are needed to compute $\bm{A}^e$ in Eq.(\ref{eq:ode}).

The integration of Eq.(\ref{eq:ode}) over all mesh elements of the $hp$-AMR hierarchy is performed explicitly in time using a high-order TVD Runge-Kutta (RK) algorithms \cite{cockburn1998runge}.
The order of the RK algorithm is chosen to match the highest spatial discretization order among the levels; for example, if a two-level mesh uses a dG$_1$ scheme and a dG$_2$ scheme, Eq.(\ref{eq:ode}) is integrated in time at both levels via a third-order RK algorithm.
The time step of the RK algorithm is chosen to be the smallest time step among the levels.
Specifically, at level $\ell$ covering a subregion $\mathscr{D}_\ell$ of the domain $\mathscr{D}$ with grids of characteristic mesh size $h_\ell$, the time step $\tau_\ell$ is subject to the CFL condition
\begin{equation}
\frac{\tau_\ell}{h_\ell} < \frac{\mathrm{C}_\ell\overline{\nu}_\ell}{\lambda_\ell},
\end{equation}
where $\lambda_\ell \equiv \mathrm{max}_{\mathscr{D}_\ell}(\sqrt{v_kv_k}+a)$ is the maximum wave speed on $\mathscr{D}_\ell$ given by the sum of the velocity magnitude $\sqrt{v_kv_k}$ and the speed of sound $a$, $\overline{\nu}_\ell$ is the volume fraction threshold that triggers the merging process for the small cells, and $\mathrm{C}_\ell = 0.3$ if the level $\ell$ uses a FV scheme or $\mathrm{C}_\ell = 1/(2p+1)$ if the level $\ell$ uses a dG$_p$ scheme.
Then, the time step of the Runge-Kutta algorithm is $\tau \equiv \min_{\ell}\tau_\ell$.
All the computations use $\gamma = 1.4$ and all reported quantities are non-dimensional.

\subsection{Supersonic vortex}\label{ssec:SUPERSONIC VORTEX}
\begin{figure}
\centering
\includegraphics[width = \textwidth]{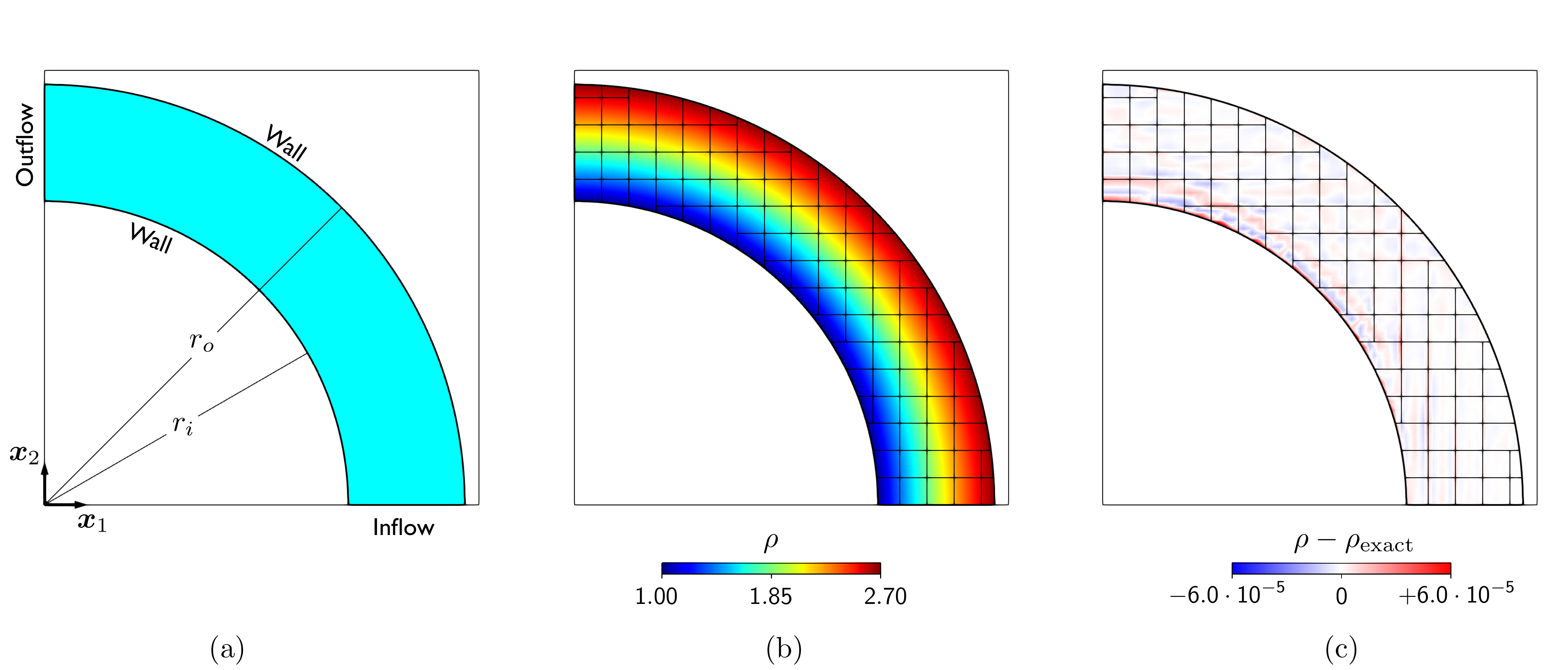}
\caption
{
    (a) Geometry of supersonic vortex problem defined in the background square $[0,1.43]^2$.
    (b) Density distribution and (c) density error obtained using a $16\times 16$ background mesh and a dG$_3$ scheme.
}
\label{fig:RESULTS - SUPERSONIC VORTEX - GEOMETRY-MESH-ERROR}
\end{figure}

In this test, we consider a two-dimensional supersonic vortex problem with an exact smooth solution that has been employed in several studies to investigate the accuracy of numerical methods for inviscid gas dynamics, see e.g.~\cite{krivodonova2006high,berger2021state}.
The problem consists of isentropic flow in a circular annulus of inner radius $r_i = 1$ and outer radius $r_o = 1.384$.
The exact solution for density $\rho_\mathrm{exact}$, pressure $p_\mathrm{exact}$, and velocity components $v_{1\mathrm{exact}}$ and $v_{2\mathrm{exact}}$ are given by the following expressions
\begin{equation}\label{eq:supersonic vortex problem - exact solution}
\rho_\mathrm{exact} = \rho_i\left[1+\frac{\gamma-1}{2}M_i^2\left(1-\frac{r_i^2}{r^2}\right)\right]^{\frac{1}{\gamma-1}},
\quad
p_\mathrm{exact} = \frac{\rho_\mathrm{exact}^\gamma}{\gamma},
\quad
v_{1\mathrm{exact}} = -v_\theta \sin{\theta}
\quad\mathrm{and}\quad
v_{2\mathrm{exact}} = v_\theta \cos{\theta},
\end{equation}
where $v_\theta = a_iM_i\frac{r_i}{r}$, $r$ is the distance between a point inside the annulus and the center of the annulus, and $\rho_i = 1$, $a_i = 1$ and $M_i = 2.25$ are the density, sound speed and Mach number, respectively, at $r = r_i$.

The problem is defined in the background square $\mathscr{R} = [0,1.43]^2$ where the domain $\mathscr{D}$ and its boundaries are represented by the level set function
\begin{equation}\label{eq:supersonic vortex problem - level set function}
\Phi(\bm{x}) \equiv
\left\{
\begin{array}{cc}
r_i^2-r^2&\mathrm{if}~r \le (r_i+r_o)/2 \\
r^2-r_o^2&\mathrm{otherwise}
\end{array}
\right.,
\quad\mathrm{with}\quad
r^2 = x_kx_k.
\end{equation}
Figure (\ref{fig:RESULTS - SUPERSONIC VORTEX - GEOMETRY-MESH-ERROR}a) shows the geometry of the problem and indicates the boundary conditions.
The exact solution given in Eq.(\ref{eq:supersonic vortex problem - exact solution}) is used to define the initial conditions for Eq.(\ref{eq:ode}).
The system is integrated in time until a steady state is reached.
We consider two error measures
\begin{equation}\label{eq:supersonic vortex problem - error measures}
e_{L_2}(\rho,\rho_\mathrm{exact}) \equiv \frac{||\rho-\rho_\mathrm{exact}||_2}{||\rho_\mathrm{exact}||_2}
\quad\mathrm{and}\quad
e_{L_\infty}(\rho,\rho_\mathrm{exact}) \equiv \frac{||\rho-\rho_\mathrm{exact}||_\infty}{||\rho_\mathrm{exact}||_\infty},
\end{equation}
where $||\cdot||_2$ and $||\cdot||_\infty$ are the standard $L_2(\mathscr{D})$ norm and $L_\infty(\mathscr{D})$ norm, respectively, defined on the domain $\mathscr{D}$, $\rho$ is the value computed using the present approach and $\rho_\mathrm{exact}$ is given by Eq.(\ref{eq:supersonic vortex problem - exact solution}).
The steady-state solution is assumed to be reached when the relative error between the evaluation of the norms given in Eq.(\ref{eq:supersonic vortex problem - error measures}) between two consecutive time steps is less than 10$^{-5}$.
For all considered simulations, this occurred at $t > 5.0$.

\begin{figure}
\centering
\begin{subfigure}{0.45\textwidth}
\includegraphics[width = \textwidth]{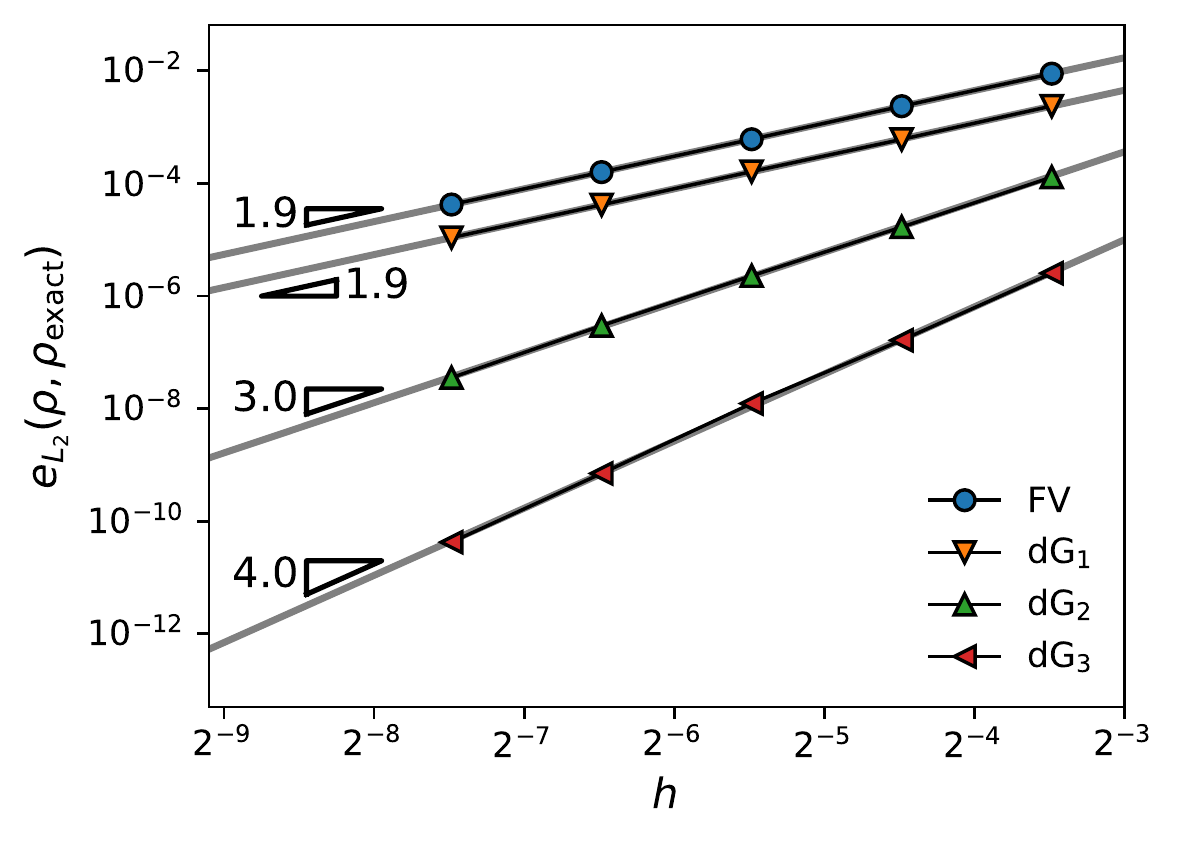}
\caption{}
\end{subfigure}
\
\begin{subfigure}{0.45\textwidth}
\includegraphics[width = \textwidth]{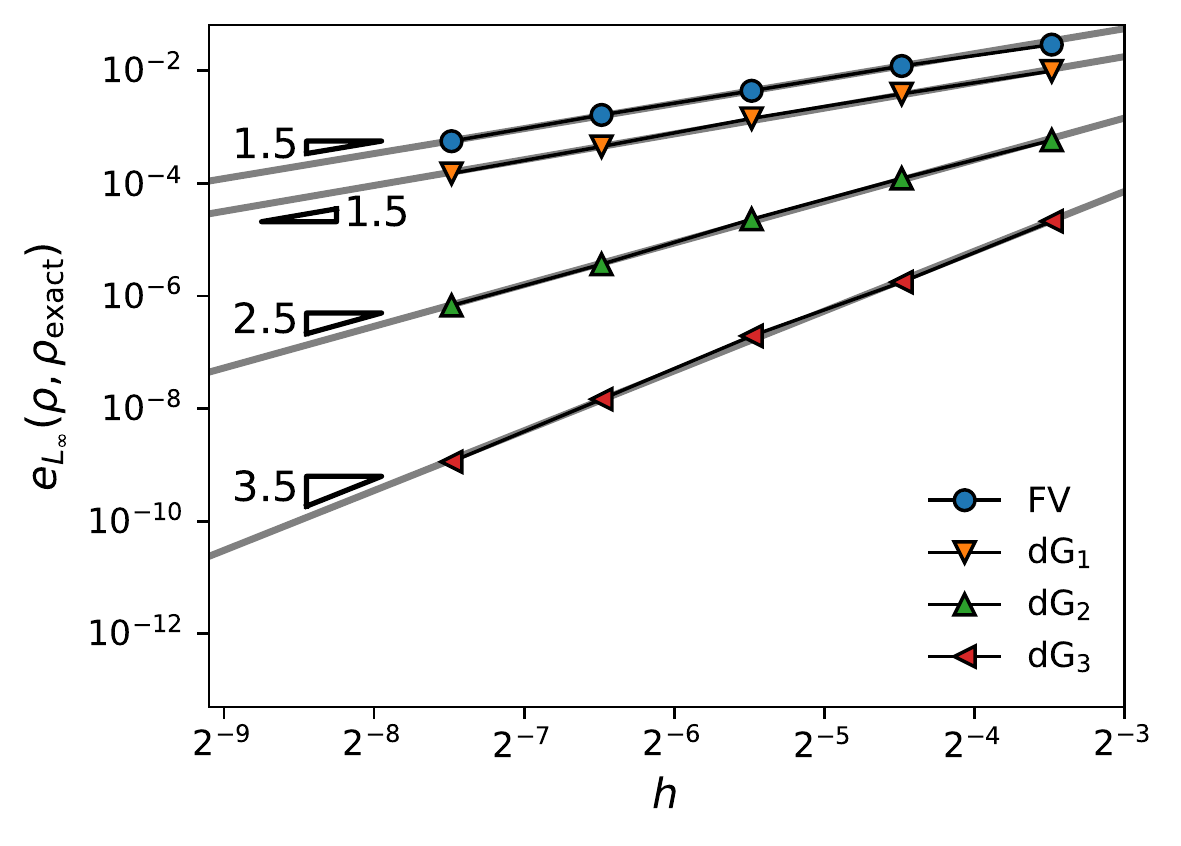}
\caption{}
\end{subfigure}
\caption
{
    Convergence in (a) the $L_2$ norm and (b) the $L_\infty$ norm of the density error for the supersonic vortex problem.
}
\label{fig:RESULTS - SUPERSONIC VORTEX - HP-CONVERGENCE}
\end{figure}

The simulations are performed using a single-level mesh, with $h \equiv 1.43/n$, where $n$ is the number of elements per side of the background square, and where the volume fraction threshold $\overline{\nu} = 0.3$.
Figure (\ref{fig:RESULTS - SUPERSONIC VORTEX - GEOMETRY-MESH-ERROR}b) shows the density distribution obtained using a $16\times 16$ background mesh and a dG$_3$ scheme; the corresponding density error distribution is reported in fig.(\ref{fig:RESULTS - SUPERSONIC VORTEX - GEOMETRY-MESH-ERROR}c), which also shows that the maximum error is observed near the inner boundary of the annulus.
The error measures defined in Eq.(\ref{eq:supersonic vortex problem - error measures}) are shown in Fig.(\ref{fig:RESULTS - SUPERSONIC VORTEX - HP-CONVERGENCE}) as functions of the numerical scheme and the mesh size $h$.
As seen in the figures, the dG schemes show an $\mathcal{O}(h^{p+1})$ convergence rate in the $L_2$ norm and an $\mathcal{O}(h^{p+1/2})$ convergence rate in the $L_\infty$ norm, while the FV scheme using the reconstruction (without limiting) presented in Sec.(\ref{ssec:RECONSTRUCTION}) obtains the same convergence rate as the dG$_1$ scheme, albeit with a slightly larger error.

\subsection{Embedded Sod's shock tube}\label{ssec:EMBEDDED SODS TUBE}
\begin{figure}
\centering
\includegraphics[width = \textwidth]{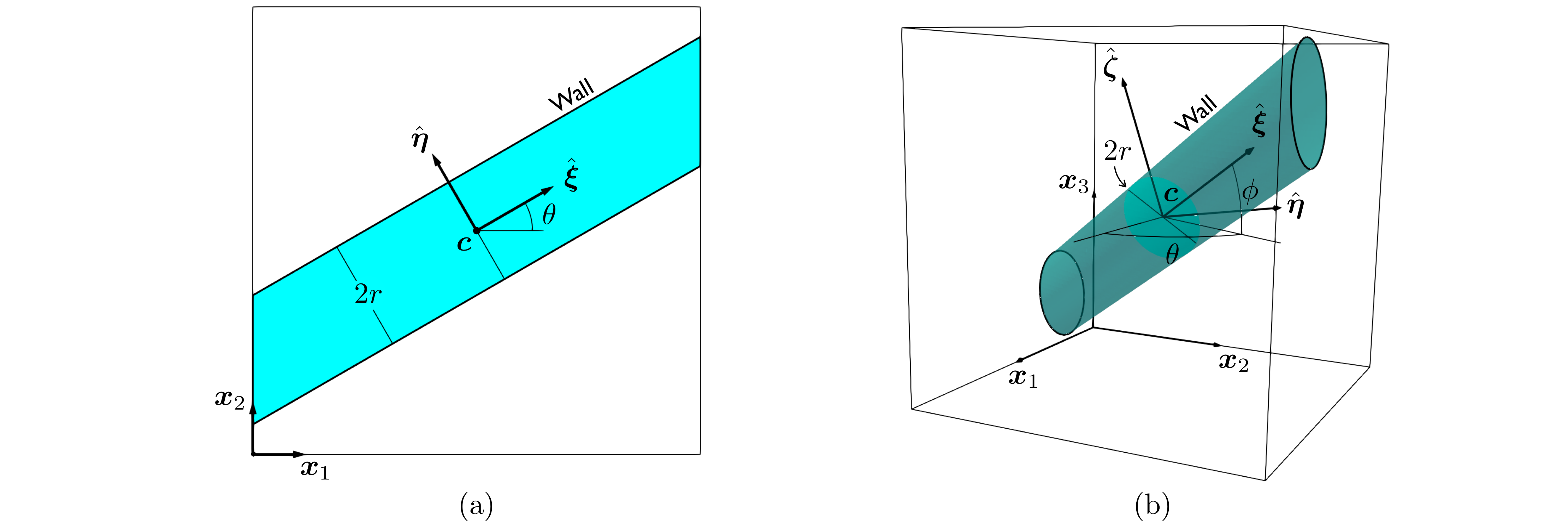}
\caption
{
    (a) Two-dimensional and (b) three-dimensional geometry for the tilted Sod's tube problem defined in the background square $[0,1]^2$ and in the background cube $[0,1]^3$, respectively.
}
\label{fig:RESULTS - SODS TUBE - GEOMETRY}
\end{figure}

In the second test, we consider Sod's shock tube problem in embedded geometries.
These simulations assess the ability of the implicitly-defined mesh and the reconstruction and limiting technique presented in Sec.(\ref{ssec:RECONSTRUCTION}) to reproduce a one-dimensional flow within a geometry that is not aligned with the background grid.

The problem is defined in the background unit cube $\mathscr{R} = [0,1]^d$ as shown in Fig.(\ref{fig:RESULTS - SODS TUBE - GEOMETRY}).
In 2D, the tube is a strip of width $2r$ centered at $\bm{c} \equiv \{0.5,0.5\}$ whose orientation is defined by the parameter $\theta$, which denotes the inclination of the tube's centerline.
In 3D, the tube is a cylinder with a circular cross section of radius $r$ centered at  $\bm{c} \equiv \{0.5,0.5,0.5\}$ whose orientation is defined by the two parameters $\theta$ and $\phi$, which denote the azimuth and the elevation, respectively, of the tube's centerline.
For convenience, we introduce a rotated coordinate system
$\hat{\bm{\xi}}\hat{\bm{\eta}}$ (in 2D) or $\hat{\bm{\xi}}\hat{\bm{\eta}}\hat{\bm{\zeta}}$ (in 3D) centered in $\bm{c}$, where $\hat{\bm{\xi}}$ is aligned with the tube's centerline and the unit vectors have the following components
\begin{equation}\label{eq:Sod's tube problem - local RS}
\{\hat{\bm{\xi}},\hat{\bm{\eta}}\} =
\left\{
\begin{array}{c}
\cos\theta \\
\sin\theta
\end{array},
\begin{array}{c}
-\sin\theta \\
\cos\theta
\end{array}
\right\}
\quad\mathrm{or}\quad
\{\hat{\bm{\xi}},\hat{\bm{\eta}},\hat{\bm{\zeta}}\} =
\left\{
\begin{array}{c}
\cos\theta\cos\phi \\
\sin\theta\cos\phi \\
\sin\phi
\end{array},
\begin{array}{c}
-\sin\theta \\
\cos\theta \\
0
\end{array},
\begin{array}{c}
-\cos\theta\sin\phi \\
-\sin\theta\sin\phi \\
\cos\phi
\end{array}
\right\}
\end{equation}
in 2D or 3D, respectively.
Then, a convenient way to represent the domain and its boundaries is through the level set function $\bm{\Phi}(\bm{x}) \equiv \eta^2-r^2$ in 2D, or $\bm{\Phi}(\bm{x}) \equiv \eta^2+\zeta^2-r^2$ in 3D, where $\xi$, $\eta$ and $\zeta$ are the components of $\bm{x}$ in the rotated coordinate system.

The initial conditions are given by
\begin{equation}
\bm{U}(t = 0,\bm{x}) =
\left\{
\begin{array}{cc}
\bm{U}_0^L&\mathrm{if}~(x_l-c_l)\hat{\xi}_l \le 0 \\
\bm{U}_0^R&\mathrm{if}~(x_l-c_l)\hat{\xi}_l > 0
\end{array}
\right.,
\quad\mathrm{where}\quad
\bm{U}_0^L =
\left\{
\begin{array}{c}
1\\
\bm{0}\\
\frac{1}{\gamma-1}
\end{array}
\right\}
\quad\mathrm{and}\quad
\bm{U}_0^R =
\left\{
\begin{array}{c}
0.125\\
\bm{0}\\
\frac{0.1}{\gamma-1}
\end{array}
\right\},
\end{equation}
and the final time is $T = 0.2$.
Wall boundary conditions are prescribed along the embedded boundaries.
At the intersection between the domain and the background cube's boundary, ghost states with the initial conditions are prescribed at the quadrature points.

\begin{figure}
\centering
\includegraphics[width = \textwidth]{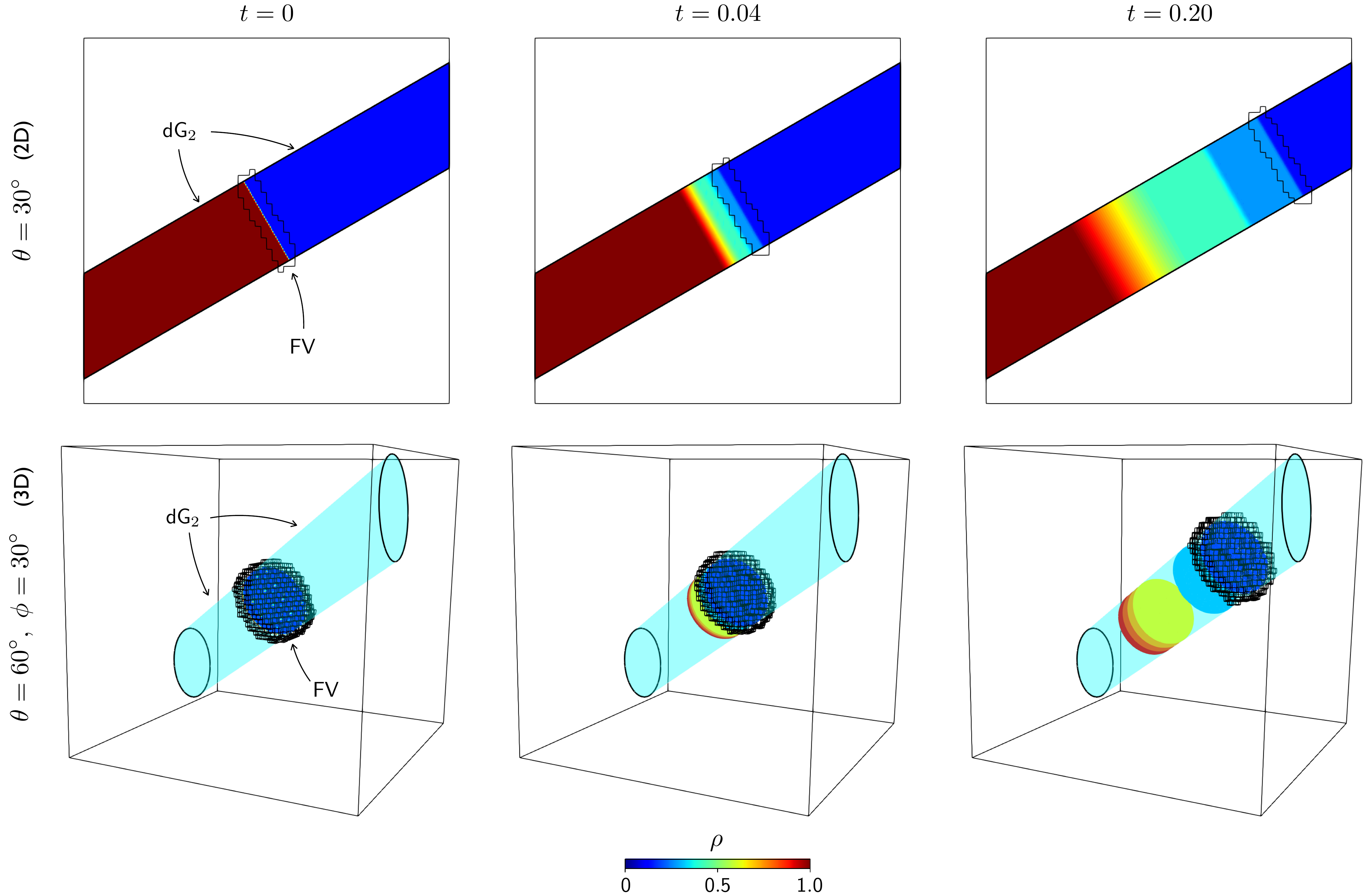}
\caption
{
    Density distribution at times $t = 0$, $t = 0.04$ and $t = 0.20$ for the 2D embedded Sod's shock tube problem (top row) and the 3D embedded shock tube problem (bottom row).
    The orientation angles $\theta = 30^\circ$, in 2D, and $\theta = 60^\circ$ and $\phi = 30^\circ$, in 3D, are those defined in Fig.(\ref{fig:RESULTS - SODS TUBE - GEOMETRY}).
    For the 3D case, the density distribution is represented by isosurfaces of $\rho = \{0.15, 0.3, 0.6, 0.7, 0.8, 0.9\}$.
}
\label{fig:RESULTS - SODS TUBE - DENSITY SCREENSHOTS}
\end{figure}

The solution to Sod's shock tube problem consists of three waves: a rarefaction wave, a contact discontinuity and a shock wave.
To account for the presence of the discontinuities, a two-level mesh is employed.
Level $\ell = 0$ is obtained by subdividing the background domain into $64^d$ cells and uses a dG$_2$ scheme.
Level $\ell = 1$ is constructed using a refinement ratio of 4 and is dynamically updated from level 0 using the shock sensor described in Sec.(\ref{ssec:CELL TAGGING}).
Specifically, at time $t$, the smoothness indicator $s^e$ given in Eq.(\ref{eq:smoothness indicator}) is computed for each element of the implicitly-defined mesh at level 0; then, the elements satisfying the condition (\ref{eq:shock sensor}) are tagged for refinement to level 1.
In Eq.(\ref{eq:shock sensor}), $\kappa_s = -2.00$ for the 2D problem and $\kappa_s = -2.45$ for the 3D problem.
Finally, the implicitly-defined mesh is generated using volume fraction thresholds $\overline{\nu}_0 = \overline{\nu}_1 = 0.3$ in 2D and $\overline{\nu}_0 = \overline{\nu}_1 = 0.15$ in 3D.
It is worth pointing out that the implicitly-defined elements are simple polygons in 2D, whereas they have curved surfaces in 3D.
However, as shown in the results, this will have very little effect on the computed numerical solutions.

The top row of images in Fig.(\ref{fig:RESULTS - SODS TUBE - DENSITY SCREENSHOTS}) shows the density distribution at times $t = 0$, $t = 0.04$ and $t = 0.20$ for the 2D tube corresponding to the orientation $\theta = 30^\circ$.
The bottom row of images in Fig.(\ref{fig:RESULTS - SODS TUBE - DENSITY SCREENSHOTS}) shows  the density distribution at the same times for the 3D tube corresponding to the orientation $\theta = 60^\circ$ and $\phi = 30^\circ$.
In 3D, the density distribution is represented by isosurfaces of $\rho = \{0.15, 0.3, 0.6, 0.7, 0.8, 0.9\}$, which allow us to see the rarefaction wave, the contact discontinuity and the shock along the tube.
Figure (\ref{fig:RESULTS - SODS TUBE - DENSITY SCREENSHOTS}) also displays the evolution of the mesh level using the FV scheme, whose location is identified by the stepped black lines.
The figure shows that the FV level tracks the shock discontinuity but not the contact discontinuity.
Although it is possible to adjust the constant $\kappa_s$ in (\ref{eq:shock sensor}) in such a way as to make the FV level track the contact discontinuities, numerical experiments show that these discontinuities are only temporarily tracked by the FV scheme and, once they are sufficiently smoothed, they stop satisfying the condition (\ref{eq:shock sensor}) and start being tracked by the unlimited dG scheme.
This is an expected behavior since shocks steepen while contact discontinuities do not.

\begin{figure}
\centering
\includegraphics[width = \textwidth]{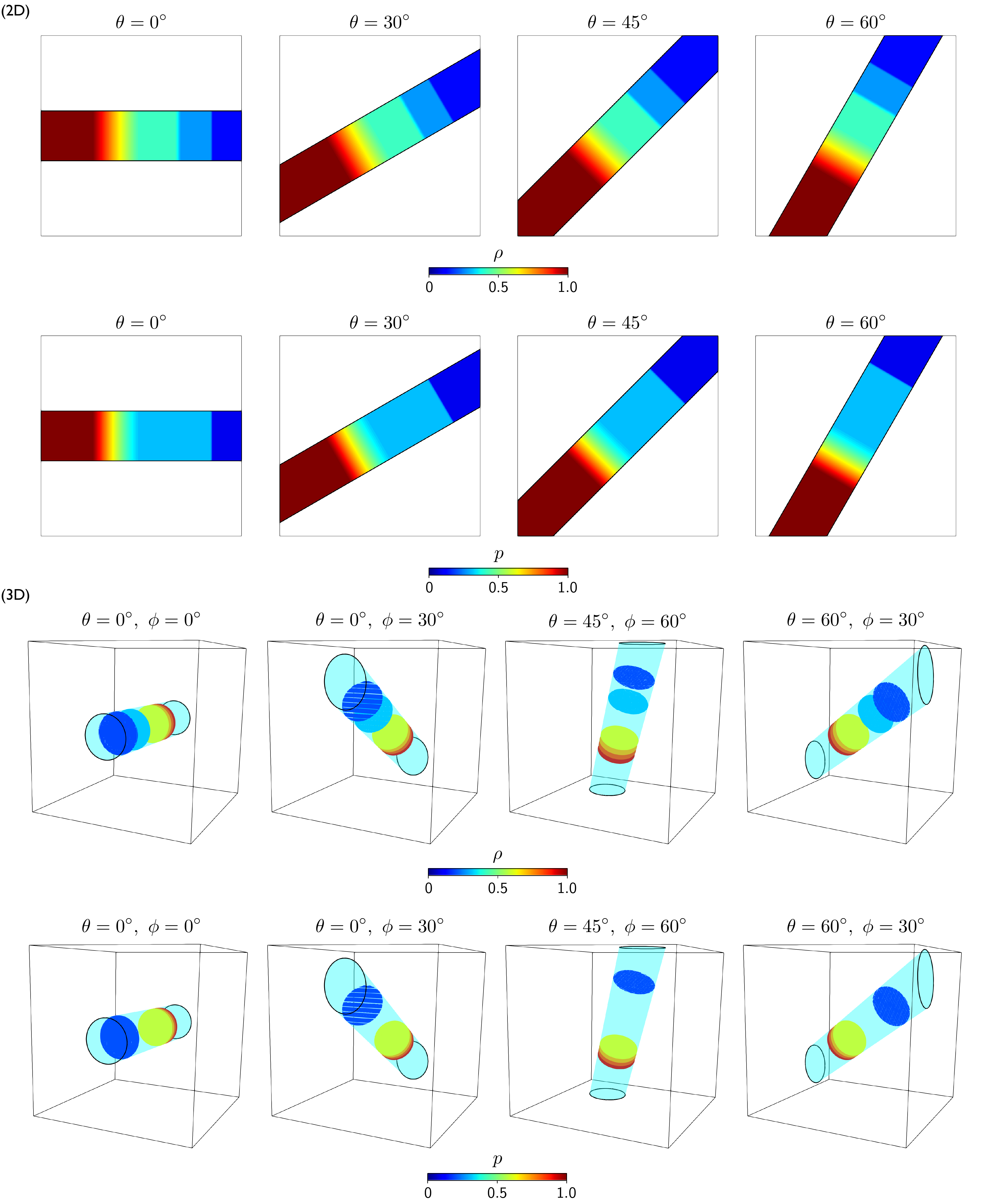}
\caption
{
    Density and pressure distribution for the eight configurations of the embedded shock tube problem at the final time $t = 0.2$.
    From the top: (first row) density distribution in 2D, (second row) pressure distribution in 2D, (third row) density distribution in 3D and (fourth row) pressure distribution in 3D.
    For the 3D cases, the density distribution is represented by isosurfaces of $\rho = \{0.15, 0.3, 0.6, 0.7, 0.8, 0.9\}$ and the pressure distribution is represented by isosurfaces of $p = \{0.15, 0.6, 0.7, 0.8, 0.9\}$.
}
\label{fig:RESULTS - SODS TUBE 2D - CONTOURS - DENSITY AND PRESSURE}
\end{figure}

\begin{figure}
\centering
\begin{subfigure}{0.235\textwidth}
\includegraphics[width=\textwidth]{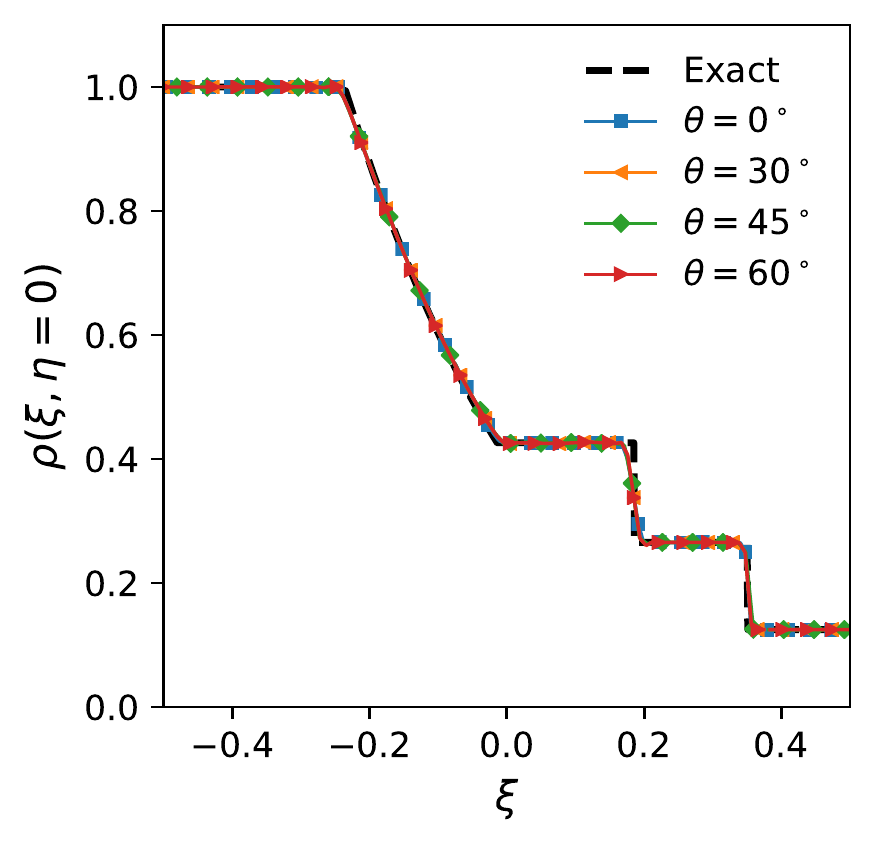}
\caption{}
\end{subfigure}
\
\begin{subfigure}{0.235\textwidth}
\includegraphics[width=\textwidth]{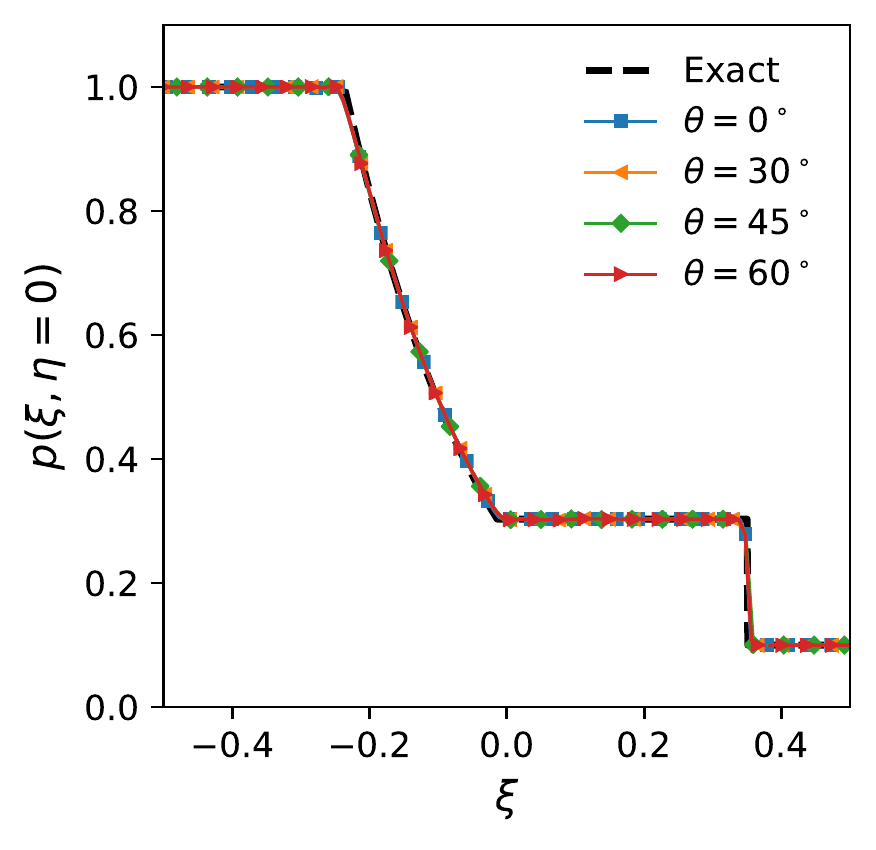}
\caption{}
\end{subfigure}
\
\begin{subfigure}{0.235\textwidth}
\includegraphics[width=\textwidth]{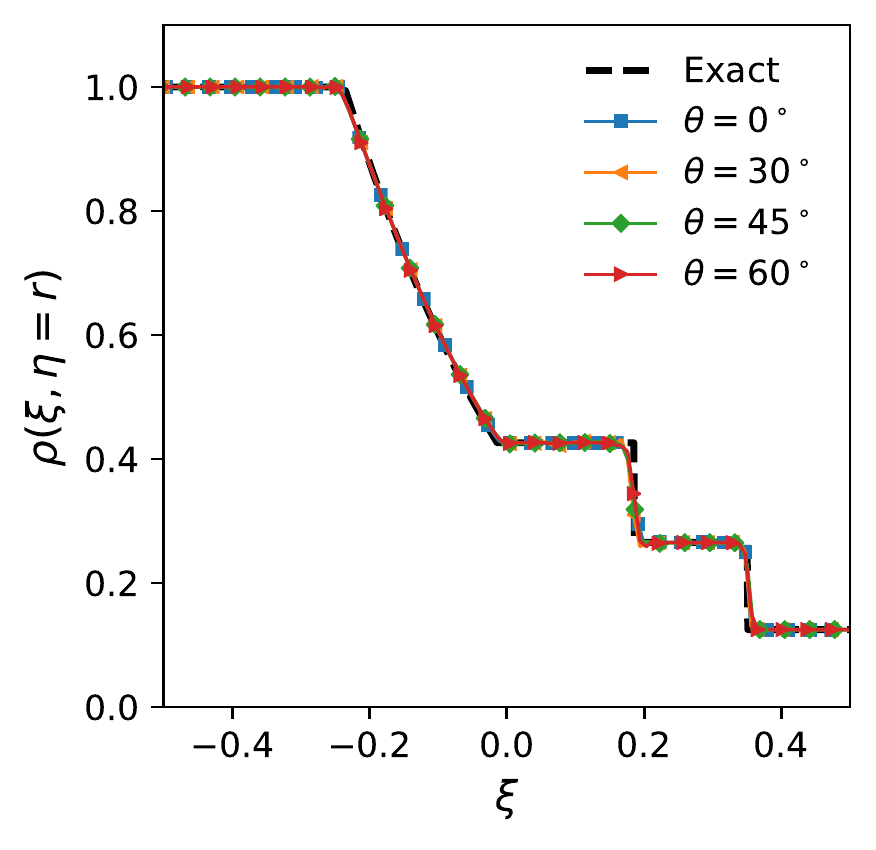}
\caption{}
\end{subfigure}
\
\begin{subfigure}{0.235\textwidth}
\includegraphics[width=\textwidth]{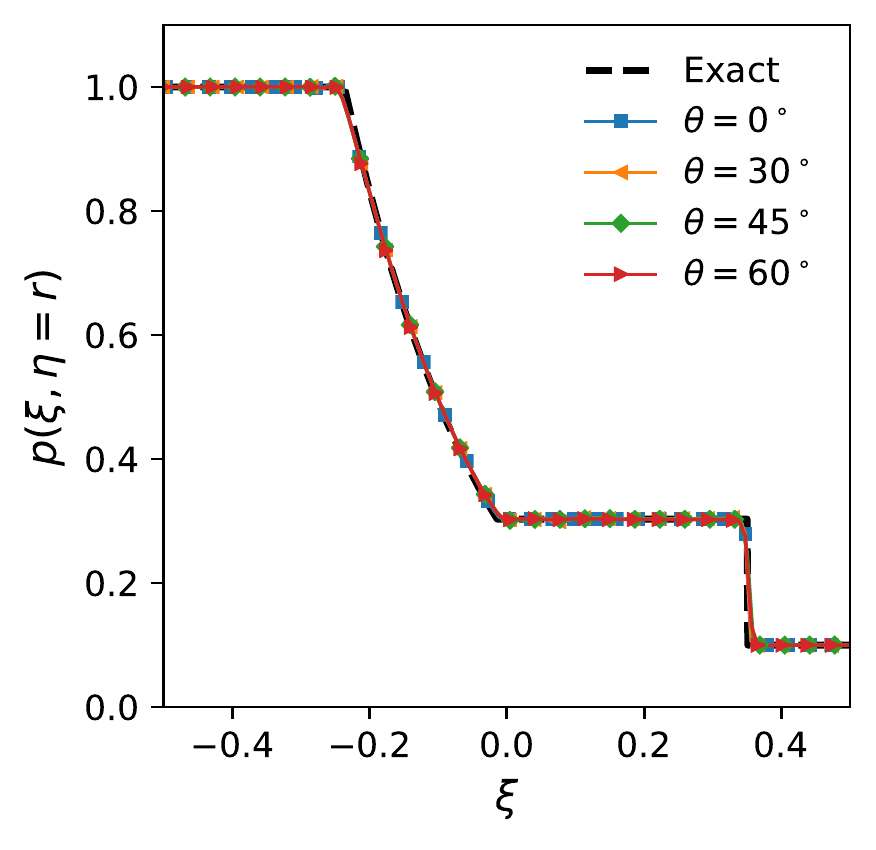}
\caption{}
\end{subfigure}
\caption
{
    Distribution of density and pressure at the final time $t = 0.2$ as functions of the orientation of the 2D Sod's shock tube;
    (a) density along the centerline, (b) pressure along the centerline, (c) density along the wall and (d) pressure along the wall.
}
\label{fig:RESULTS - SODS TUBE 2D - LINE PLOTS - DENSITY AND PRESSURE}
\end{figure}

\begin{figure}
\centering
\begin{subfigure}{0.235\textwidth}
\includegraphics[width=\textwidth]{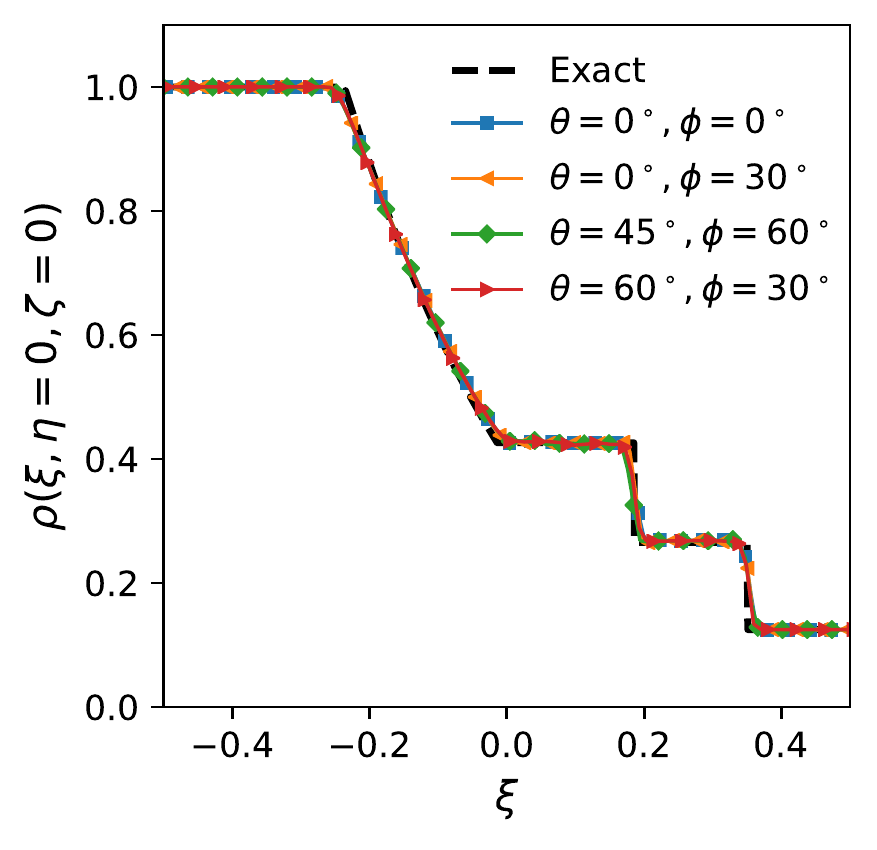}
\caption{}
\end{subfigure}
\
\begin{subfigure}{0.235\textwidth}
\includegraphics[width=\textwidth]{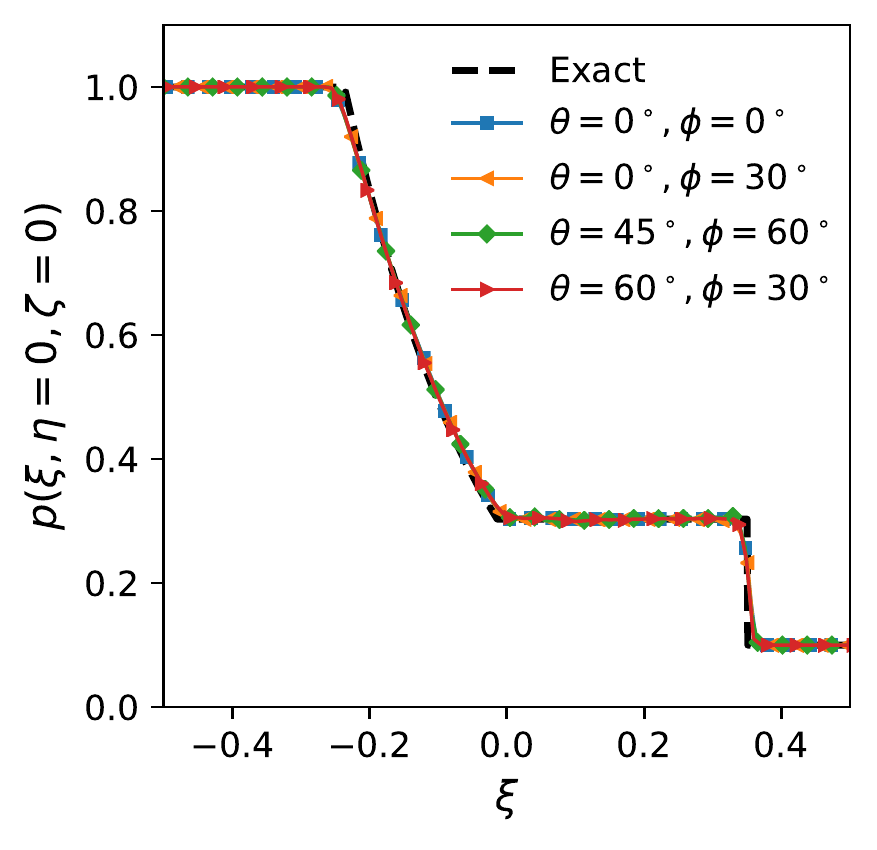}
\caption{}
\end{subfigure}
\
\begin{subfigure}{0.235\textwidth}
\includegraphics[width=\textwidth]{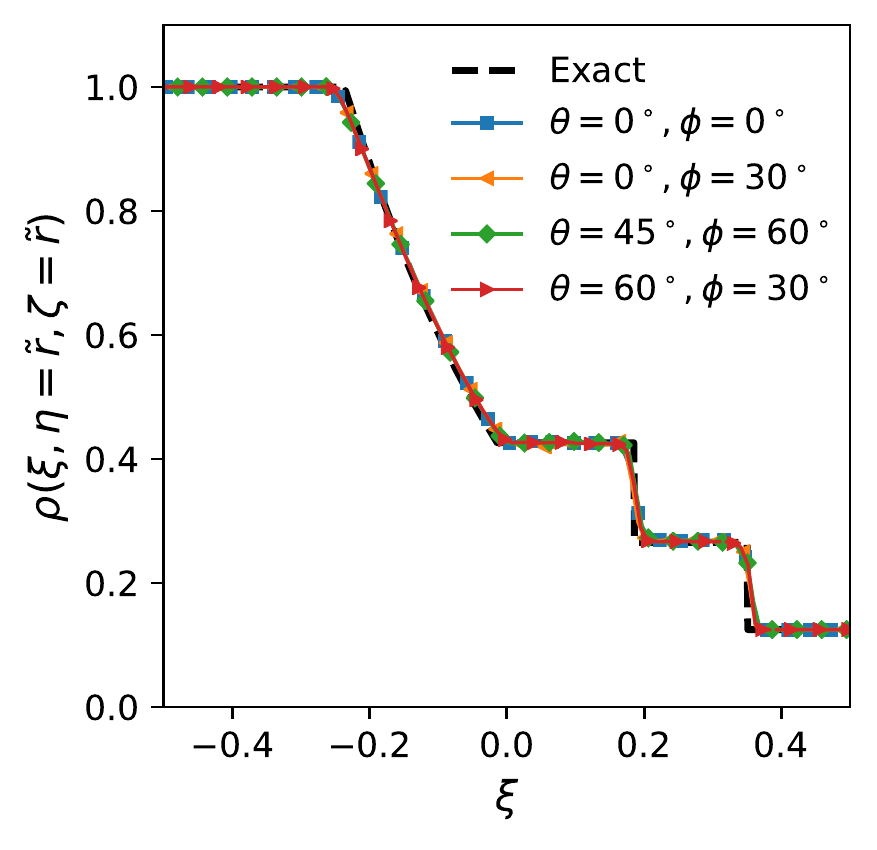}
\caption{}
\end{subfigure}
\
\begin{subfigure}{0.235\textwidth}
\includegraphics[width=\textwidth]{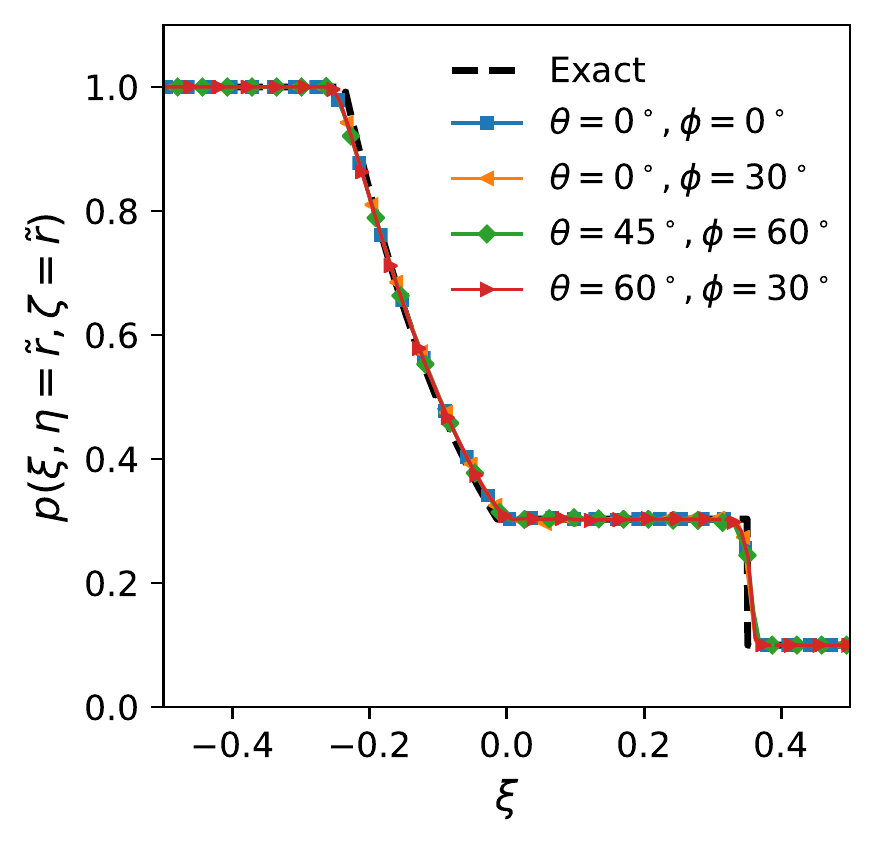}
\caption{}
\end{subfigure}
\caption
{
    Distribution of density and pressure at the final time $t = 0.2$ as functions of the orientation of the 3D Sod's shock tube;
    (a) density along the centerline, (b) pressure along the centerline, (c) density along the wall and (d) pressure along the wall.
}
\label{fig:RESULTS - SODS TUBE 3D - LINE PLOTS - DENSITY AND PRESSURE}
\end{figure}

The robustness of the proposed approach is tested by considering four orientations each in 2D and 3D while keeping the same settings for the implicitly-defined mesh and the shock sensor.
Figure (\ref{fig:RESULTS - SODS TUBE 2D - CONTOURS - DENSITY AND PRESSURE}) shows the results at the final time $t = 0.2$ for orientation angles $\theta = 0^\circ$, $30^\circ$, $45^\circ$ and $60^\circ$, in 2D, and for orientation angles $\{\theta,\phi\} = \{0^\circ,0^\circ\}$, $\{0^\circ,30^\circ\}$, $\{45^\circ,60^\circ\}$ and $\{60^\circ,30^\circ\}$, in 3D.
The first row shows the density distribution for the four 2D configurations, the second row shows the corresponding pressure distribution, the third row shows the density distribution for the four 3D configurations and the fourth row shows the corresponding pressure.
For the 3D cases, the density distribution is represented by isosurfaces of $\rho = \{0.15, 0.3, 0.6, 0.7, 0.8, 0.9\}$ while the pressure distribution is represented by isosurfaces of $p = \{0.15, 0.6, 0.7, 0.8, 0.9\}$.
For all configurations, the coupled dG-FV scheme is able to reproduce the one-dimensional flow in embedded curved geometries that are not aligned with the background grid.
This can be seen more clearly by looking at line plots of the density and pressure distribution along the centerline and the wall for the embedded shock tubes at the final time $t = 0.2$.
Figures (\ref{fig:RESULTS - SODS TUBE 2D - LINE PLOTS - DENSITY AND PRESSURE}a), (\ref{fig:RESULTS - SODS TUBE 2D - LINE PLOTS - DENSITY AND PRESSURE}b), (\ref{fig:RESULTS - SODS TUBE 2D - LINE PLOTS - DENSITY AND PRESSURE}c) and (\ref{fig:RESULTS - SODS TUBE 2D - LINE PLOTS - DENSITY AND PRESSURE}d) show the density along the centerline, the pressure along the centerline, the density along the wall and the pressure along the wall, respectively, as functions of the local reference system coordinate $\xi$ and the 2D tube's orientation.
Similarly, Figs.(\ref{fig:RESULTS - SODS TUBE 3D - LINE PLOTS - DENSITY AND PRESSURE}a), (\ref{fig:RESULTS - SODS TUBE 3D - LINE PLOTS - DENSITY AND PRESSURE}b), (\ref{fig:RESULTS - SODS TUBE 3D - LINE PLOTS - DENSITY AND PRESSURE}c) and (\ref{fig:RESULTS - SODS TUBE 3D - LINE PLOTS - DENSITY AND PRESSURE}d) show the same quantities as functions of the local reference system coordinate $\xi$ and\ the 3D tube's orientation.
As seen in the figures, the rarefaction wave, the contact discontinuity and the shock wave are all well captured by the proposed approach and match the one-dimensional solution well, regardless of the orientation of the embedded geometry.
Moreover, the numerical results in the tilted geometries are nearly indistinguishable from the numerical results in the aligned geometry, which can be viewed as a reference solution where the implicitly-defined mesh has little effect.

Finally, we note in Fig.(\ref{fig:RESULTS - SODS TUBE - GEOMETRY}), the wall values in 2D are those computed at $\eta = r$, i.e.~at the upper embedded boundary of the tube, whereas the wall values in 3D are those computed at $\eta = \zeta = \tilde{r} \equiv r/\sqrt{2}$, i.e.~at one line on the embedded boundary of the tube.
Similar results are observed if the wall values are evaluated at other boundary locations.

\subsection{Shock reflection}\label{ssec:SHOCK REFLECTION}
\begin{figure}
\centering
\includegraphics[width = \textwidth]{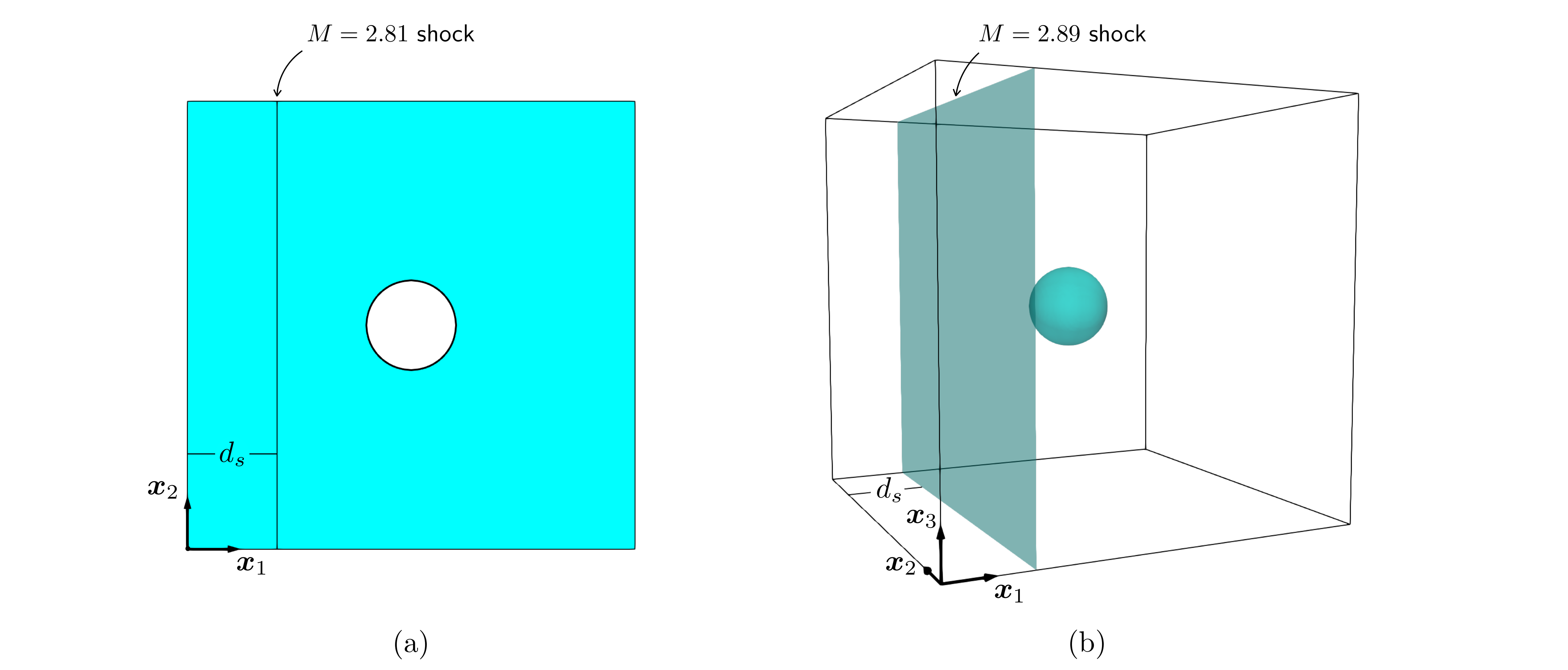}
\caption
{
    Illustration of (a) the problem of a shock reflecting from a cylinder and (b) the problem of a shock reflecting from a sphere.
}
\label{fig:RESULTS - SHOCK REFLECTION - GEOMETRY}
\end{figure}

In the third test, we consider a shock reflection from a cylinder in 2D and a shock reflection from a sphere in 3D.
These problems have been investigated experimentally \cite{bryson1961diffraction} and allow us to assess the capabilities of numerical methods to capture the shock structure as well as the details of the smooth flow patterns.

Shock reflection from a cylinder in 2D is defined in the background unit square $\mathscr{R} = [0,1]^2$.
The geometry is represented by the level set function $\Phi(\bm{x}) \equiv r^2-(x_k-c_k)(x_k-c_k)$ where $\bm{c} =\{0.5,0.5\}$ and $r = 0.1$ are the center and the radius of the cylinder, respectively as shown in Fig.(\ref{fig:RESULTS - SHOCK REFLECTION - GEOMETRY}a).
The shock is initially located at $x_1 = d_s = 0.2$ and travels towards the cylinder at a mach number $M = 2.81$.
The initial conditions are given by
\begin{equation}\label{eq:shock reflection problem - ICs}
\bm{U}(t = 0,\bm{x}) =
\left\{
\begin{array}{cc}
\bm{U}_0^L&\mathrm{if}~x_1 \le d_s \\
\bm{U}_0^R&\mathrm{if}~x_1 > d_s
\end{array}
\right.,
\quad\mathrm{where}\quad
\bm{U}_0^R =
\left\{
\begin{array}{c}
1\\
\bm{0}\\
\frac{1}{\gamma-1}
\end{array}
\right\}
\end{equation}
where $\bm{U}_0^L$ is derived from the shock jump relations \cite{toro2013riemann}.
Wall boundary conditions are prescribed on the boundary of the cylinder and at $x_2 = 0$ and $x_2 = 1$; inflow boundary conditions are prescribed at $x_1 = 0$ and outflow boundary conditions are prescribed at $x_2 = 1$.
The final time of the simulation is $T = 0.2$.

We use a four-level mesh with a refinement ratio of four at each level.
Level $\ell = 0$ uses a dG$_1$ scheme and a mesh size $h_0 = 1/64$;
level $\ell = 1$ uses a dG$_3$ scheme and is generated from level 0 using $\kappa_0^\rho = 0.5$ in (\ref{eq:density gradient refinement});
level $\ell = 2$ uses a dG$_3$ scheme and is generated from level 1 using $\kappa_1^\rho = 2.0$ in (\ref{eq:density gradient refinement});
level $\ell = 3$ uses the FV scheme and is generated from level 2 using $\kappa^s = -3.25$ in (\ref{eq:shock sensor}).
Finally, cell merging is triggered by the volume fraction thresholds $\overline{\nu}_0 = \overline{\nu}_1 = \overline{\nu}_2 = \overline{\nu}_3 = 0.3$.

\begin{figure}
\centering
\includegraphics[width = \textwidth]{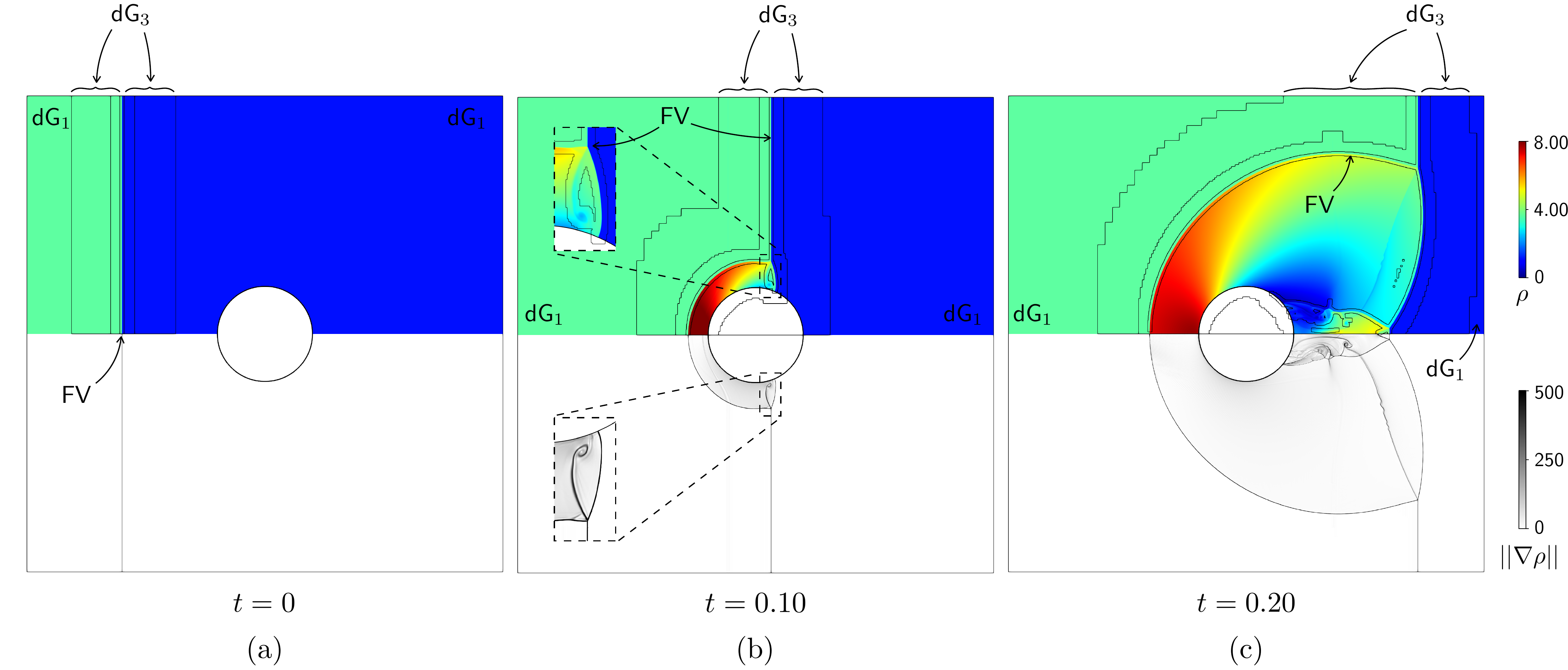}
\caption
{
    Distribution of the density $\rho$ and the density gradient magnitude $||\nabla\rho||$ at times (a) $t = 0$, (b) $t = 0.1$ and (c) $t = 0.2$ for a shock reflecting from a cylinder.
    The inset in figure (b) shows a close-up of the contact discontinuity and the Mach stem.
    In the images, the stepped black lines denote the outline of the AMR levels.
}
\label{fig:RESULTS - SHOCK REFLECTION - DENSITY - 2D}
\end{figure}

\begin{figure}
\centering
\includegraphics[width = \textwidth]{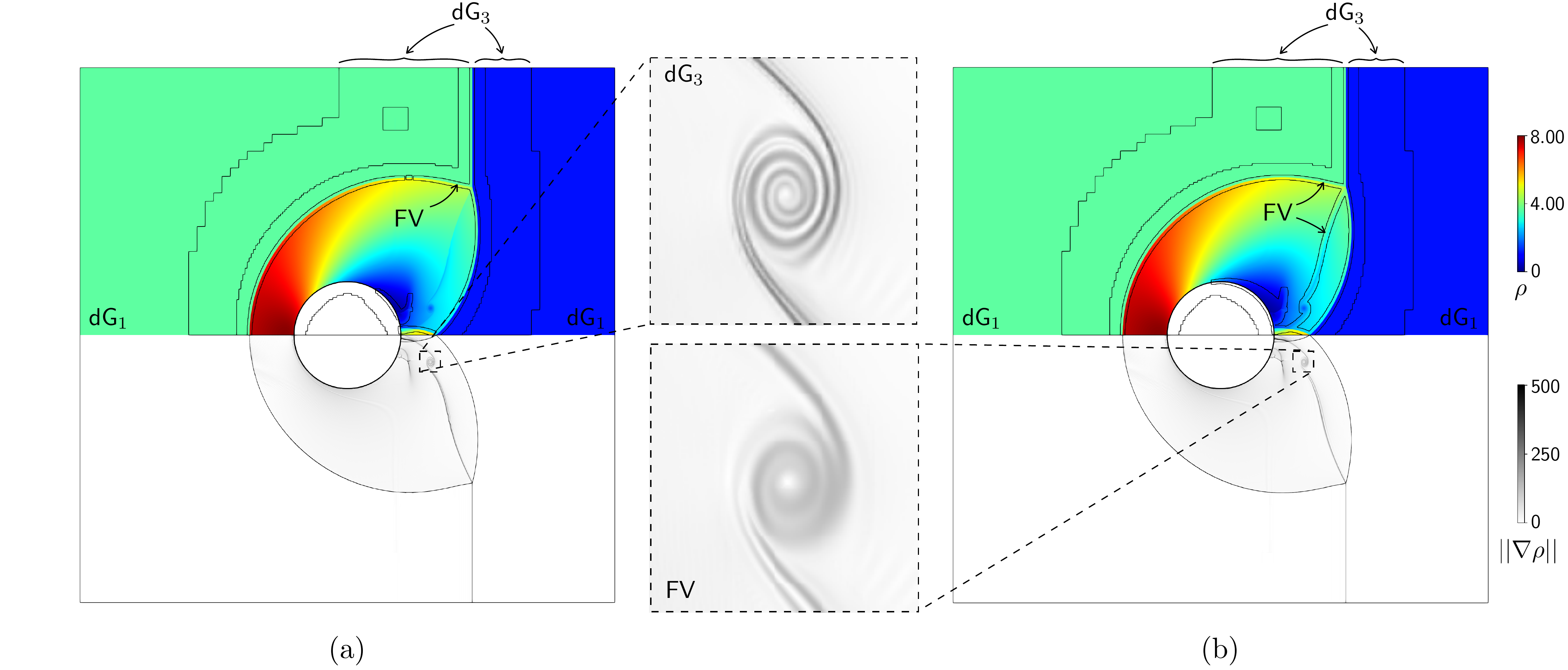}
\caption
{
    Distribution of the density $\rho$ and the density gradient magnitude $||\nabla\rho||$ computed at time $t = 0.16$ using different values of the threshold parameter $\kappa^s$ entering the shock-sensor condition (\ref{eq:shock sensor}):
    (a) $\kappa^s = -3.25$ and
    (b) $\kappa^s = -4.00$.
    The insets show a close-up on the vortex of the contact discontinuity.
}
\label{fig:RESULTS - SHOCK REFLECTION - DENSITY - DG VS FV}
\end{figure}

Figure (\ref{fig:RESULTS - SHOCK REFLECTION - DENSITY - 2D}) shows contour plots of the density $\rho$ and the density gradient magnitude $||\nabla\rho||$ at $t = 0$, $t = 0.1$ and $t = 0.2$.
The figure also show the configurations of the AMR levels and how they adapt to the location of the travelling discontinuities.
By looking at Fig.(\ref{fig:RESULTS - SHOCK REFLECTION - DENSITY - 2D}b), one can see how the FV level tracks the shock reflecting from the cylinder as well as the shock attached to the surface.
Similarly, once the incident shock has moved past the cylinder (Fig.(\ref{fig:RESULTS - SHOCK REFLECTION - DENSITY - 2D}c)), the FV level remains attached to the cylinder's surface to track the shock in the wake.

Figure (\ref{fig:RESULTS - SHOCK REFLECTION - DENSITY - 2D}) also illustrates the difference in how shocks and contact discontinuities are resolved during the simulation. 
In particular, by comparing Fig.(\ref{fig:RESULTS - SHOCK REFLECTION - DENSITY - 2D}b) and Fig.(\ref{fig:RESULTS - SHOCK REFLECTION - DENSITY - 2D}c), one can see that shocks are always resolved by the FV scheme, whereas the contact discontinuity is initially tracked by the FV scheme but, once smoothed, it starts being tracked by the dG schemes.
This is consistent to the results of the embedded shock tube test presented in Sec.(\ref{ssec:EMBEDDED SODS TUBE}).
To illustrate the benefits of using high-order schemes versus low-order schemes, Fig.(\ref{fig:RESULTS - SHOCK REFLECTION - DENSITY - DG VS FV}a) shows the density and the density gradient magnitude at time $t = 0.16$ for the simulation described above, whereas Fig.(\ref{fig:RESULTS - SHOCK REFLECTION - DENSITY - DG VS FV}b) shows the density and the density gradient magnitude at time $t = 0.16$ but with $\kappa^s = -4.00$.
A lower value of $\kappa^s$ increase the number elements where a discontinuity is sensed and, for this test case, causes the contact discontinuity to be tracked by the FV scheme.
By comparing the insets of Fig.(\ref{fig:RESULTS - SHOCK REFLECTION - DENSITY - DG VS FV}a) and Fig.(\ref{fig:RESULTS - SHOCK REFLECTION - DENSITY - DG VS FV}b), it is clear that the dG$_3$ scheme provides better resolution of the rollup of the contact discontinuity than the solution obtained with the FV scheme.

\begin{figure}
\centering
\includegraphics[width = \textwidth]{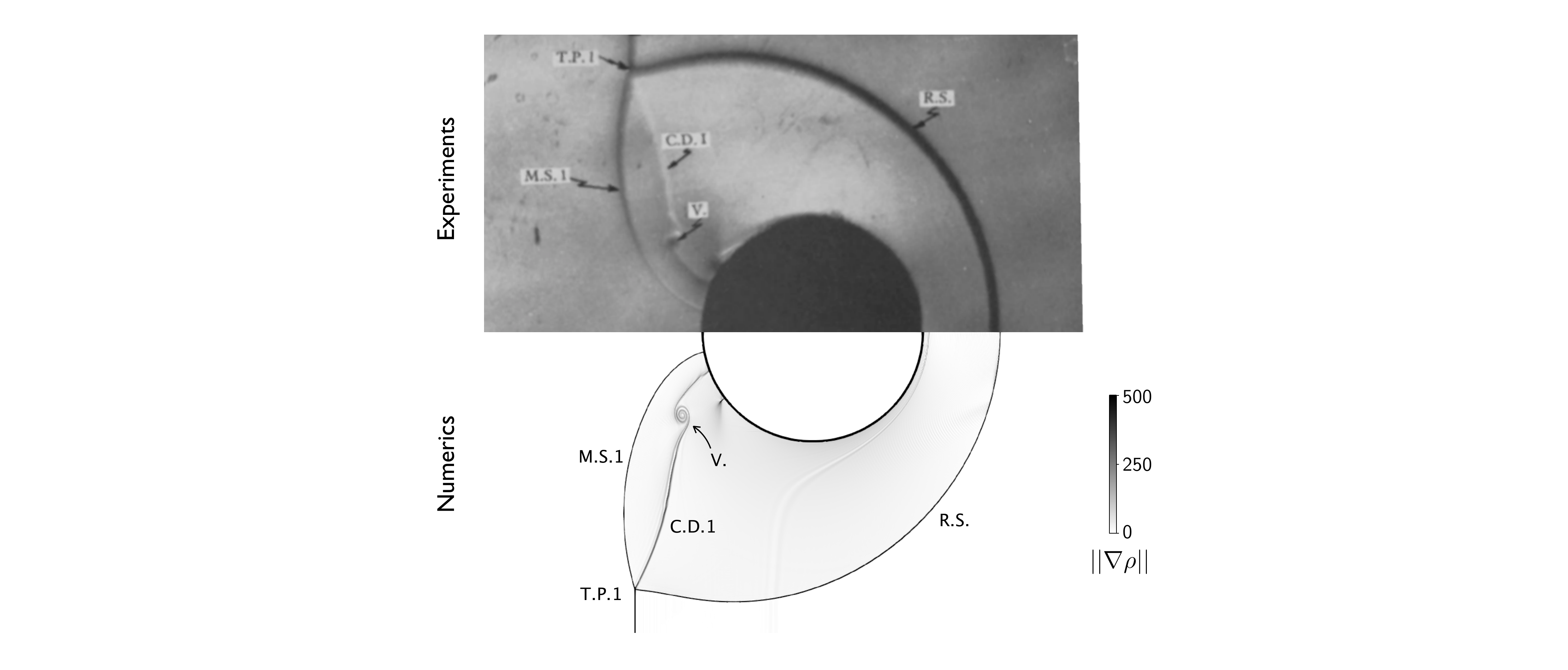}
\caption
{
    Comparison between the experimental Schlieren photograph of Ref.\cite{bryson1961diffraction} and the numerical solution at $t = 0.14$.
    All flow features observed experimentally, such as the reflected shock (R.S.), the primary Mach stem (M.S.1), the primary contact discontinuity (C.D.1), the primary triple point (T.P.1) and the vortex (V.), are well captured.
}
\label{fig:RESULTS - SHOCK REFLECTION - COMPARISON WITH EXPERIMENTS - 2D}
\end{figure}

Finally, Fig.(\ref{fig:RESULTS - SHOCK REFLECTION - COMPARISON WITH EXPERIMENTS - 2D}) shows the comparison between the Schlieren photograph from Ref.\cite{bryson1961diffraction} of a shock reflection from a cylinder and the results computed at $t = 0.14$.
As seen in the figure, the $hp$-AMR strategy allows us to capture all the details of the flow structures observed in the experiments.

\begin{figure}
\centering
\includegraphics[width = \textwidth]{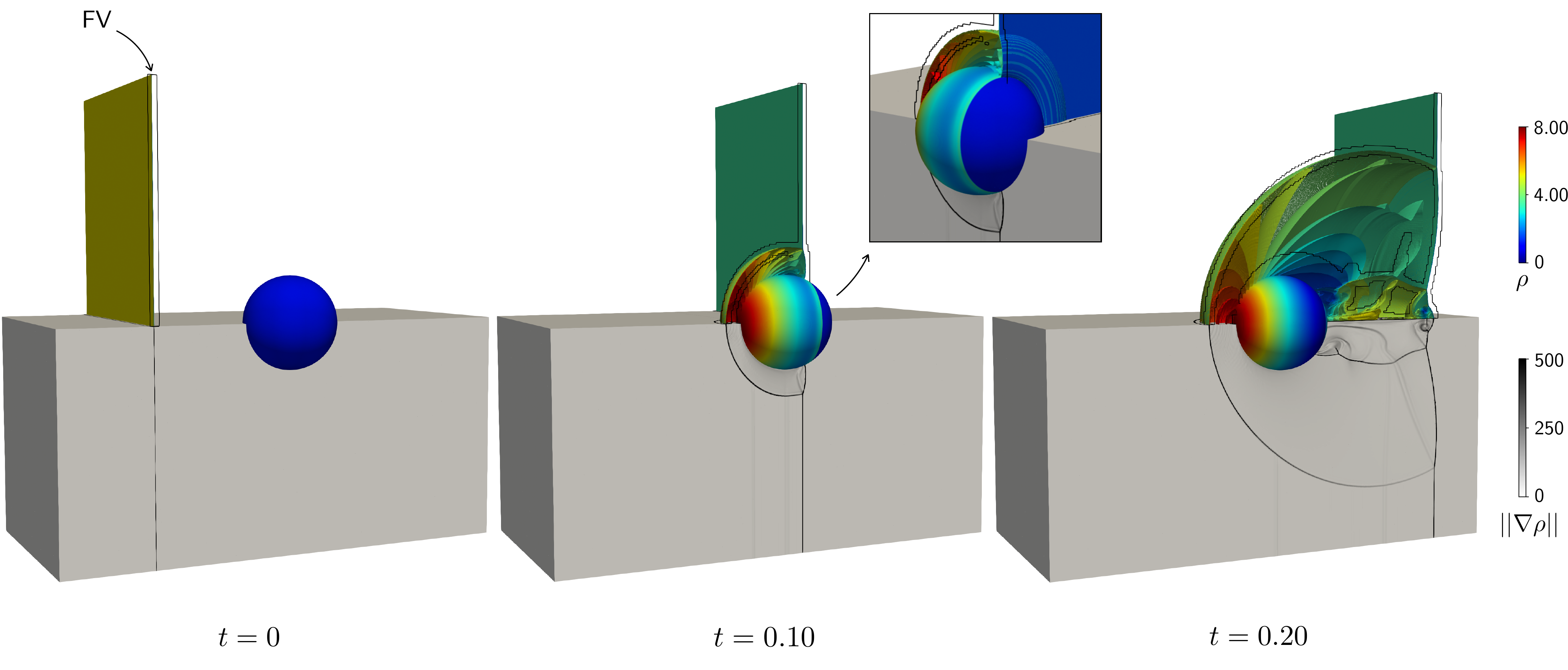}
\caption
{
    Distribution of the density $\rho$ and the density gradient magnitude $||\nabla\rho||$ at the times $t = 0$, $t = 0.1$ and $t = 0.2$ for a shock reflecting from a sphere.
    The inset shows a closeup on the shock propagating on the sphere surface at $t = 0.1$.
}
\label{fig:RESULTS - SHOCK REFLECTION - DENSITY - 3D}
\end{figure}

The case of a shock reflection from a sphere is modelled in 3D as shown in Fig.(\ref{fig:RESULTS - SHOCK REFLECTION - GEOMETRY}b).
The problem is defined in the background unit cube $\mathscr{R} = [0,1]^3$.
The geometry is represented by the level set function $\Phi(\bm{x}) \equiv r^2-(x_k-c_k)(x_k-c_k)$ where $\bm{c} =\{0.5,0.5,0.5\}$ and $r = 0.1$ are the center and the radius of the sphere, respectively.
The shock is initially located at $x_1 = d_s = 0.2$ and travels towards the sphere at a mach number $M = 2.89$.
The initial conditions are given by Eq.(\ref{eq:shock reflection problem - ICs}).
Wall boundary conditions are prescribed on the boundary of the sphere and at $x_2 = 0$, $x_2 = 1$, $x_3 = 0$ and $x_3 = 1$; inflow boundary conditions are prescribed at $x_1 = 0$ and outflow boundary conditions are prescribed at $x_2 = 1$.
The final simulation time is $T = 0.2$.

We use a four-level mesh with a refinement ratio of four between the levels.
Level $\ell = 0$ uses a dG$_1$ scheme and and a mesh size $h_0 = 1/8$;
level $\ell = 1$ uses a dG$_3$ scheme and is generated from level 0 using $\kappa_0^\rho = 0.5$ in (\ref{eq:density gradient refinement});
level $\ell = 2$ uses a dG$_3$ scheme and is generated from level 1 using $\kappa_1^\rho = 2.0$ in (\ref{eq:density gradient refinement});
level $\ell = 3$ uses the FV scheme and is generated from level 2 using $\kappa^s = -3.5$ in (\ref{eq:shock sensor}).
Cell merging is triggered by the volume fraction thresholds $\overline{\nu}_0 = \overline{\nu}_1 = \overline{\nu}_2 = \overline{\nu}_3 = 0.15$.

Figure (\ref{fig:RESULTS - SHOCK REFLECTION - DENSITY - 3D}) shows the density $\rho$ and the density gradient magnitude $||\nabla\rho||$ at $t = 0$, $t = 0.1$ and $t = 0.2$.
In the figure, the density on the surface of the sphere is represented by a continuous contour plot and, within the fluid domain, by isosurfaces evaluated at $\rho = \{1.5, 2.0, \dots, 7.5\}$.
The figure also shows the trace of the FV level on the planes $x_2 = 0.5$ and $x_3 = 0.5$; for the sake of visualization the AMR levels $\ell = 0$, $1$ and $2$ are not displayed.
Similar to 2D, one can see how the FV level tracks the reflected shock.

\begin{figure}
\centering
\includegraphics[width = \textwidth]{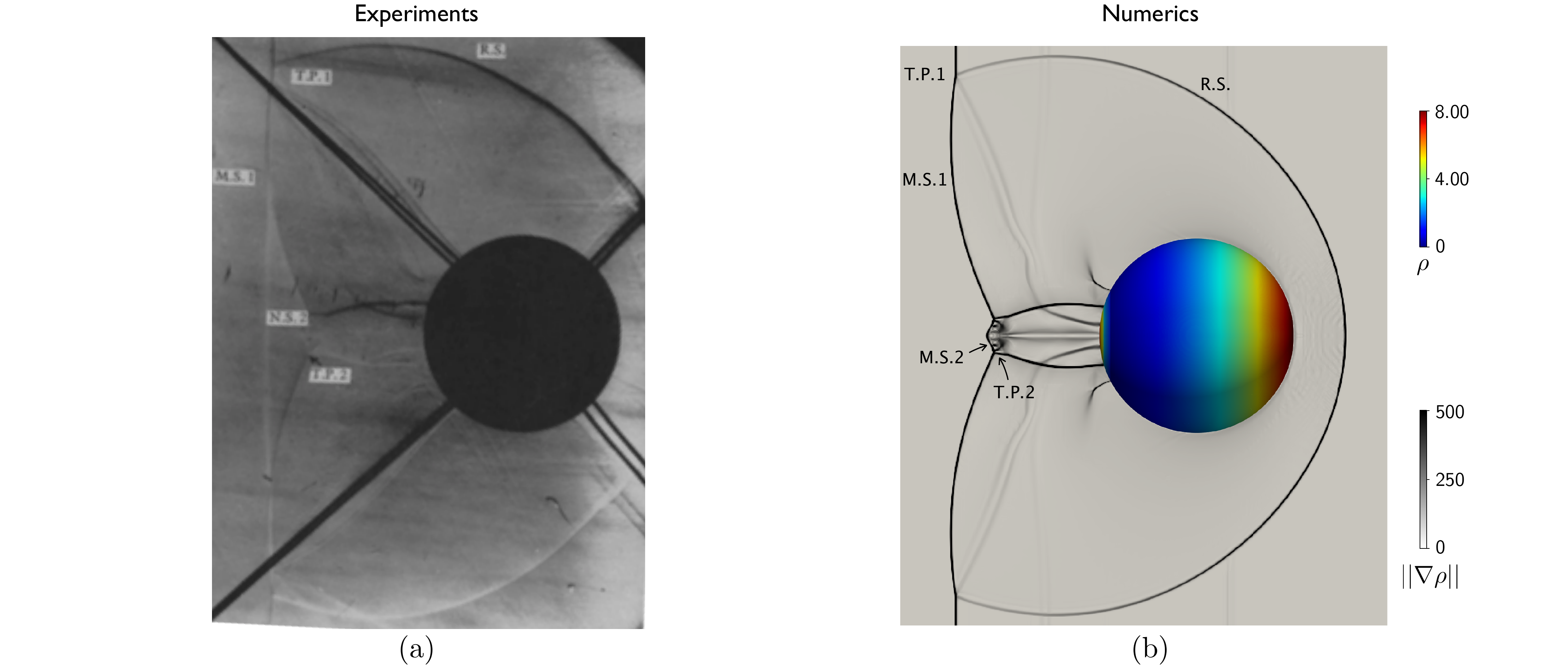}
\caption
{
    Comparison between the experimental Schlieren photograph from Ref.\cite{bryson1961diffraction} and the numerical solution at $t = 0.16$.
    The flow features observed experimentally, such as the reflected shock (R.S.), the primary Mach stem (M.S.1), the secondary Mach stem (M.S.2), the primary triple point (T.P.1) and the secondary triple point (T.P.2), are well captured.
}
\label{fig:RESULTS - SHOCK REFLECTION - COMPARISON WITH EXPERIMENTS - 3D}
\end{figure}

The comparison between experimental observations and the numerical results is shown in Figs.(\ref{fig:RESULTS - SHOCK REFLECTION - COMPARISON WITH EXPERIMENTS - 3D}a) and (\ref{fig:RESULTS - SHOCK REFLECTION - COMPARISON WITH EXPERIMENTS - 3D}b), which show the Schlieren photograph fom Ref.\cite{bryson1961diffraction} of shock reflection from a sphere and the results at $t = 0.16$.
Although small differences likely due to viscous effects can be seen inside the wake behind the sphere, the $hp$-AMR strategy captures all the details of the flow structures observed in the experiments.

\subsection{Shock reflection from a convex-concave surface}\label{ssec:SHOCK REFLECTION CC}
\begin{figure}
\centering
\includegraphics[width = \textwidth]{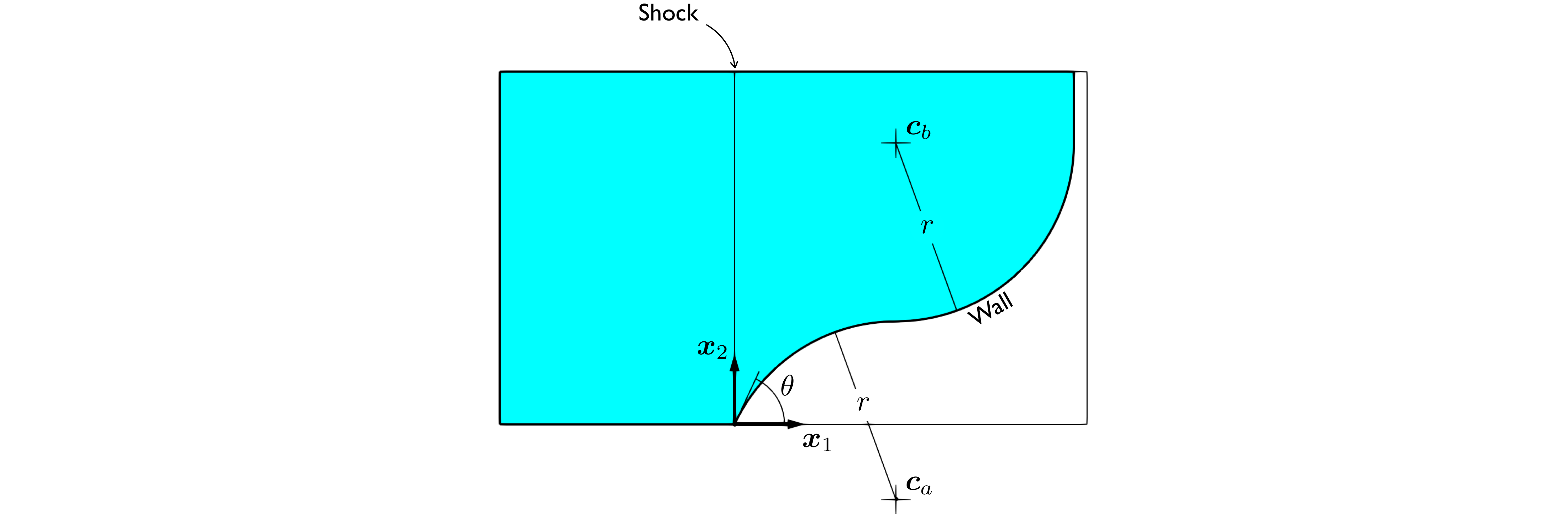}
\caption
{
    Illustration of the problem of a shock reflecting from a convex-concave wall.
}
\label{fig:RESULTS - SHOCK REFLECTION CC - GEOMETRY}
\end{figure}

The final example is another shock reflection problem consisting of a shock reflecting from a convex-concave wall as illustrated in Fig.(\ref{fig:RESULTS - SHOCK REFLECTION CC - GEOMETRY}).
Shocks reflecting from curved surfaces have been recently investigated experimentally and numerically \cite{ram2015high,soni2017shock,koronio2020similarity} because of their importance in shock-focusing applications.
Here, we model the experimental test of Ram et al.\cite{ram2015high}, who provided highly resolved measurements of the transition from regular to Mach reflection for the problem depicted in Fig.(\ref{fig:RESULTS - SHOCK REFLECTION CC - GEOMETRY}).

The problem is defined on the background rectangle $\mathscr{R} = [-80,122.5]\times[0,120]$.
The height of the numerical domain is larger than the height of the observation region in the experiments to avoid any reflection from the top boundary.
The convex-concave wall consists of two surfaces generated by two circular arcs and a flat surface.
With reference to Fig.(\ref{fig:RESULTS - SHOCK REFLECTION CC - GEOMETRY}), the radius of the cylinders is $r = 60.6$, while the location of their centers is given by $\bm{c}_{a} \equiv \{r\sin{\theta},-r\cos{\theta}\}$ and $\bm{c}_{b} \equiv \{r\sin{\theta},r(2-\cos{\theta})\}$ with $\theta = 65^\circ$.
Then, the geometry is defined by the level set function
\begin{equation}
\bm{\Phi}(\bm{x}) \equiv 
\left\{
\begin{array}{cc}
d_2(\bm{x})&\mathrm{if}~x_1 > c_{b1} \wedge x_2 < c_{b2}\\
\max\{d_1(\bm{x}),d_3(\bm{x})\}&\mathrm{otherwise}
\end{array}
\right.,
\end{equation}
where $d_1(\bm{x}) \equiv r-\sqrt{(x_k-c_{ak})(x_k-c_{ak})}$, $d_2(\bm{x}) \equiv \sqrt{(x_k-c_{bk})(x_k-c_{bk})}-r$ and $d_3(\bm{x}) \equiv x_1-c_{b1}-r$.

The shock is initially located at $x_1 = 0$.
Following the experiments \cite{ram2015high}, we consider three shock Mach numbers, $M = 1.19$, $M = 1.30$ and $M = 1.41$.
The initial conditions are given by Eq.(\ref{eq:shock reflection problem - ICs}) where $d_s = 0$.
Wall boundary conditions are prescribed on the convex-concave surface and at the top boundary of the background rectangle, whereas inflow boundary conditions are prescribed at $x_1 = -80$.
The final time of the simulation is chosen to make sure the shock does not reach the flat wall.

We use a four-level mesh with a refinement ratio of four between levels.
Level $\ell = 0$ uses a dG$_1$ scheme with mesh size $h = 3.75$;
level $\ell = 1$ uses a dG$_3$ scheme and is generated from level 0 using $\kappa_0^\rho = 0.005$ in (\ref{eq:density gradient refinement});
level $\ell = 2$ uses a dG$_3$ scheme and is generated from level 1 using $\kappa_1^\rho = 0.02$;
level $\ell = 3$ uses the FV scheme and is generated from level 2 using $\kappa^s = -4.00$ in (\ref{eq:shock sensor}).
Finally, cell merging is triggered by the volume fraction thresholds $\overline{\nu}_0 = \overline{\nu}_1 = \overline{\nu}_2 = \overline{\nu}_3 = 0.3$.

\begin{figure}
\centering
\includegraphics[width = \textwidth]{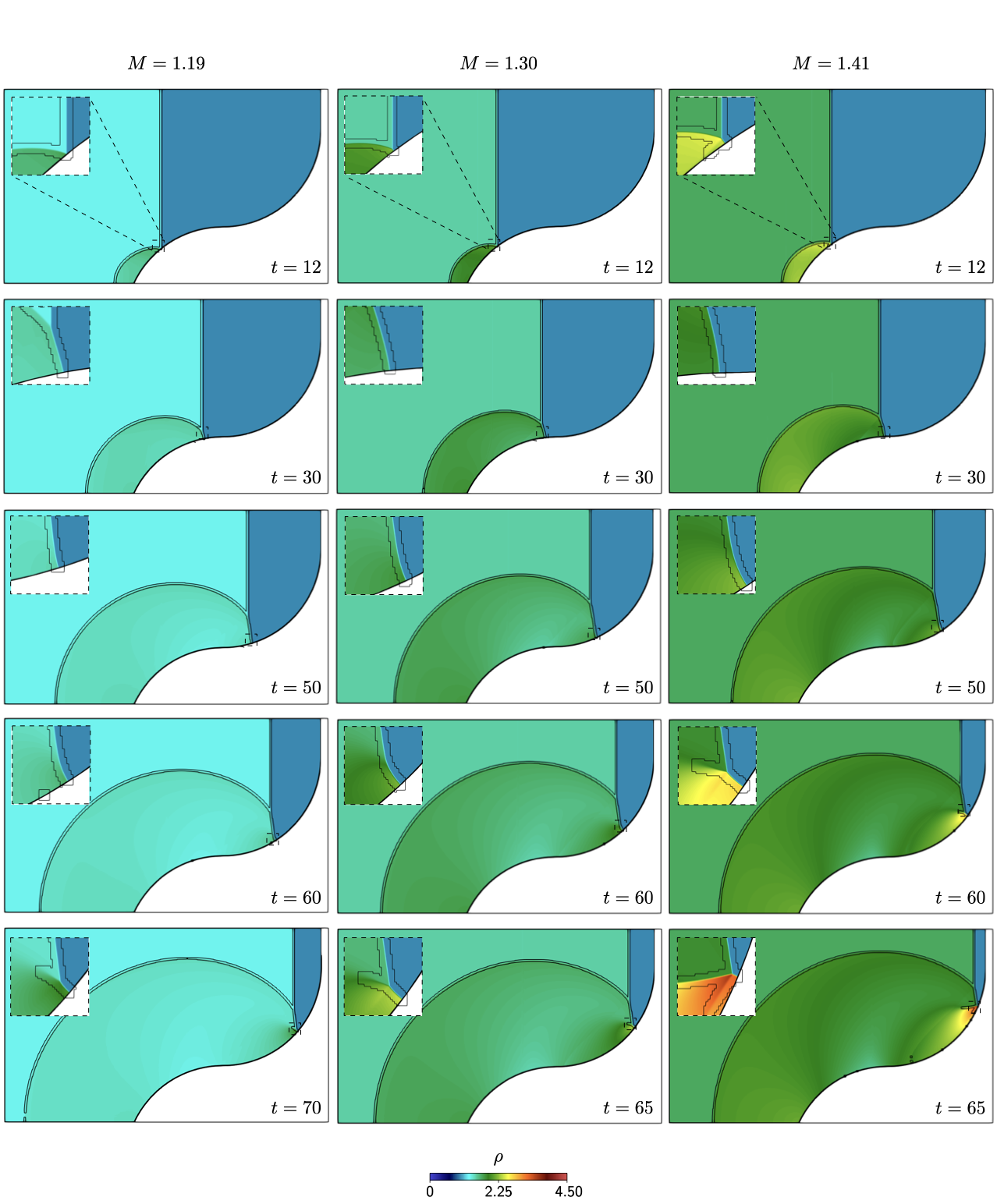}
\caption
{
    Density distribution at selected times for shock reflection from a convex-concave wall at different Mach numbers: (left column) $M = 1.19$, (center column) $M = 1.30$ and (right column) $M = 1.41$.
    In each image, the inset shows a close-up of the region where the shock is attached to the curved surface.
    The time is reported in the bottom right corner of each image.
}
\label{fig:RESULTS - SHOCK REFLECTION CC - DENSITY}
\end{figure}

\begin{figure}
\centering
\includegraphics[width = \textwidth]{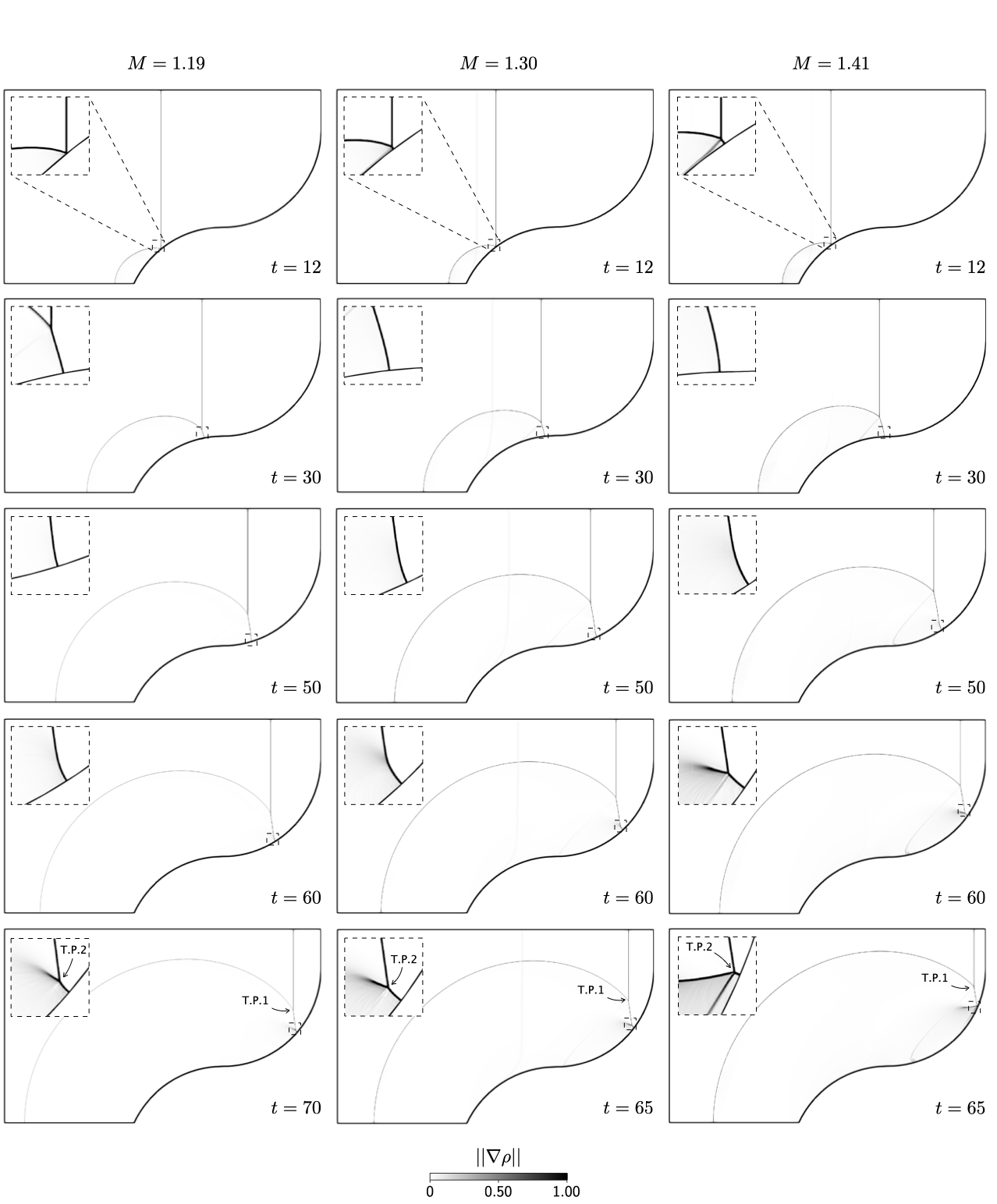}
\caption
{
    Density gradient distribution at selected times for shock reflection from a convex-concave wall at different Mach numbers: (left column) $M = 1.19$, (center column) $M = 1.30$ and (right column) $M = 1.41$.
    In each image, the inset shows a close-up of the region where the shock is attached to the curved surface.
    The time is reported in the bottom right corner of each image.
}
\label{fig:RESULTS - SHOCK REFLECTION CC - DENSITY GRADIENT}
\end{figure}

Figure (\ref{fig:RESULTS - SHOCK REFLECTION CC - DENSITY}) show the distribution of the density $\rho$ at selected times for the different Mach numbers.
In the figure, the inset of each image shows a close-up of the region where the shock is attached to the curved surface and allows us to distinguish the different shock structures.
For instance, for the $M = 1.19$ case reported in the left column of Fig.(\ref{fig:RESULTS - SHOCK REFLECTION CC - DENSITY}), the reflection structure transitions from regular reflection at $t = 12$ to Mach reflection at $t = 30$.
Similarly, the last row of Fig.(\ref{fig:RESULTS - SHOCK REFLECTION CC - DENSITY}) highlights the formation of the secondary triple point.
The images also show the outline of the FV level, which is able to track the shock in all its difference configurations.

The different shock structures can be seen more clearly in Fig.(\ref{fig:RESULTS - SHOCK REFLECTION CC - DENSITY GRADIENT}), which shows the density gradient magnitude $||\nabla\rho||$ at the same times as shown in Fig.(\ref{fig:RESULTS - SHOCK REFLECTION CC - DENSITY}), allowing us to see the effect of the Mach number on the flow patterns.
For instance, at $t = 12$, the flow with $M = 1.19$ is characterized by regular reflection, the flow with $M = 1.41$ is characterized by Mach reflection and the flow with $M = 1.30$ is transitioning between the two configurations.

\begin{figure}
\centering
\begin{subfigure}{0.32\textwidth}
\includegraphics[width=\textwidth]{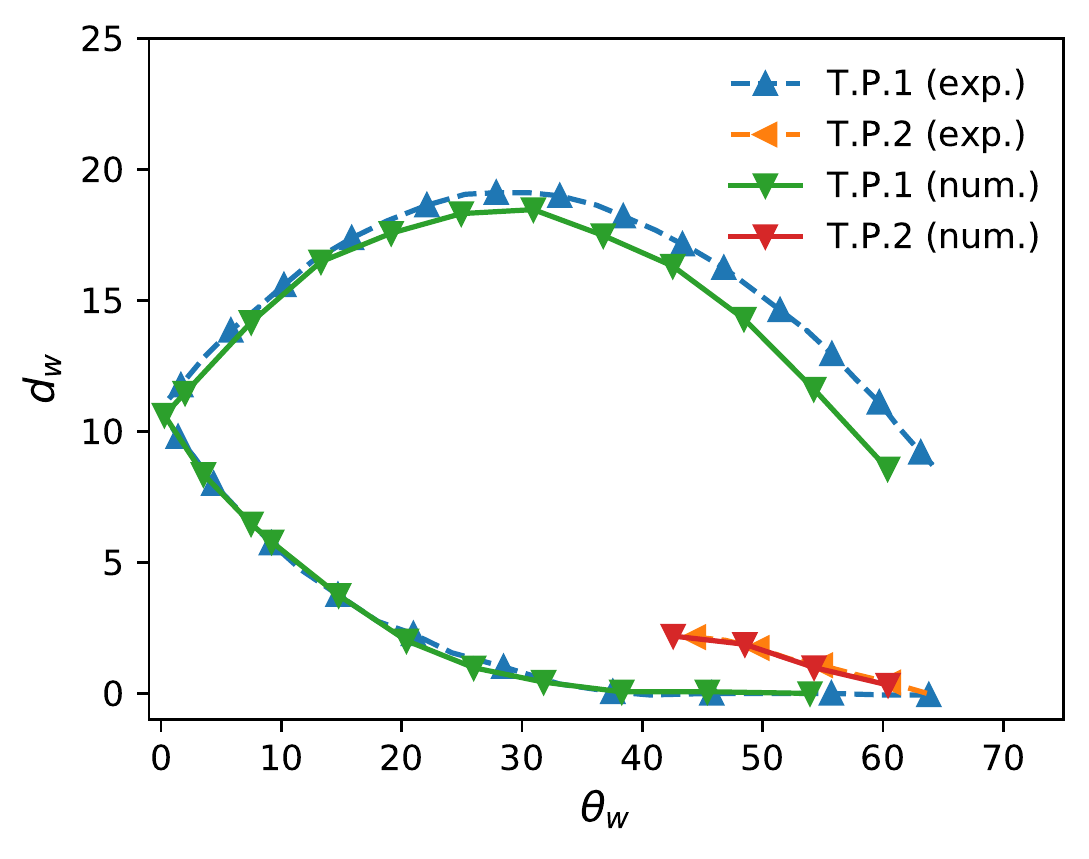}
\caption{}
\end{subfigure}
\
\begin{subfigure}{0.32\textwidth}
\includegraphics[width=\textwidth]{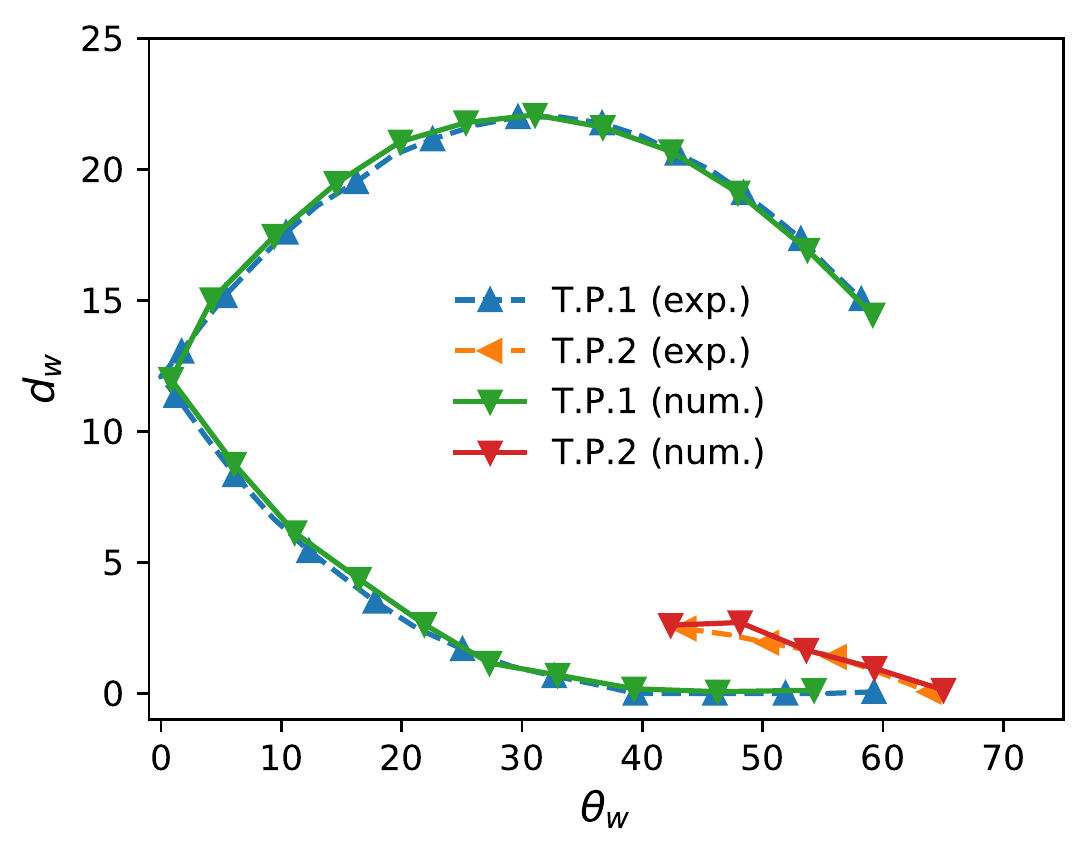}
\caption{}
\end{subfigure}
\
\begin{subfigure}{0.32\textwidth}
\includegraphics[width=\textwidth]{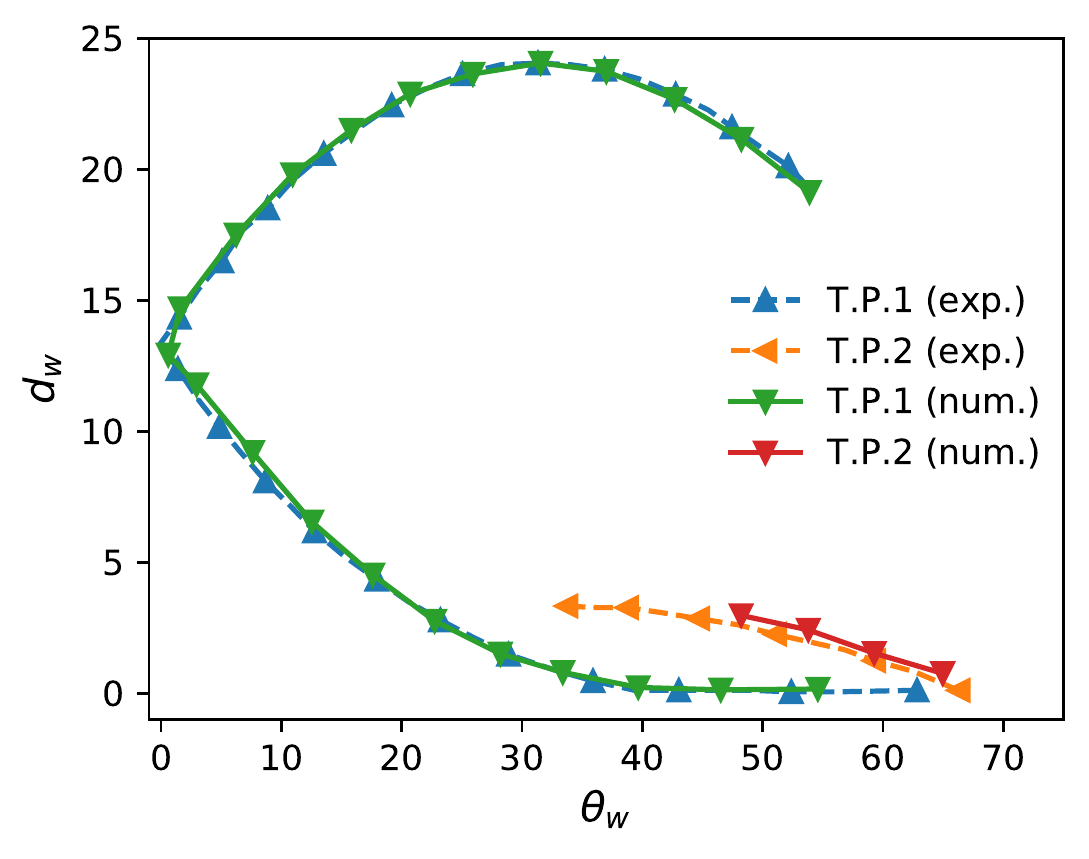}
\caption{}
\end{subfigure}
\caption
{
    Plot of the distance $d_w$ of the triple points from the wall versus the inclination $\theta_w$ of wall surface at the foot of the Mach stem:
    (a) $M = 1.19$,
    (b) $M = 1.30$, and
    (c) $M = 1.41$.
    In the legend of the figures, exp.~denotes the experimental measurements and num.~denotes the numerical results, while T.P.1 and T.P.2 refer to the results of the primary and secondary triple points, respectively, as indicated in the last row of Fig.(\ref{fig:RESULTS - SHOCK REFLECTION CC - DENSITY GRADIENT}).
    Experimental results are taken from Ref.\cite{ram2015high}.
}
\label{fig:RESULTS - SHOCK REFLECTION CC - COMPARISON WITH EXPERIMENTS}
\end{figure}

Finally, using a sequence of snapshots such as those presented in Fig.(\ref{fig:RESULTS - SHOCK REFLECTION CC - DENSITY GRADIENT}), one can quantitatively compare the numerical results with the experimental measurements.
Specifically, in Ref.\cite{ram2015high}, the Authors measured the distance $d_w$ between from both the primary triple point and the secondary triple points to the reflecting surface and plotted $d_w$ versus the angle of the reflecting surface at the foot of the Mach stem, $\theta_w$.
The comparison between the experiments and the numerics is shown in Figs.(\ref{fig:RESULTS - SHOCK REFLECTION CC - COMPARISON WITH EXPERIMENTS}a), (\ref{fig:RESULTS - SHOCK REFLECTION CC - COMPARISON WITH EXPERIMENTS}b) and (\ref{fig:RESULTS - SHOCK REFLECTION CC - COMPARISON WITH EXPERIMENTS}c), which corresponds to the cases $M = 1.19$, $M = 1.30$, $M = 1.41$, respectively.
As shown by the figures, the  embedded-boundary $hp$-AMR framework is able to reproduce the experimental measurements with remarkable accuracy.

\section{Conclusions and future work}\label{sec:CONCLUSIONS}
We have presented a methodology for solving the equations of inviscid gas dynamics on structured grids with an embedded boundary representation of complex geometry.
The novel features of the approach are the coupling of a discontinuous Galerkin method that provides high-order accuracy in regions of smooth flow with a Finite Volume method that provides robust shock-tracking capabilities of solution discontinuities through an $hp$ adaptive mesh refinement strategy and an implicit mesh approach that enables high-order accurate resolution of the geometry.
The methodology uses a compact-stencil AMR discretization strategy, is fully explicit and does not require the introduction of additional variables (e.g.~the auxiliary flux in artificial-viscosity methods).

The methodology has been tested on embedded-geometry problems with and without solution discontinuities.
The supersonic vortex problem of Sec.(\ref{ssec:SUPERSONIC VORTEX}) demonstrates the high-order accuracy of dG schemes for smooth flows, while the remaining examples illustrate the performance of the method on problems with discontinuities in several configurations.

The method offers several avenues of improvements and further research.
First, the current AMR strategy relies on the definition of two parameters, namely $\kappa_\ell^\rho$ and $\kappa^s$, entering Eqs.(\ref{eq:density gradient refinement}) and (\ref{eq:smoothness indicator}).
Here, $\kappa_\ell^\rho$ is chosen to ensure that regions of constant fields are resolved with low-order dG schemes while the remaining parts of the domain where the flow is smooth are resolved with high-order dG schemes.
In addition, $\kappa^s$ is chosen so that contact discontinuities transition as soon as possible from being in treated the FV scheme to being treated by dG, while making sure shocks are always contained within the FV level.
In the simulations, this criterion allowed high-order accuracy in the regions of smooth flow and robustly captured the travelling shocks but led to a problem-dependent choice of the threshold parameters.
Therefore, more advanced conditions for the dynamic evolution of the $hp$-AMR levels, including more sophisticated shock sensors \cite{krivodonova2004shock,klockner2011viscous,lv2016entropy,fernandez2018physics}, or a more robust treatment of the coupling between the FV and the dG schemes, using e.g.~\emph{a posteriori} limiters \cite{clain2011high,dumbser2014posteriori,zanotti2015space}, might be beneficial.

Second, Eq.(\ref{eq:ode}) is integrated in time using high-order TVD Runge-Kutta schemes, which are
very efficient for many practical applications but are limited to fourth order accuracy.
Methods such as the spectral deferred correction methods \cite{minion2004semi,almgren2013use,minion2018higher} or ADER schemes \cite{titarev2002ader,balsara2013efficient,dumbser2014posteriori,zanotti2015space} could be implemented within the present framework to obtain space-time orders of accuracy that are higher than the ones considered here.

The numerical test cases considered in Sec.(\ref{sec:RESULTS}) used geometries that could be described by relatively simple level set functions.
This does not represent a fundamental limitation of the proposed framework, which could be employed to model more complex geometries, e.g.~including sharp edges or corners, provided that the corresponding high-order quadrature rules are available.
In fact, it is worth pointing out that even the simple three-dimensional embedded geometries considered in Sec.(\ref{sec:RESULTS}) induce very complex configurations of cut and merged elements, which are robustly handled by present approach in both static and dynamic AMR configurations.
However, implicitly-defined meshes where cells are split by embedded boundaries that are smaller than the cell size would required doubly-value cells, which are currently not supported.

The methodology was implemented using \texttt{AMReX} \cite{zhang2019amrex}, an exascale-ready software framework for massively parallel, block-structured adaptive mesh refinement applications.
The computations were performed using classic MPI parallelization and future research will provide a comprehensive scalability analysis involving also the use of modern accelerators, such as general-purpose graphical processing units.

Future lines of research will also include the addition of viscous and reaction terms to the governing equations (\ref{eq:governing equations - vectors}) and the implementation of the corresponding dG and/or FV schemes. 
This will provide high-order numerical schemes for more complex phenomena, such as combustion, where high-order resolution of curved geometries would be desirable.
Finally, since the methodology already features dynamic adaptivity and dynamic generation of the implicit mesh, the extension to moving geometry would be relatively straightforward.
The moving geometry could either be a prescribed motion or it could be computed as part of the evolution of the system by evolving different dynamics on either side of the interface.

\section*{Acknowledgements}
We thank Dr.~Robert Saye for the support in the use of \texttt{Algoim} (\url{https://github.com/algoim/algoim}) and Dr.~Weiqun Zhang for the support in the use of \texttt{AMReX} (\url{https://amrex-codes.github.io/amrex/}).
The work was supported by the Exascale Computing Project (17-SC-20-SC), a collaborative effort of the U.S. Department of Energy Office of Science and the National Nuclear Security Administration, through
U.S.~Department of Energy, Office of Science,
Office of Advanced Scientific Computing Research,
under contract DE-AC02-05CH11231.
This research used resources of the National Energy Research Scientific Computing Center, a DOE Office of Science User Facility supported by the Office of Science of the U.S. Department of Energy under Contract No. DE-AC02-05CH11231; 
and resources of the Oak Ridge Leadership Computing Facility at the Oak Ridge National Laboratory, which is supported by the Office of Science of the U.S. Department of Energy under Contract No. DE-AC05-00OR22725.


\bibliography{bibliography}

\begin{thebibliography}{80}
\providecommand{\natexlab}[1]{#1}
\providecommand{\url}[1]{\texttt{#1}}
\expandafter\ifx\csname urlstyle\endcsname\relax
  \providecommand{\doi}[1]{doi: #1}\else
  \providecommand{\doi}{doi: \begingroup \urlstyle{rm}\Url}\fi

\bibitem[Aftosmis et~al.(2000)Aftosmis, Berger, and
  Adomavicius]{aftosmis2000parallel}
Aftosmis, M., Berger, M., and Adomavicius, G.
\newblock A parallel multilevel method for adaptively refined cartesian grids
  with embedded boundaries.
\newblock In \emph{38th Aerospace Sciences Meeting and Exhibit}, page 808,
  2000.

\bibitem[Peskin(2002)]{peskin2002immersed}
Peskin, C.~S.
\newblock The immersed boundary method.
\newblock \emph{Acta numerica}, 11:\penalty0 479--517, 2002.

\bibitem[Mokbel et~al.(2018)Mokbel, Abels, and Aland]{mokbel2018phase}
Mokbel, D., Abels, H., and Aland, S.
\newblock A phase-field model for fluid--structure interaction.
\newblock \emph{Journal of Computational Physics}, 372:\penalty0 823--840,
  2018.

\bibitem[Kemm et~al.(2020)Kemm, Gaburro, Thein, and Dumbser]{kemm2020simple}
Kemm, F., Gaburro, E., Thein, F., and Dumbser, M.
\newblock A simple diffuse interface approach for compressible flows around
  moving solids of arbitrary shape based on a reduced baer--nunziato model.
\newblock \emph{Computers \& Fluids}, 204:\penalty0 104536, 2020.

\bibitem[Mo{\"e}s et~al.(1999)Mo{\"e}s, Dolbow, and Belytschko]{moes1999finite}
Mo{\"e}s, N., Dolbow, J., and Belytschko, T.
\newblock A finite element method for crack growth without remeshing.
\newblock \emph{International journal for numerical methods in engineering},
  46\penalty0 (1):\penalty0 131--150, 1999.

\bibitem[Qin and Krivodonova(2013)]{qin2013discontinuous}
Qin, R. and Krivodonova, L.
\newblock A discontinuous galerkin method for solutions of the euler equations
  on cartesian grids with embedded geometries.
\newblock \emph{Journal of Computational Science}, 4\penalty0 (1-2):\penalty0
  24--35, 2013.

\bibitem[Alauzet et~al.(2016)Alauzet, Fabr{\`e}ges, Fern{\'a}ndez, and
  Landajuela]{alauzet2016nitsche}
Alauzet, F., Fabr{\`e}ges, B., Fern{\'a}ndez, M.~A., and Landajuela, M.
\newblock Nitsche-xfem for the coupling of an incompressible fluid with
  immersed thin-walled structures.
\newblock \emph{Computer Methods in Applied Mechanics and Engineering},
  301:\penalty0 300--335, 2016.

\bibitem[Saye(2017{\natexlab{a}})]{saye2017implicit}
Saye, R.
\newblock Implicit mesh discontinuous galerkin methods and interfacial gauge
  methods for high-order accurate interface dynamics, with applications to
  surface tension dynamics, rigid body fluid--structure interaction, and free
  surface flow: Part i.
\newblock \emph{Journal of Computational Physics}, 344:\penalty0 647--682,
  2017{\natexlab{a}}.

\bibitem[Saye(2017{\natexlab{b}})]{saye2017implicitII}
Saye, R.
\newblock Implicit mesh discontinuous galerkin methods and interfacial gauge
  methods for high-order accurate interface dynamics, with applications to
  surface tension dynamics, rigid body fluid--structure interaction, and free
  surface flow: Part ii.
\newblock \emph{Journal of Computational Physics}, 344:\penalty0 683--723,
  2017{\natexlab{b}}.

\bibitem[Milazzo et~al.(2018)Milazzo, Benedetti, and
  Gulizzi]{milazzo2018extended}
Milazzo, A., Benedetti, I., and Gulizzi, V.
\newblock An extended ritz formulation for buckling and post-buckling analysis
  of cracked multilayered plates.
\newblock \emph{Composite Structures}, 201:\penalty0 980--994, 2018.

\bibitem[Berger(2017)]{berger2017cut}
Berger, M.
\newblock Cut cells: Meshes and solvers.
\newblock In \emph{Handbook of Numerical Analysis}, volume~18, pages 1--22.
  Elsevier, 2017.

\bibitem[LeVeque et~al.(2002)]{leveque2002finite}
LeVeque, R.~J. et~al.
\newblock \emph{Finite volume methods for hyperbolic problems}, volume~31.
\newblock Cambridge university press, 2002.

\bibitem[Versteeg and Malalasekera(2007)]{versteeg2007introduction}
Versteeg, H.~K. and Malalasekera, W.
\newblock \emph{An introduction to computational fluid dynamics: the finite
  volume method}.
\newblock Pearson education, 2007.

\bibitem[Clarke et~al.(1986)Clarke, Salas, and Hassan]{clarke1986euler}
Clarke, D.~K., Salas, M., and Hassan, H.
\newblock Euler calculations for multielement airfoils using cartesian grids.
\newblock \emph{AIAA journal}, 24\penalty0 (3):\penalty0 353--358, 1986.

\bibitem[Gaffney and Hassan(1987)]{gaffney1987euler}
Gaffney, JR, R. and Hassan, H.
\newblock Euler calculations for wings using cartesian grids.
\newblock In \emph{25th AIAA Aerospace Sciences Meeting}, page 356, 1987.

\bibitem[Hartmann et~al.(2008)Hartmann, Meinke, and
  Schr{\"o}der]{hartmann2008adaptive}
Hartmann, D., Meinke, M., and Schr{\"o}der, W.
\newblock An adaptive multilevel multigrid formulation for cartesian
  hierarchical grid methods.
\newblock \emph{Computers \& Fluids}, 37\penalty0 (9):\penalty0 1103--1125,
  2008.

\bibitem[Ji et~al.(2010)Ji, Lien, and Yee]{ji2010numerical}
Ji, H., Lien, F.-S., and Yee, E.
\newblock Numerical simulation of detonation using an adaptive cartesian
  cut-cell method combined with a cell-merging technique.
\newblock \emph{Computers \& fluids}, 39\penalty0 (6):\penalty0 1041--1057,
  2010.

\bibitem[Hartmann et~al.(2011)Hartmann, Meinke, and
  Schr{\"o}der]{hartmann2011strictly}
Hartmann, D., Meinke, M., and Schr{\"o}der, W.
\newblock A strictly conservative cartesian cut-cell method for compressible
  viscous flows on adaptive grids.
\newblock \emph{Computer Methods in Applied Mechanics and Engineering},
  200\penalty0 (9-12):\penalty0 1038--1052, 2011.

\bibitem[Cecere and Giacomazzi(2014)]{cecere2014immersed}
Cecere, D. and Giacomazzi, E.
\newblock An immersed volume method for large eddy simulation of compressible
  flows using a staggered-grid approach.
\newblock \emph{Computer Methods in Applied Mechanics and Engineering},
  280:\penalty0 1--27, 2014.

\bibitem[Muralidharan and Menon(2016)]{muralidharan2016high}
Muralidharan, B. and Menon, S.
\newblock A high-order adaptive cartesian cut-cell method for simulation of
  compressible viscous flow over immersed bodies.
\newblock \emph{Journal of Computational Physics}, 321:\penalty0 342--368,
  2016.

\bibitem[Pember et~al.(1995)Pember, Bell, Colella, Curtchfield, and
  Welcome]{pember1995adaptive}
Pember, R.~B., Bell, J.~B., Colella, P., Curtchfield, W.~Y., and Welcome, M.~L.
\newblock An adaptive cartesian grid method for unsteady compressible flow in
  irregular regions.
\newblock \emph{Journal of computational Physics}, 120\penalty0 (2):\penalty0
  278--304, 1995.

\bibitem[Almgren et~al.(1997)Almgren, Bell, Colella, and
  Marthaler]{almgren1997cartesian}
Almgren, A.~S., Bell, J.~B., Colella, P., and Marthaler, T.
\newblock A cartesian grid projection method for the incompressible euler
  equations in complex geometries.
\newblock \emph{SIAM Journal on Scientific Computing}, 18\penalty0
  (5):\penalty0 1289--1309, 1997.

\bibitem[Colella et~al.(2006)Colella, Graves, Keen, and
  Modiano]{colella2006cartesian}
Colella, P., Graves, D.~T., Keen, B.~J., and Modiano, D.
\newblock A cartesian grid embedded boundary method for hyperbolic conservation
  laws.
\newblock \emph{Journal of Computational Physics}, 211\penalty0 (1):\penalty0
  347--366, 2006.

\bibitem[Graves et~al.(2013)Graves, Colella, Modiano, Johnson, Sjogreen, and
  Gao]{graves2013cartesian}
Graves, D., Colella, P., Modiano, D., Johnson, J., Sjogreen, B., and Gao, X.
\newblock A cartesian grid embedded boundary method for the compressible
  navier--stokes equations.
\newblock \emph{Communications in Applied Mathematics and Computational
  Science}, 8\penalty0 (1):\penalty0 99--122, 2013.

\bibitem[Helzel et~al.(2005)Helzel, Berger, and LeVeque]{helzel2005high}
Helzel, C., Berger, M.~J., and LeVeque, R.~J.
\newblock A high-resolution rotated grid method for conservation laws with
  embedded geometries.
\newblock \emph{SIAM Journal on Scientific Computing}, 26\penalty0
  (3):\penalty0 785--809, 2005.

\bibitem[Berger and Helzel(2012)]{berger2012simplified}
Berger, M. and Helzel, C.
\newblock A simplified h-box method for embedded boundary grids.
\newblock \emph{SIAM Journal on Scientific Computing}, 34\penalty0
  (2):\penalty0 A861--A888, 2012.

\bibitem[Gokhale et~al.(2018)Gokhale, Nikiforakis, and
  Klein]{gokhale2018dimensionally}
Gokhale, N., Nikiforakis, N., and Klein, R.
\newblock A dimensionally split cartesian cut cell method for hyperbolic
  conservation laws.
\newblock \emph{Journal of Computational Physics}, 364:\penalty0 186--208,
  2018.

\bibitem[Berger and Giuliani(2021)]{berger2021state}
Berger, M. and Giuliani, A.
\newblock A state redistribution algorithm for finite volume schemes on cut
  cell meshes.
\newblock \emph{Journal of Computational Physics}, 428:\penalty0 109820, 2021.

\bibitem[Cockburn et~al.(1990)Cockburn, Hou, and Shu]{cockburn1990runge}
Cockburn, B., Hou, S., and Shu, C.-W.
\newblock The runge-kutta local projection discontinuous galerkin finite
  element method for conservation laws. iv. the multidimensional case.
\newblock \emph{Mathematics of Computation}, 54\penalty0 (190):\penalty0
  545--581, 1990.

\bibitem[Cockburn and Shu(1998)]{cockburn1998runge}
Cockburn, B. and Shu, C.-W.
\newblock The runge--kutta discontinuous galerkin method for conservation laws
  v: multidimensional systems.
\newblock \emph{Journal of Computational Physics}, 141\penalty0 (2):\penalty0
  199--224, 1998.

\bibitem[Arnold et~al.(2002)Arnold, Brezzi, Cockburn, and
  Marini]{arnold2002unified}
Arnold, D.~N., Brezzi, F., Cockburn, B., and Marini, L.~D.
\newblock Unified analysis of discontinuous galerkin methods for elliptic
  problems.
\newblock \emph{SIAM journal on numerical analysis}, 39\penalty0 (5):\penalty0
  1749--1779, 2002.

\bibitem[Cockburn(2018)]{cockburn2018discontinuous}
Cockburn, B.
\newblock Discontinuous galerkin methods for computational fluid dynamics.
\newblock \emph{Encyclopedia of Computational Mechanics Second Edition}, pages
  1--63, 2018.

\bibitem[Cangiani et~al.(2017)Cangiani, Dong, and Georgoulis]{cangiani2017hp}
Cangiani, A., Dong, Z., and Georgoulis, E.~H.
\newblock hp-version space-time discontinuous galerkin methods for parabolic
  problems on prismatic meshes.
\newblock \emph{SIAM Journal on Scientific Computing}, 39\penalty0
  (4):\penalty0 A1251--A1279, 2017.

\bibitem[Antonietti and Pennesi(2019)]{antonietti2019v}
Antonietti, P.~F. and Pennesi, G.
\newblock V-cycle multigrid algorithms for discontinuous galerkin methods on
  non-nested polytopic meshes.
\newblock \emph{Journal of Scientific Computing}, 78\penalty0 (1):\penalty0
  625--652, 2019.

\bibitem[Bastian and Engwer(2009)]{bastian2009unfitted}
Bastian, P. and Engwer, C.
\newblock An unfitted finite element method using discontinuous galerkin.
\newblock \emph{International journal for numerical methods in engineering},
  79\penalty0 (12):\penalty0 1557--1576, 2009.

\bibitem[Johansson and Larson(2013)]{johansson2013high}
Johansson, A. and Larson, M.~G.
\newblock A high order discontinuous galerkin nitsche method for elliptic
  problems with fictitious boundary.
\newblock \emph{Numerische Mathematik}, 123\penalty0 (4):\penalty0 607--628,
  2013.

\bibitem[Brandstetter and Govindjee(2015)]{brandstetter2015high}
Brandstetter, G. and Govindjee, S.
\newblock A high-order immersed boundary discontinuous-galerkin method for
  poisson's equation with discontinuous coefficients and singular sources.
\newblock \emph{International Journal for Numerical Methods in Engineering},
  101\penalty0 (11):\penalty0 847--869, 2015.

\bibitem[G{\"u}rkan and Massing(2019)]{gurkan2019stabilized}
G{\"u}rkan, C. and Massing, A.
\newblock A stabilized cut discontinuous galerkin framework for elliptic
  boundary value and interface problems.
\newblock \emph{Computer Methods in Applied Mechanics and Engineering},
  348:\penalty0 466--499, 2019.

\bibitem[Lew and Buscaglia(2008)]{lew2008discontinuous}
Lew, A.~J. and Buscaglia, G.~C.
\newblock A discontinuous-galerkin-based immersed boundary method.
\newblock \emph{International Journal for Numerical Methods in Engineering},
  76\penalty0 (4):\penalty0 427--454, 2008.

\bibitem[Heimann et~al.(2013)Heimann, Engwer, Ippisch, and
  Bastian]{heimann2013unfitted}
Heimann, F., Engwer, C., Ippisch, O., and Bastian, P.
\newblock An unfitted interior penalty discontinuous galerkin method for
  incompressible navier--stokes two-phase flow.
\newblock \emph{International Journal for Numerical Methods in Fluids},
  71\penalty0 (3):\penalty0 269--293, 2013.

\bibitem[Kummer(2017)]{kummer2017extended}
Kummer, F.
\newblock Extended discontinuous galerkin methods for two-phase flows: the
  spatial discretization.
\newblock \emph{International Journal for Numerical Methods in Engineering},
  109\penalty0 (2):\penalty0 259--289, 2017.

\bibitem[Krause and Kummer(2017)]{krause2017incompressible}
Krause, D. and Kummer, F.
\newblock An incompressible immersed boundary solver for moving body flows
  using a cut cell discontinuous galerkin method.
\newblock \emph{Computers \& Fluids}, 153:\penalty0 118--129, 2017.

\bibitem[Saye(2020)]{saye2020fast}
Saye, R.
\newblock Fast multigrid solution of high-order accurate multiphase stokes
  problems.
\newblock \emph{Communications in Applied Mathematics and Computational
  Science}, 15\penalty0 (2):\penalty0 147--196, 2020.

\bibitem[Gulizzi et~al.(2020{\natexlab{a}})Gulizzi, Benedetti, and
  Milazzo]{gulizzi2020implicit}
Gulizzi, V., Benedetti, I., and Milazzo, A.
\newblock An implicit mesh discontinuous galerkin formulation for higher-order
  plate theories.
\newblock \emph{Mechanics of Advanced Materials and Structures}, 27\penalty0
  (17):\penalty0 1494--1508, 2020{\natexlab{a}}.

\bibitem[Gulizzi et~al.(2020{\natexlab{b}})Gulizzi, Benedetti, and
  Milazzo]{gulizzi2020high}
Gulizzi, V., Benedetti, I., and Milazzo, A.
\newblock A high-resolution layer-wise discontinuous galerkin formulation for
  multilayered composite plates.
\newblock \emph{Composite Structures}, 242:\penalty0 112137,
  2020{\natexlab{b}}.

\bibitem[M{\"u}ller et~al.(2017)M{\"u}ller, Kr{\"a}mer-Eis, Kummer, and
  Oberlack]{muller2017high}
M{\"u}ller, B., Kr{\"a}mer-Eis, S., Kummer, F., and Oberlack, M.
\newblock A high-order discontinuous galerkin method for compressible flows
  with immersed boundaries.
\newblock \emph{International Journal for Numerical Methods in Engineering},
  110\penalty0 (1):\penalty0 3--30, 2017.

\bibitem[Geisenhofer et~al.(2019)Geisenhofer, Kummer, and
  M{\"u}ller]{geisenhofer2019discontinuous}
Geisenhofer, M., Kummer, F., and M{\"u}ller, B.
\newblock A discontinuous galerkin immersed boundary solver for compressible
  flows: Adaptive local time stepping for artificial viscosity--based
  shock-capturing on cut cells.
\newblock \emph{International Journal for Numerical Methods in Fluids},
  91\penalty0 (9):\penalty0 448--472, 2019.

\bibitem[Fidkowski and Darmofal(2007)]{fidkowski2007triangular}
Fidkowski, K.~J. and Darmofal, D.~L.
\newblock A triangular cut-cell adaptive method for high-order discretizations
  of the compressible navier--stokes equations.
\newblock \emph{Journal of Computational Physics}, 225\penalty0 (2):\penalty0
  1653--1672, 2007.

\bibitem[Xiao et~al.(2019)Xiao, Febrianto, Zhang, and Cirak]{xiao2019immersed}
Xiao, H., Febrianto, E., Zhang, Q., and Cirak, F.
\newblock An immersed discontinuous galerkin method for compressible
  navier-stokes equations on unstructured meshes.
\newblock \emph{International Journal for Numerical Methods in Fluids},
  91\penalty0 (10):\penalty0 487--508, 2019.

\bibitem[Persson and Peraire(2006)]{persson2006sub}
Persson, P.-O. and Peraire, J.
\newblock Sub-cell shock capturing for discontinuous galerkin methods.
\newblock In \emph{44th AIAA Aerospace Sciences Meeting and Exhibit}, page 112,
  2006.

\bibitem[Chaudhuri et~al.(2017)Chaudhuri, Jacobs, Don, Abbassi, and
  Mashayek]{chaudhuri2017explicit}
Chaudhuri, A., Jacobs, G.~B., Don, W.-S., Abbassi, H., and Mashayek, F.
\newblock Explicit discontinuous spectral element method with entropy
  generation based artificial viscosity for shocked viscous flows.
\newblock \emph{Journal of Computational Physics}, 332:\penalty0 99--117, 2017.

\bibitem[Krivodonova(2007)]{krivodonova2007limiters}
Krivodonova, L.
\newblock Limiters for high-order discontinuous galerkin methods.
\newblock \emph{Journal of Computational Physics}, 226\penalty0 (1):\penalty0
  879--896, 2007.

\bibitem[Zahr et~al.(2020)Zahr, Shi, and Persson]{zahr2020implicit}
Zahr, M.~J., Shi, A., and Persson, P.-O.
\newblock Implicit shock tracking using an optimization-based high-order
  discontinuous galerkin method.
\newblock \emph{Journal of Computational Physics}, 410:\penalty0 109385, 2020.

\bibitem[Saye(2015)]{saye2015high}
Saye, R.
\newblock High-order quadrature methods for implicitly defined surfaces and
  volumes in hyperrectangles.
\newblock \emph{SIAM Journal on Scientific Computing}, 37\penalty0
  (2):\penalty0 A993--A1019, 2015.

\bibitem[Colella and Glaz(1985)]{colella1985efficient}
Colella, P. and Glaz, H.~M.
\newblock Efficient solution algorithms for the riemann problem for real gases.
\newblock \emph{Journal of Computational Physics}, 59\penalty0 (2):\penalty0
  264--289, 1985.

\bibitem[Toro(2013)]{toro2013riemann}
Toro, E.~F.
\newblock \emph{Riemann solvers and numerical methods for fluid dynamics: a
  practical introduction}.
\newblock Springer Science \& Business Media, 2013.

\bibitem[Van~Leer(1979)]{van1979towards}
Van~Leer, B.
\newblock Towards the ultimate conservative difference scheme. v. a
  second-order sequel to godunov's method.
\newblock \emph{Journal of computational Physics}, 32\penalty0 (1):\penalty0
  101--136, 1979.

\bibitem[Berger et~al.(2005)Berger, Aftosmis, and Muman]{berger2005analysis}
Berger, M., Aftosmis, M., and Muman, S.
\newblock Analysis of slope limiters on irregular grids.
\newblock In \emph{43rd AIAA Aerospace Sciences Meeting and Exhibit}, page 490,
  2005.

\bibitem[Barth and Jespersen(1989)]{barth1989design}
Barth, T. and Jespersen, D.
\newblock The design and application of upwind schemes on unstructured meshes.
\newblock In \emph{27th Aerospace sciences meeting}, page 366, 1989.

\bibitem[Berger and Colella(1989)]{berger1989local}
Berger, M.~J. and Colella, P.
\newblock Local adaptive mesh refinement for shock hydrodynamics.
\newblock \emph{Journal of computational Physics}, 82\penalty0 (1):\penalty0
  64--84, 1989.

\bibitem[Bell et~al.(1994)Bell, Berger, Saltzman, and Welcome]{bell1994three}
Bell, J., Berger, M., Saltzman, J., and Welcome, M.
\newblock Three-dimensional adaptive mesh refinement for hyperbolic
  conservation laws.
\newblock \emph{SIAM Journal on Scientific Computing}, 15\penalty0
  (1):\penalty0 127--138, 1994.

\bibitem[Zhang et~al.(2019)Zhang, Almgren, Beckner, Bell, Blaschke, Chan, Day,
  Friesen, Gott, Graves, et~al.]{zhang2019amrex}
Zhang, W., Almgren, A., Beckner, V., Bell, J., Blaschke, J., Chan, C., Day, M.,
  Friesen, B., Gott, K., Graves, D., et~al.
\newblock Amrex: a framework for block-structured adaptive mesh refinement.
\newblock \emph{Journal of Open Source Software}, 4\penalty0 (37):\penalty0
  1370--1370, 2019.

\bibitem[Fortunato et~al.(2019)Fortunato, Rycroft, and
  Saye]{fortunato2019efficient}
Fortunato, D., Rycroft, C.~H., and Saye, R.
\newblock Efficient operator-coarsening multigrid schemes for local
  discontinuous galerkin methods.
\newblock \emph{SIAM Journal on Scientific Computing}, 41\penalty0
  (6):\penalty0 A3913--A3937, 2019.

\bibitem[Krivodonova and Berger(2006)]{krivodonova2006high}
Krivodonova, L. and Berger, M.
\newblock High-order accurate implementation of solid wall boundary conditions
  in curved geometries.
\newblock \emph{Journal of computational physics}, 211\penalty0 (2):\penalty0
  492--512, 2006.

\bibitem[Bryson and Gross(1961)]{bryson1961diffraction}
Bryson, A. and Gross, R.
\newblock Diffraction of strong shocks by cones, cylinders, and spheres.
\newblock \emph{Journal of Fluid Mechanics}, 10\penalty0 (1):\penalty0 1--16,
  1961.

\bibitem[Ram et~al.(2015)Ram, Geva, and Sadot]{ram2015high}
Ram, O., Geva, M., and Sadot, O.
\newblock High spatial and temporal resolution study of shock wave reflection
  over a coupled convex--concave cylindrical surface.
\newblock \emph{Journal of Fluid Mechanics}, 768:\penalty0 219--239, 2015.

\bibitem[Soni et~al.(2017)Soni, Hadjadj, Chaudhuri, and Ben-Dor]{soni2017shock}
Soni, V., Hadjadj, A., Chaudhuri, A., and Ben-Dor, G.
\newblock Shock-wave reflections over double-concave cylindrical reflectors.
\newblock \emph{Journal of Fluid Mechanics}, 813:\penalty0 70, 2017.

\bibitem[Koronio et~al.(2020)Koronio, Ben-Dor, Sadot, and
  Geva]{koronio2020similarity}
Koronio, E., Ben-Dor, G., Sadot, O., and Geva, M.
\newblock Similarity in mach stem evolution and termination in unsteady
  shock-wave reflection.
\newblock \emph{Journal of Fluid Mechanics}, 902, 2020.

\bibitem[Krivodonova et~al.(2004)Krivodonova, Xin, Remacle, Chevaugeon, and
  Flaherty]{krivodonova2004shock}
Krivodonova, L., Xin, J., Remacle, J.-F., Chevaugeon, N., and Flaherty, J.~E.
\newblock Shock detection and limiting with discontinuous galerkin methods for
  hyperbolic conservation laws.
\newblock \emph{Applied Numerical Mathematics}, 48\penalty0 (3-4):\penalty0
  323--338, 2004.

\bibitem[Kl{\"o}ckner et~al.(2011)Kl{\"o}ckner, Warburton, and
  Hesthaven]{klockner2011viscous}
Kl{\"o}ckner, A., Warburton, T., and Hesthaven, J.~S.
\newblock Viscous shock capturing in a time-explicit discontinuous galerkin
  method.
\newblock \emph{Mathematical Modelling of Natural Phenomena}, 6\penalty0
  (3):\penalty0 57--83, 2011.

\bibitem[Lv et~al.(2016)Lv, See, and Ihme]{lv2016entropy}
Lv, Y., See, Y.~C., and Ihme, M.
\newblock An entropy-residual shock detector for solving conservation laws
  using high-order discontinuous galerkin methods.
\newblock \emph{Journal of Computational Physics}, 322:\penalty0 448--472,
  2016.

\bibitem[Fernandez et~al.(2018)Fernandez, Nguyen, and
  Peraire]{fernandez2018physics}
Fernandez, P., Nguyen, C., and Peraire, J.
\newblock A physics-based shock capturing method for unsteady laminar and
  turbulent flows.
\newblock In \emph{2018 AIAA Aerospace Sciences Meeting}, page 0062, 2018.

\bibitem[Clain et~al.(2011)Clain, Diot, and Loub{\`e}re]{clain2011high}
Clain, S., Diot, S., and Loub{\`e}re, R.
\newblock A high-order finite volume method for systems of conservation
  laws—multi-dimensional optimal order detection (mood).
\newblock \emph{Journal of computational Physics}, 230\penalty0 (10):\penalty0
  4028--4050, 2011.

\bibitem[Dumbser et~al.(2014)Dumbser, Zanotti, Loub{\`e}re, and
  Diot]{dumbser2014posteriori}
Dumbser, M., Zanotti, O., Loub{\`e}re, R., and Diot, S.
\newblock A posteriori subcell limiting of the discontinuous galerkin finite
  element method for hyperbolic conservation laws.
\newblock \emph{Journal of Computational Physics}, 278:\penalty0 47--75, 2014.

\bibitem[Zanotti et~al.(2015)Zanotti, Fambri, Dumbser, and
  Hidalgo]{zanotti2015space}
Zanotti, O., Fambri, F., Dumbser, M., and Hidalgo, A.
\newblock Space--time adaptive ader discontinuous galerkin finite element
  schemes with a posteriori sub-cell finite volume limiting.
\newblock \emph{Computers \& Fluids}, 118:\penalty0 204--224, 2015.

\bibitem[Minion(2004)]{minion2004semi}
Minion, M.~L.
\newblock Semi-implicit projection methods for incompressible flow based on
  spectral deferred corrections.
\newblock \emph{Applied numerical mathematics}, 48\penalty0 (3-4):\penalty0
  369--387, 2004.

\bibitem[Almgren et~al.(2013)Almgren, Aspden, Bell, and Minion]{almgren2013use}
Almgren, A.~S., Aspden, A., Bell, J.~B., and Minion, M.~L.
\newblock On the use of higher-order projection methods for incompressible
  turbulent flow.
\newblock \emph{SIAM Journal on Scientific Computing}, 35\penalty0
  (1):\penalty0 B25--B42, 2013.

\bibitem[Minion and Saye(2018)]{minion2018higher}
Minion, M.~L. and Saye, R.
\newblock Higher-order temporal integration for the incompressible
  navier--stokes equations in bounded domains.
\newblock \emph{Journal of Computational Physics}, 375:\penalty0 797--822,
  2018.

\bibitem[Titarev and Toro(2002)]{titarev2002ader}
Titarev, V.~A. and Toro, E.~F.
\newblock Ader: Arbitrary high order godunov approach.
\newblock \emph{Journal of Scientific Computing}, 17\penalty0 (1):\penalty0
  609--618, 2002.

\bibitem[Balsara et~al.(2013)Balsara, Meyer, Dumbser, Du, and
  Xu]{balsara2013efficient}
Balsara, D.~S., Meyer, C., Dumbser, M., Du, H., and Xu, Z.
\newblock Efficient implementation of ader schemes for euler and
  magnetohydrodynamical flows on structured meshes--speed comparisons with
  runge--kutta methods.
\newblock \emph{Journal of Computational Physics}, 235:\penalty0 934--969,
  2013.

\end{thebibliography}

\end{document}